\documentclass{amsart}
\usepackage{amsmath}
\usepackage{amsthm}
\usepackage{amssymb}
\usepackage{xypic}
\newcommand{\setmin}{\setminus}
\newcommand{\alef}{{\mathbf\aleph}}
\newcommand{\omom}{{{}^\omega\omega}}
\newcommand{\twoom}{{{}^\omega 2}}
\newcommand{\twolom}{{{}^{<\omega} 2}}

\newcommand{\tles }{\preceq}
\newcommand{\borel}{{\mathsf   {BOREL}}}
\newcommand{\comp}{\circ}
\newcommand{\nor}[1]{|\!|#1|\!|}
\newcommand{\K}{{\mathsf K}}
\newcommand{\Mon}{{\mathsf M}}
\newcommand{\FF}{{\mathsf F}}
\newcommand{\G}{{\mathsf {G}}}

\newcommand{\gb}{{\mathfrak b}}
\newcommand{\gd}{{\mathfrak d}}
\newcommand{\m}{{\mathfrak {ma}}}
\newcommand{\w}{{\mathfrak {am}}}

\newcommand{\ZFCa}{{{\mathsf {ZFC}}}}
\newcommand{\CH}{{\mathsf {CH}}}
\newcommand{\bv}[1]{[\![ #1 ]\!]}
\newcommand{\reals}{{\mathbb R}}
\newcommand{\rationals}{{\mathbb Q}}
\newcommand{\clopen}{{\mathbb C}}
\newcommand{\rest}{{\mathord{\restriction}}}
\newcommand{\ADD}{{\mathsf   {ADD}}}
\newcommand{\COV}{{\mathsf   {COV}}}
\newcommand{\UNIF}{{\mathsf   {NON}}}
\newcommand{\COF}{{\mathsf   {COF}}}
\newcommand{\B}{{\mathsf   {B}}}
\newcommand{\D}{{\mathsf   {D}}}
\newcommand{\gS}{{{\mathbf {S}}}_{g, g^\star}}
\newcommand{\dom}{{{\mathsf {dom}}}}
\newcommand{\stem}{{{\mathsf   {stem}}}}
\newcommand{\supp}{{{\mathsf    {dom}}}}
\newcommand{\suppp}{{{\mathsf    {supp}}}}

\newcommand{\suc}{{{\mathsf    {succ}}}}

\newcommand{\forces}{\Vdash}
\newcommand{\F}{{\mathcal F}}
\newcommand{\V}{{\mathbf V}}
\newcommand{\eL}{{\mathbf V}_0}
\newcommand{\HH}{{\mathbf H}}
\newcommand{\<}{\langle}
\renewcommand{\>}{\rangle}
\newcommand{\thinks}{\models}
\newcommand{\lft}[2]{\mathopen\ifcase#1{}\oo\or
                        \big#2\or\Big#2\else\oo\fi} 
\newcommand{\rgt}[2]{\mathclose\ifcase#1{}\oo\or
                        \big#2\or\Big#2\else\oo\fi} 
\newcommand{\Power}{{\mathbf P}}

\def\onesetup#1 #2 #3 #4 #5 {
\gdef\pone{#1}
\gdef\ptwo{#2}
\gdef\pthree{#3}
\gdef\pfour{#4}
\gdef\pfive{#5}}
\def\twosetup#1 #2 #3 #4 #5 {%
\gdef\psix{#1}
\gdef\pseven{#2}
\gdef\peight{#3}
\gdef\pnine{#4}
\gdef\pten{#5}}

\newcommand{\add}{{\mathsf   {add}}}
\newcommand{\cov}{{\mathsf  {cov}}}
\newcommand{\unif}{{\mathsf  {non}}}
\newcommand{\cof}{{\mathsf   {cof}}}
\newcommand{\cf}{{{\mathsf   {cf}}}}
\newcommand{\N}{{\mathcal N}}
\newcommand{\M}{{\mathcal M}}

\theoremstyle{plain}
\newtheorem{theorem}{Theorem}[section]
\theoremstyle{plain}
\newtheorem{lemma}[theorem]{Lemma}

\newtheorem{definition}[theorem]{Definition}
\newtheorem{example}[theorem]{Example}
\begin{document}
\title{Invariants of Measure and Category}
\author{Tomek Bartoszynski}
\address{Department of Mathematics and Computer Science\\
Boise State University\\
Boise, Idaho 83725 U.S.A.}
\thanks{Author  partially supported by 
NSF grant DMS 95-05375 and Alexander von Humboldt Foundation}
\email{tomek@math.idbsu.edu, http://math.idbsu.edu/\char 126 tomek}

\maketitle

%\CompileMatrices

\section{Introduction}
The purpose of this chapter is to discuss various results concerning
the relationship between measure and category. We are mostly
interested in set-theoretic properties of these ideals,
particularly, their cardinal characteristics.
This is a very large area, and it was necessary to make some
choices. We decided to present several new results and new approaches
to old problems. In most cases we do not present the optimal
result, but  a simpler theorem that still carries most
of the weight of that original result. 
For example, we construct Borel morphisms in the Cicho\'n diagram
while continuous ones can be constructed.
We believe however that the reader should have no problems upgrading
the material presented here to the current state of the art. 
The standard reference for this subject is \cite{BJbook}, and this
chapter updates it as most
of the material presented here was proved after \cite{BJbook} was published.

Measure and category have been studied for about a century. The beautiful
book \cite{Oxt80Mea} contains a lot of classical results, mostly from
analysis and topology, that involve these notions. 
The role played by Lebesgue measure and the Baire category in these
results is more or less identical. There are, of course, theorems
indicating lack of complete symmetry but they do not seem very
significant. For example, Kuratowski's theorem (cf. Theorem \ref{classic1}) asserts
that for every Borel function $f: \twoom \longrightarrow \twoom$ there
exists a meager set $F \subseteq \twoom$  such that $f \rest (\twoom
\setmin F)$ is continuous.
The dual proposition stating that for every Borel function $f: \twoom
\longrightarrow \twoom$ there exists a measure one set $G \subseteq
\twoom$ such that $ f \rest G$ is continuous is false. We only have a
theorem of Luzin which guarantees that a such $G$'s can have 
measure arbitrarily close to one.

The last 15 years have brought a wealth of results indicating that
hypotheses relating to measure are often stronger than the analogous 
ones relating to category.
This chapter contains several examples of this phenomenon. Before we
delve into this subject let us give a little historical background.
The first result of this kind is due to Shelah \cite{Sh176}. He
showed that
\begin{itemize}
\item If all projective sets are measurable then there exists an inner
  model with an inaccessible cardinal.
\item It is consistent with $\ZFCa$ that all projective sets have the
  property of Baire.
\end{itemize}

In 1984 Bartoszynski \cite{Bar84Add} and Raisonnier and Stern
in their \cite{RaiSte85Str} showed that additivity of measure is not greater
than additivity of category, while Miller \cite{Mil81Som} showed
that it can be strictly smaller.
In subsequent years several more results of that kind were found. Let
us mention one more (cf. \cite{BarGolJudMea}) concerning filters on $
\omega $ (treated as subsets of $ \twoom$):
\begin{itemize}
\item There exists a measurable filter that does not have the Baire
  property. In fact, every filter that has measure zero can be
  extended to a measure zero filter that does not have the Baire
  property. 
\item It is consistent with $\ZFCa$ that every filter that has the
  Baire property  is measurable.
\end{itemize}
All these results as well as many others concerning measurability and
the Baire property of projective sets, connections with forcing and
others can be found in \cite{BJbook}.
\section{Tukey connections}
The starting point for our considerations is the following list of cardinal
invariants of an ideal.
For an ideal $ {\mathcal J} $ of subsets of a set $X$ define
\begin{enumerate}
\item $\add({\mathcal J}) = \min\left\{ |{\mathcal A}| : {\mathcal A} \subseteq {\mathcal J}
\hbox{ and } \bigcup {\mathcal A} \not \in {\mathcal J}\right\}$,
\item $\cov({\mathcal J}) = \min\left\{ |{\mathcal A}| : {\mathcal A}
    \subseteq {\mathcal J} 
\hbox{ and } \bigcup {\mathcal A} = X\right\}$,
\item $\unif({\mathcal J}) = \min\left\{ |Y| : Y \subseteq X
\hbox{ and } Y \not \in {\mathcal J}\right\}$,
\item $\cof({\mathcal J}) = \min\left\{|{\mathcal A}| : {\mathcal A}
    \subseteq {\mathcal J} 
\hbox{ and } \forall B \in {\mathcal J} \ \exists A \in {\mathcal A} \ B
\subseteq A \right\}$  .
\end{enumerate}

\begin{definition}
  Suppose that $ {\mathcal P} $ and $ {\mathcal Q} $ are partial
  orderings. We say that $ {\mathcal P} \tles {\mathcal Q} $ if there
  is function $f: {\mathcal P} \longrightarrow {\mathcal Q} $ such
  that for every bounded set $X \subseteq {\mathcal Q} $, $f^{-1}(X)$
  is bounded in $ {\mathcal P} $. Such a function $f$ is called Tukey
  embedding.  Define $ {\mathcal P} \equiv {\mathcal Q} $ if $ {\mathcal P} \tles {\mathcal Q} $ and $ {\mathcal  Q} \tles {\mathcal P} $.
\end{definition}
Note that if $f: {\mathcal P} \longrightarrow {\mathcal Q} $ is a
Tukey embedding then there is an associated function $f^\star:
{\mathcal Q} \longrightarrow {\mathcal P} $ defined such that 
$f^\star(q)$ is a bound of the set $f^{-1}\lft1(\{(p): p \leq q\}\rgt1)$.
Observe that $f$ maps every set unbounded in $ {\mathcal P} $ onto a
set unbounded in $ {\mathcal Q} $ and $f^\star$ maps every set cofinal
in $ {\mathcal Q} $ onto a set cofinal in $ {\mathcal P} $.

\begin{lemma}\label{tukey}
Suppose that  ${\mathcal I}$ and ${\mathcal J}$ are two ideals. If
${\mathcal 
I} \tles {\mathcal J}$, then $\add({\mathcal I}) \geq
\add({\mathcal J})$ 
and $\cof({\mathcal I}) \leq \cof({\mathcal J})$.
\end{lemma}
\begin{proof}
Suppose that $f : {\mathcal I} \longrightarrow {\mathcal J}$ is a Tukey
function.

Let ${\mathcal A} \subseteq {\mathcal I}$ be a family of size $<
\add({\mathcal J})$. 
Find a set $B \in {\mathcal J}$ such that $\bigcup_{A \in {\mathcal A}}
f(A) \subseteq B$. It follows that $\bigcup {\mathcal A} \subseteq
f^\star(B)$.

Similarly, if ${\mathcal B} \subseteq {\mathcal J}$ is a basis for
${\mathcal J}$, then $\left\{f^\star(B) : B \in {\mathcal B}\right\}$
is a 
basis for ${\mathcal I}$.
\end{proof}

We will need a slightly stronger definition which will encompass both
cardinal invariants and Tukey embeddings. 

\begin{definition}\label{1.3}
  Suppose that ${\mathbf A}=(A_-,A_+,A)$, where $A$ is a binary relation
  between $A_-$ and $A_+$.
Let 
$$ {\mathfrak d}({\mathbf A})=\{Z: Z \subseteq A_+ \ \& \ \forall x \in A_- \
\exists z \in Z \ A(x,z)\}.$$
$$ {\mathfrak b}({\mathbf A})=\{Z: Z \subseteq A_- \ \& \ \forall y \in A_+ \
\exists z \in Z \ \neg A(z,y)\}.$$
$$\nor{{\mathbf A}}=\min\{|Z|: Z \in {\mathfrak d}({\mathbf A})\}.$$

Define ${\mathbf A}^\perp = (A_+,A_-,A^\perp)$, where
$A^\perp=\{(z,x): \neg A(x,z)\}$.
Note that ${\mathfrak b}({\mathbf A})= {\mathfrak d}({\mathbf A}^\perp)$.
\end{definition}

Note that $\nor{{\mathbf A}}$ is the smallest size of the
``dominating'' family in $A_+$ and $\nor{{\mathbf A}^\perp}$ is the
smallest size of the 
``unbounded'' family in $A_-$.
Virtually all cardinal characteristics of the continuum can be
expressed in this framework. 
For two ideals $ {\mathcal I} \subseteq {\mathcal J} $ of subsets of
$X$ we have: 
\begin{itemize}
\item $\cof({\mathcal J})=\nor{({\mathcal J},{\mathcal
      J},\subseteq)}$,
\item $\add({\mathcal J})=\nor{({\mathcal J},{\mathcal
      J},\subseteq)^\perp}=\nor{({\mathcal J},{\mathcal J},
      \not\supseteq)}$,
  \item $\cov({\mathcal J})=\nor{(X,{\mathcal J},\in)}$,
  \item $\unif({\mathcal J})=\nor{(X,{\mathcal
        J},\in)^\perp}=\nor{({\mathcal J},X,\not\ni)}$,
%  \item $\cof({\mathcal I}, {\mathcal J})=\nor{({\mathcal I},{\mathcal
%        J},\subseteq)}$,
%  \item $\add({\mathcal I},{\mathcal J})=\nor{({\mathcal I},{\mathcal
%        J},\subseteq)^\perp}=\nor{({\mathcal J},{\mathcal
%        I},\not\supseteq)}$,
  \item ${\mathfrak d}=\nor{(\omom, \omom, \leq^\star
      )}$,
    \item ${\mathfrak b}=\nor{(\omom, \omom, \leq^\star
      )^\perp}=\nor{(\omom, \omom, \not\geq^\star
      )}$, where for $f,g \in \omom$ we define $f \leq^\star g$ if
    $f(n) \leq g(n)$ for all but finitely many $n \in \omega $.
\end{itemize}

The notion of Tukey embedding generalizes to the following:
\begin{definition}
  A morphism $ \varphi $ between $ {\mathbf A}$ and ${\mathbf B}$ is a
  pair of functions
$ \varphi_-:A_- \longrightarrow B_-$ and $ \varphi_+: B_+
  \longrightarrow A_+$ such that for each $a \in A_-$ and $b \in B_+$,
$$A\lft1(a, \varphi_+(b)\rgt1), \hbox{ whenever } B\lft1(\varphi_-(a),b\rgt1).$$
If there is a morphism between $ {\mathbf A}$ and ${\mathbf B}$, we say
that  $ {\mathbf A}\tles {\mathbf B}$.
\end{definition}
Note that  if a pair of functions $f,f^\star$ witnesses that $ {\mathcal
  P} \tles {\mathcal Q} $, then $ \varphi =(f,f^\star)$ is a morphism
  between $({\mathcal P},{\mathcal P},\leq)$ and $({\mathcal
  Q},{\mathcal Q},\leq)$.
\begin{lemma}
  \begin{enumerate}
  \item  $ {\mathbf A} \leq {\mathbf B} \iff   {\mathbf A}^\perp \geq
    {\mathbf B}^\perp$,
  \item If $ {\mathbf A} \leq {\mathbf B}$, then
$\nor{{\mathbf A}} \leq \nor{{\mathbf B}}$ and 
$\nor{{\mathbf A}^\perp} \geq \nor{{\mathbf B}^\perp}$.
  \end{enumerate}
\end{lemma}
\begin{proof}
(1) If $ \varphi=(\varphi_-,\varphi_+)$ is a morphism between ${\mathbf
  A}$ and ${\mathbf B}$, then 
$ \varphi^\perp=(\varphi_+,\varphi_-)$ is a morphism between ${\mathbf
  B}^\perp$ and ${\mathbf A}^\perp$.

(2) Suppose that $Z \in {\mathfrak d}({\mathbf B})$ is
such that 
$|Z|=\nor{{\mathbf B}}$. Then $\{\varphi_+(z): z\in Z\}$ is cofinal in
$A_+$. In other words, $\nor{{\mathbf A}}\leq |Z|$.
\end{proof}

For two Polish spaces $X,Y$ (i.e. metric, separable with no isolated
points) define $\borel(X,Y)$ to be the space of all Borel functions
from $X$ to $Y$.

Given relation ${\mathbf A}$ and assuming that both $A_-$ and $A_+$
are Polish spaces we define families of small sets of reals
as:
$$\D({\mathbf A})=\left\{X \subseteq \reals: \forall f \in
\borel(\reals,A_+) \ f{``}(X) \not \in {\mathfrak d}({\mathbf A})\right\}$$
and
$$\B({\mathbf A})=\left\{X \subseteq \reals: \forall f \in
\borel(\reals,A_-) \ f{``}(X) \not \in {\mathfrak b}({\mathbf A})\right\}$$
In other words, $\D({\mathbf A})$
consists of sets of reals whose Borel images are not ``dominating''
and $\B({\mathbf A})$ consists of sets whose Borel images are
``bounded''. 

\begin{lemma}
  \begin{enumerate}
  \item $\unif\lft1(\D({\mathbf A})\rgt1)=\nor{{\mathbf A}}$ and
    $\unif\lft1(\B({\mathbf A})\rgt1)=\nor{{\mathbf A}^\perp}$. 
  \item If there exists a Borel morphism from ${\mathbf A}$ to ${\mathbf B}$, then
$\B({\mathbf B}) \subseteq \B({\mathbf A})$ and $\D({\mathbf A})
\supseteq \D({\mathbf B})$. 
  \end{enumerate}
\end{lemma}
\begin{proof}
  (1) Clearly $\unif(\D({\mathbf A}) \geq \nor{{\mathbf A}}$. To show
  the other inequality notice that there is a Borel mapping from
  $\reals$ onto $A_+$.

\bigskip

(2) Suppose that $X \not\in \B({\mathbf A})$ and let $f:\reals \longrightarrow
A_-$ be a Borel map such that $f{``}(X)\in {\mathfrak b}({\mathbf A})$.
It follows that $ \varphi_- \comp f{``}(X)\in {\mathfrak b}({\mathbf
  B})$. Since $ \varphi_-\comp f$  is a Borel mapping it
follows that $X \not \in \B({\mathbf B})$.
\end{proof}

For cardinals $ \kappa = \nor{{\mathbf A}}$ and $ \lambda
=\nor{{\mathbf B}}$ the question whether the inequality $ \kappa  \leq
\lambda $ is provable in $\ZFCa$ leads naturally to the question
whether $ {\mathbf A} \tles {\mathbf B}$ and $\D({\mathbf A}) \subseteq
\D({\mathbf B})$. 
Even though these questions are more general, in most cases the proof
that $ \kappa \leq \lambda $ yields 
$ {\mathbf A} \tles {\mathbf B}$. Moreover, the existence of a Borel
morphism witnessing that $ {\mathbf A} \tles {\mathbf B}$ uncovers the
combinatorial aspects of these problems.

\bigskip

{\bf Historical remarks}
Tukey embeddings were defined in \cite{Tukeyart} and further studied in
\cite{Isbell}. 
In context of the orderings considered here see \cite{Fre91Fam},
\cite{FrePa} and \cite{LoVe}. 

The framework used in the definition \ref{1.3} is due to Vojt\'a\v s
\cite{Vojtas};
the particular formulation used here comes from 
\cite{Blass96}.

\section{Inequalities provable in ZFC}
The notions defined in the previous section 
 are quite general. The focus of this chapter is on the ideal
 of meager sets ($\M$) and measure zero (null)
sets ($\N$) with respect to the standard
product measure on $\mu$ on $\twoom $ or the Lebesgue measure $\mu$
on $\reals$.

For an ideal $ {\mathcal J} $, we  identify a Borel
mapping $H: \reals 
\longrightarrow {\mathcal J}$  with a Borel set $H \subseteq \reals
\times 
\reals$ in such a way that 
\begin{enumerate}
\item $H(x)= (H)_x=\lft1\{y: (x,y) \in H\rgt1\}.$
\item $H$ is a Borel ${\mathcal
   J}$-set, that is, $(H)_x \in {\mathcal J}$ for all $x \in \reals.$
\end{enumerate}

Using this terminology we define the following classes of small sets:

 \begin{multline*}\COF(\N)=\D(\N,\N,\subseteq)=\\
\lft1\{X \subseteq \reals: $ for
  every Borel $\N$-set $H$, $\{(H)_x: x \in X\}$ is not a basis of
  $\N\rgt1\},
\end{multline*}
\begin{multline*}
\ADD(\N)=\B(\N,\N,\subseteq)=\\
\lft1\{X \subseteq \reals: $ for
  every Borel $\N$-set  $H$,  $\bigcup_{x \in X} (H)_x \in \N\rgt1\},\end{multline*}
  \begin{multline*}\COV(\N)=\D(\reals,\N,\in)=\\
\lft1\{X \subseteq \reals: $ for every Borel $\N$-set $H$,
$\bigcup_{x \in X} (H)_x \neq \reals\rgt1\},\end{multline*} 
  
$\UNIF(\N)=\B(\reals,\N,\in)=
\lft1\{X \subseteq \reals:$ every
    Borel image of $X$ is in $\N \rgt1\},$

    $\D=\D(\omom, \omom, \leq^\star)$,

     $\B=\B(\omom, \omom, \leq^\star)$.

In the same way we define $\ADD(\M)$, $\COV(\M)$, etc.

Instead of dealing with all null and meager sets we need to consider
only  suitably chosen cofinal families.

\begin{enumerate}
\item $A \in \N$ if and only if there exists a family of basic open
  sets $\{C_n: n \in \omega\}$ such that 
$\sum_{n=0}^\infty \mu(C_n) < \infty$ and
$A \subseteq \bigcap_{n \in \omega} \bigcup_{m>n} C_m$,
\item $A \in \M$ if and only if there is a family  of $\{F_n: n \in
  \omega\}$ of closed nowhere dense sets such that $A \subseteq
  \bigcup_{n \in \omega} F_n$.
\end{enumerate}
In particular every null set can be covered by a null set of type
$G_\delta$ and every meager set can be covered by a meager set of type
$F_\sigma$.

\begin{definition}
 $\clopen=\{C^{{n}}_m: n,m\in \omega \}$ is a family of clopen
subsets of $\twoom $ such that $\mu(C^n_m)=2^{-n}$ for each $m$.
\end{definition}

\begin{lemma}\label{charnull}
$A \in \N \iff  \exists f \in \omom \ \left(A \subseteq
    \bigcap_m \bigcup_{n>m} C^n_{f(n)}\right)$.
\end{lemma}
\begin{proof}
  ($ \leftarrow$) Note that the set 
$ \bigcup_{n>m} C^n_{f(n)}$ has measure at most $2^{-m}$.

\bigskip

($ \rightarrow$) 
For an open set $U \subseteq \twoom $ let 
$$\widetilde{U}=\lft2\{t \in \twolom : [t] \subseteq U \ \& \ \forall s
\subsetneq t \ \lft1([s] \not \subseteq U\rgt1)\rgt2\}.$$
Note that $\widetilde{U}$ is a canonical representation of $U$ as a
union of disjoint basic intervals.

Find open sets $\{G_n: n \in \omega \}$ covering $A$
such that $\mu(G_n) \leq 2^{-n}$.
Let $\{t_n: n \in \omega\}$ be the lexicographic enumeration  of
$\bigcup_{n \in \omega} \widetilde{G}_n$.
Define for $ n \in \omega $,
$$h(n+1)=\min\left\{k> h(n):\sum_{j=k}^\infty \mu\lft1([t_j]\rgt1) \leq
  2^{-n}\right\},$$
and let 
$$D_n=\bigcup_{j=h(n)}^{h(n+1)} [t_j].$$
Let $f \in \omom $ be such that $D_n=C^n_{f(n)}$ for each $n$.
\end{proof}

\begin{definition}\label{gooddef}
Let $\{U_n: n \in \omega \}$ be a basis in $\twoom $ and let
${\mathcal S}=\{S^n_m: n,m \in \omega \}$ be any family of clopen sets.
We say that  ${\mathcal S}$ is {\em good} if
\begin{enumerate}
\item $S^n_m \cap U_n\neq \emptyset$ for $ m\in \omega $,
\item For any open dense set $U \subseteq \twoom $ and $ n \in
\omega $ there is $m$ such that $S^n_m \subseteq U$.
\end{enumerate}
\end{definition}

\begin{lemma}\label{two}
Suppose that the family ${\mathcal S}=\{S^n_m: n,m\in \omega\}$ is  good. Then  
 $$A \in \M \iff  \exists f \in \omom \ \left(A \subseteq
    \twoom \setmin \bigcap_m \bigcup_{n>m} S^n_{f(n)}\right).$$
\end{lemma}
\begin{proof}
  ($\leftarrow$)
Note that the set $\bigcup_{n>m} S^n_{f(n)}$ is open and dense for
every $m$.

\bigskip

($\rightarrow$) 
Let $\<F_n: n \in \omega\>$ be an increasing sequence of closed
nowhere dense  sets
covering $A$.
For each $n$ let
$$f(n)=\min\{m: U_n \cap S^n_m \cap F_n = \emptyset\}.$$
It is clear that 
$\bigcup_n F_n \cap \bigcap_m \bigcup_{n>m} S^n_{f(n)}=\emptyset$.
\end{proof}

Define master sets $N,M \subseteq \omom \times \twoom $ by
$$N= \bigcap_m \bigcup_{n>m} \bigcup_{f \in \omom }
\{f\}\times C^n_{f(n)}$$
and
$$M= (\omom \times \twoom) \setmin \bigcap_m \bigcup_{n>m}
\bigcup_{f \in \omom } 
\{f\}\times S^n_{f(n)}$$

Note that $N$ is a $G_\delta $ set while $M$ is an
$F_\sigma$ set. Moreover, $\left\{(N)_f: f \in \omom\right\}$ is cofinal in
$\N$ and $\left\{(M)_f: f \in \omom\right\}$ is cofinal in
$\M$. 

The following lemma shows that the representation of meager sets
 does not depend on the choice of good family:
\begin{lemma}
  Suppose that ${\mathcal S}=\{S^n_m: n,m\in \omega\}$ and ${\mathcal
  T}=\{T^n_m: n,m\in \omega\}$ are   good and $M$ and $\overline{M}$
  are associated master sets.
Then there are Borel mappings $ \varphi_-,\varphi_+: \omom
  \longrightarrow \omom $ such that 
$$(M)_f \subseteq (M)_{\varphi_+(g)}, \hbox{ whenever }
(\overline{M})_{\varphi_-(f)} \subseteq (\overline{M})_g.$$
\end{lemma}
\begin{proof}
  For $f,g \in \omom $ and $ n \in \omega $ define
$$ \varphi_-(f)(n)=\min\left\{k: T^n_k \subseteq \bigcup_{m\geq n}
  S^m_{f(m)} \right\}$$
and
$$ \varphi_+(g)(n)=\min\left\{k: S^n_k \subseteq \bigcup_{m\geq n}
  T^m_{g(m)} \right\}.$$
We leave it to the reader to verify that these mappings have the
  required properties.
\end{proof}

The following two theorems will be helpful in many subsequent
constructions.
\begin{theorem}\label{classic}
  Suppose that $H \subseteq \twoom \times \twoom $ is a Borel set.
  \begin{enumerate}
  \item $\{x: (H)_x \in \N\}$ is Borel,
  \item $\{x: (H)_x \in \M\}$ is Borel,
  \item If $U$ is open and $(H)_x$ is compact for every $x$, then
$\{x: U \cap (H)_x = \emptyset\}$ is Borel,
\item If for every $x$, $(H)_x$ is ``large,'' where large is either
  ``of positive measure'' or ``nonmeager'',  then there exists a Borel
  function $f: \twoom \longrightarrow \twoom $ such that 
  for every $x$, $f(x) \in (H)_x$.
  \end{enumerate}
\end{theorem}
\begin{proof}
  See \cite{Kechris95} 16.A for (1) and (2), 18.B for (4).
\end{proof}
\begin{theorem}\label{classic1}
  If $X$ and $Y$ are Polish spaces and $f: X \longrightarrow Y$ is a Borel
  mapping then there is a dense $G_\delta $ set $G \subseteq X$ such
  that $f\rest G$ is continuous.
\end{theorem}
\begin{proof}
  This is a special case 
of a theorem of Kuratowski; see \cite{Kechris95} 8.I.
\end{proof}
Lemmas  \ref{charnull} and \ref{two} have their two-dimensional
analogs. 
\begin{lemma}\label{firstlem}
  The following conditions are equivalent for a Borel set $H \subseteq
  \twoom \times \twoom $:
  \begin{enumerate}
  \item  $\forall x \ \lft1((H)_x \in \N\rgt1)$,
  \item For every $\varepsilon>0$
there exists a Borel set $B \subseteq \twoom \times
\twoom$ such that 
\begin{enumerate}
\item $H \subseteq B$,
\item for every $x$, $(B)_x$ is an open set of measure $< \varepsilon$.
\end{enumerate}
\item There exists a Borel function $x \rightsquigarrow f_x$ such that 
$$ \forall x \ \lft1((H)_x \subseteq (N)_{f_x}\rgt1).$$
  \end{enumerate}
\end{lemma}
\begin{proof}
  (2)$\rightarrow$(3) Let $\{B_n: n \in \omega \}$ be a family of
  Borel sets such that 
\begin{enumerate}
\item $H \subseteq \bigcap_n B_n$,
\item For every $x$, $(B_n)_x$ is an open set of measure $< 2^{-n}$.
\end{enumerate}
Look at the proof of the Lemma~\ref{charnull} to see that for each $x$, $(B)_x
=(N)_{f_x}$ and the mapping $x \rightsquigarrow f_x$ is Borel.

\bigskip

(3) $\rightarrow$(1) is obvious.

\bigskip

(1)$\rightarrow$(2) By induction on complexity we show that for every
$ \varepsilon > 0$ and  a Borel set $H \subseteq \twoom \times
\twoom $ there exists a Borel set $B \supseteq H$ such that for
every $x$, 
$(B)_x$ is open  and  $\mu\lft1((B\setmin
H)_x\rgt1)<\varepsilon$.
The only nontrivial part is to show that if the theorem holds for sets in
$ { \Sigma}^0_\alpha $, then  it holds for any set $A \in 
 {  \Pi}^0_\alpha$.
To see this write $A=\bigcap_n A_n$ where $\<A_n: n \in \omega \>$ is
a descending sequence of sets in $ { \Sigma}^0_\alpha $.
For each $n$ let $B_n$ be the set obtained from induction hypothesis
for $A_n$ and $ \varepsilon/2$.
Let $K^n=\left\{x: \mu\lft1((A_n\setmin A)_x\rgt1) < \varepsilon/2
\right\}$. Each set $K^n$ is Borel.
Now define
$$B= B_0\cap (K_0\times \twoom) \cup \bigcup_{n \in \omega} B_{n+1} \cap \left((K^{n+1} \setmin K^{n})\times
\twoom\right).$$ 
\end{proof}
\begin{lemma}
  The following conditions are equivalent for a Borel set $H \subseteq
  \twoom \times \twoom $:
  \begin{enumerate}
  \item $ \forall x \ \lft1((H)_x \in \M\rgt1)$,
\item  There
exists a family of Borel sets $\{G_n :n \in \omega\} \subseteq
\twoom \times \twoom $ such that
\begin{enumerate}
  \item $(G_n)_x$ is a closed nowhere dense set for all $x \in
\twoom$,
\item $H \subseteq \bigcup_{n \in \omega} G_n$.
\end{enumerate}
\item There exists a Borel function $x \rightsquigarrow f_x$ such that 
$$ \forall x \ \lft1((H)_x \subseteq (M)_{f_x}\rgt1).$$
  \end{enumerate}
\end{lemma}
\begin{proof}
   (1)$\rightarrow$(2) By induction on the complexity we show
that for any Borel set $H \subseteq \twoom \times \twoom $ there
are Borel sets $B$ and $\{F_n: n \in \omega \}$ such that 
\begin{enumerate}
\item $(B)_x$ is open for every $x$,
\item $(F_n)_x$ is closed nowhere dense for every $x$ and $n$,
\item $H \triangle B \subseteq B \cup \bigcup_n F_n.$
\end{enumerate}
As before the nontrivial part is to show the theorem for the class
${\Pi}^0_\alpha$ given that it holds for ${\Sigma}^0_\alpha $. 
Suppose that $A \in {\Sigma}^0_\alpha $ and $B$ is the set obtained by
applying the inductive hypothesis to $A$.
Let $\<U_n: n \in \omega\>$ be an enumeration of basic sets in
$\twoom $.
Define for $ n \in \omega $,
$$Z_n =\{x : U_n \cap (B)_x =\emptyset\}.$$
Note that sets $Z_n$ are Borel.
Let 
$B'=\bigcup_{n} Z_n\times U_n.$
The vertical sections of the 
set $F=\twoom \times \twoom \setmin (B \cup B')$ are closed
nowhere dense and $(\twoom \setmin A) \triangle B' \subseteq F$,
which ends the proof.

\bigskip

(2) $ \rightarrow$ (3)
For $ x \in \twoom $ let 
$$f_x(n)=\min\rgt2\{k: \forall i \leq n  \ \lft1(S^n_k \cap (G_i)_x =
\emptyset\rgt1)\rgt2\}.$$
\end{proof}

From these two lemmas it follows that:
\begin{lemma}\label{before8}
  Let $ {\mathcal I}$ be $\N$ or $\M$ and let $I$ be the associated
  master set. Then for $X \subseteq \reals$:
\begin{enumerate}
\item $X \in \ADD({\mathcal I}) \iff \forall F \in \borel(\reals, \omom)
  \ \exists f \in \omom \ \forall x \in X \ \lft1((I)_{F(x)}
  \subseteq (I)_f\rgt1),$
\item $X \in \COF({\mathcal I}) \iff \forall F \in \borel(\reals, \omom)
  \ \exists  f \in \omom \ \forall  x \in X\  \lft1((I)_{f} \not
  \subseteq (I)_{F(x)}\rgt1),$
\item $X \in \COV({\mathcal I}) \iff \forall F \in \borel(\reals, \omom)
  \ \exists z \ \forall x \in X \ \lft1(z \not\in (I)_{F(x)}\rgt1)$,
\item $X \in \UNIF({\mathcal I}) \iff \forall F \in \borel(\reals, \omom)
  \ \exists f \ \forall x \in X \ \lft1(F(x) \in (I)_{f}\rgt1)$.
\end{enumerate}
\end{lemma}

The goal of this section is to establish:
\begin{theorem}\label{8}
\onesetup {(\reals,\N,\in)} {(\M,\reals, \not\ni)} {(\M,\M,\subseteq)}
{(\N,\N, \subseteq)} {(\omom, \omom, \not\geq^\star)}

\twosetup {(\omom, \omom, \leq^\star)}  {(\N,\N,
\not \supseteq)} {(\M,\M, \not \supseteq)} {(\reals,\M,\in)}
{(\N,\reals,\not\ni)}

$$\diagram
\pone \rto^\tles & \ptwo \rto^\tles & \pthree \rto^\tles& \pfour\\
& \pfive \rto^\tles \uto^\tles&\psix \uto^\tles& \\
\pseven\uuto^\tles\rto^\tles&\peight \rto^\tles\uto^\tles&\pnine
\rto^\tles\uto^\tles&\pten\uuto^\tles 
\enddiagram$$

As a consequence we will get the following two diagrams:

\onesetup \COV(\N) \UNIF(\M) \COF(\M) \COF(\N) {\B}

\twosetup {\D} \ADD(\N) \ADD(\M) \COV(\M) \UNIF(\N)

$$\diagram
\pone \rto^\subseteq & \ptwo \rto^\subseteq & \pthree \rto^\subseteq& \pfour\\
& \pfive \rto^\subseteq \uto^\subseteq&\psix \uto^\subseteq& \\
\pseven\uuto^\subseteq\rto^\subseteq&\peight
\rto^\subseteq\uto^\subseteq&\pnine 
\rto^\subseteq\uto^\subseteq&\pten\uuto^\subseteq 
\enddiagram$$

\onesetup \cov(\N) \unif(\M) \cof(\M) \cof(\N) {{\mathfrak b}}

\twosetup {{\mathfrak d}} \add(\N) \add(\M) \cov(\M) \unif(\N)

$$\diagram
\pone \rto^\leq & \ptwo \rto^\leq & \pthree \rto^\leq& \pfour\\
& \pfive \rto^\leq \uto^\leq&\psix \uto^\leq& \\
\pseven\uuto^\leq\rto^\leq&\peight \rto^\leq\uto^\leq&\pnine
\rto^\leq\uto^\leq&\pten\uuto^\leq 
\enddiagram$$
The last of these diagrams is called the Cicho\'n diagram.
\end{theorem}

It is enough to find the following morphisms:
\begin{enumerate}
\item ${(\reals,\N,\in)} \tles {(\M,\reals, \not\ni)}$,
\item ${(\M,\M,\subseteq)} \tles {(\N,\N,\subseteq)}$,
\item ${(\M,\M,\not\supseteq)}\tles {(\omom,
    \omom, \not\geq^\star)}$,
\item ${(\omom, \omom, \not\geq^\star)} \tles
  {(\M,\reals,\not\ni)} $,
\item ${(\N,\N,
\not \supseteq)} \tles {(\reals,\N,\in)}$. 
\end{enumerate}
The remaining morphisms are dual to those listed above. In each case
we will find a Borel morphism.
Note that thanks to master sets $M$ and $N$ defined earlier, Borel
morphisms between these structures can be interpreted as the
automorphisms of the index set i.e. $ \omom $.

\begin{theorem}\label{main}
  $\M \tles \N$; there are two Borel functions $ \varphi_-,\varphi_+:
  \omom \longrightarrow \omom $ such that 
$$(M)_f \subseteq (M)_{\varphi_+(g)}, \hbox{ whenever }
(N)_{\varphi_-(f)} \subseteq (N)_g.$$
In particular, $\ADD(\N) \subseteq \ADD(\M)$ and $\COF(\M) \subseteq
\COF(\N)$, $\add(\N) \leq \add(\M)$ and $\cof(\M) \leq \cof(\N)$.
\end{theorem}
\begin{proof}
  Let 
$$ {\mathcal C}=\left\{S \in {}^\omega (^{<\omega}[\omega]) :
  \sum_{n=1}^\infty \frac{|S(n)|}{2^n}< \infty\right\}.$$
For $S,S' \in {\mathcal C}$ define $S \subseteq^\star S'$ if for all
  but finitely many $n$, $S(n) \subseteq S'(n)$.
  \begin{lemma}\label{firstlemma}
    $\N \equiv {\mathcal C}$. 
  \end{lemma}
  \begin{proof}
    To see that $\N \tles {\mathcal C}$ define
$ \varphi_-:\omom \longrightarrow {\mathcal C}$ and $
\varphi_+: {\mathcal C} \longrightarrow \omom $ such that for
$f \in \omom $ and $S \in {\mathcal C}$ we have
$$(N)_f \subseteq (N)_{\varphi_+(S)}, \hbox{ whenever } \varphi_-(f)
\subseteq^\star S.$$
For $ f \in \omom $ put $ \varphi_-(f) =h$, where
$h(n)=\{f(2n),f(2n+1)\}$. 
If $S \in {\mathcal C}$ let $ \varphi_+(S)=g\in \omom $ be such
that 
$$C^n_{g(n)}=\bigcup_{k\in S(n)} C^{2n}_k \cup \bigcup_{k\in S(n)}
C^{2n+1}_k.$$ 
Verification that these mappings have the required properties is
straightforward.

\bigskip

To show that ${\mathcal C} \tles \N$ we will find Borel functions $
\varphi_-: {\mathcal C} \longrightarrow \omom $ and $
\varphi_+: \omom \longrightarrow {\mathcal C}$ such that for
$S \in {\mathcal C}$ and $f \in \omom $,
$$S \subseteq^\star \varphi_+(f) \hbox{ whenever } (N)_{\varphi_-(S)}
\subseteq (N)_f.$$
Let $\{G^n_m: n,m\in \omega\}$ be a family of clopen probabilistically
independent sets such that $\mu(G^n_m)=2^{-n}$.
For $S \in {\mathcal C}$ define
$ \varphi_-(S)=f \in \omom $ such that 
$$\bigcap_{m\in \omega } \bigcup_{n>m}\bigcup_{k \in S(n)} G^{n}_k
\subseteq (N)_f.$$   
First consider $H' \subseteq {\mathcal C}\times \twoom $ defined by
$(H')_S=\bigcap_{m\in \omega } \bigcup_{n>m}\bigcup_{k \in S(n)} G^{n}_k$
for $S \in {\mathcal C}$. Note that $H'$ is a Borel set and
$(H')_S $ has measure zero for every $S$. Fix a Borel isomorphism $a:
{\mathcal C} \longrightarrow \omom$ and let $H \subseteq \omom\times
\twoom$ be defined as
$(H)_{a(S)}=(H')_S$ for $S\in {\mathcal C}$. Apply \ref{firstlem} to
find a Borel mapping $x \rightsquigarrow f_x$ such that $(H)_x
\subseteq (N)_{f_x}$ and define $\varphi_-(S)=f_{a(S)}$.

To define $ \varphi_+: \omom \longrightarrow {\mathcal C}$ we
proceed as follows.
Find a Borel set  $K \subseteq \omom
\times \twoom $ such 
that
\begin{enumerate}
\item $(K)_f$ is a compact set of measure $\geq 1/2$ for all $f \in
  \omom $.
\item $N \cap K=\emptyset$.
\item For any basic open set $U \subseteq \twoom $ and $f \in
  \omom $, if $U \cap (K)_f \neq \emptyset$ then 
$U \cap (K)_f$ has positive measure.
\end{enumerate}
First use Lemma~\ref{firstlem} to find a set $K'$ satisfying the first two conditions.
Let $\<U_j: j\in \omega \>$ be an enumeration of basic open sets in
$\twoom $. 
For each $j$ let $Z_j =\{f: \mu(U_j \cap (K')_f)=0\}$. By Theorem
\ref{classic}, the sets $Z_j$ are
Borel for each $j$. 
Define $K=K' \setmin \bigcup_j (Z_j \times U_j)$.

For $f \in \omom, j,n \in \omega$ define
$$S^f_{j}(n)=\left\{i \in \omega : (K)_f \cap U_j \neq \emptyset \ \&\ (K)_f
\cap U_j \cap G^n_i = \emptyset\right\}.$$
Note that
$$0<\mu((K)_f \cap U_j)\leq \prod_n \prod_{i \in S^f_{j}(n)}
\mu(\twoom \setmin G^n_i).$$
Thus
$$0< \prod_{n=1}^\infty  \left(1-\frac{1}{2^n}\right)^{|S_j^f(n)|}.$$
It follows that
$$ \sum_{n=1}^\infty \frac{|S_j^f(n)|}{2^n} < \infty,$$
so $S^f_j \in {\mathcal C}$ for each $j$.
Moreover, the mapping $f \rightsquigarrow \<S^f_j: j\in \omega \>\in
{^\omega{\mathcal
  C}} $ is Borel (by Theorem \ref{classic}(3)).
Fix a Borel mapping from $^\omega {\mathcal C} $ to ${\mathcal C}$ such
that $ \<S_j^f: j\in \omega \> \rightsquigarrow S_\infty^f$ such that 
$$ \forall j \ \forall^\infty n \ S_j^f(n) \subseteq S_\infty^f(n).$$

Finally define $ \varphi_+$ by the formula:
$$ \varphi_+(f)(n)=S^f_\infty(n).$$
Suppose that for some $S \in {\mathcal C}$, $(N)_{\varphi_-(S)}
\subseteq (N)_f$. 
It follows that,
$$(K)_f \cap \bigcap_m \bigcup_{n>m}\bigcup_{k\in S(n)}
G^n_k =\emptyset.$$
By the Baire category theorem, there is a basic open set $U_j$ and
  $m_0 \in \omega $ such
  that $U_j \cap (K)_f \neq \emptyset$ but
$$(K)_f \cap U_j \cap \bigcup_{n>m_0}\bigcup_{k\in S(n)}
 G^n_k =\emptyset.$$
Therefore
$$ \forall^\infty n \ S(2) \subseteq
S^f_j(n) \subseteq S^f_\infty(n)=\varphi_+(f)(n),$$
which finishes the proof.
  \end{proof}
  \begin{lemma}\label{secondlemma}
    $\M \tles {\mathcal C}$; there are Borel mappings $ \varphi_-:
    \omom \longrightarrow {\mathcal C}$ and $ \varphi_+:
    {\mathcal C} \longrightarrow \omom $ such that for any $f
    \in \omom $ and $S \in {\mathcal C}$,
  $$(M)_f \subseteq (M)_{\varphi_+(S)} \hbox{ whenever } \varphi_-(f)
  \subseteq^\star S.$$
  \end{lemma}
  \begin{proof}
    We will need the following lemma:
    \begin{lemma}\label{2.10}
      There exists a good family $\{S^n_m:n,m \in \omega\}$ such that 
$$ \forall X \in {[\omega]^{\leq 2^n}} \ \left(\bigcap_{j \in X}
        S^n_j \neq \emptyset\right).$$
    \end{lemma}
    \begin{proof}
Fix $ n \in \omega $.
      Let $\< C_{m} : m \in \omega\>$ be an enumeration of all 
clopen sets. For $k \in \omega$ define
$$A_{k} = \left\{ l>k : C_{l} \cap \bigcap_{i \in I} C_{i} \cap U_n\neq
\emptyset \hbox{ whenever } I \subseteq k+1 \hbox{ and }
U_n\cap \bigcap_{i \in I} C_{i} \neq \emptyset\right\} .$$
Consider family 
$${\mathcal S}_n = \left\{ \bigcup_{i \leq 2^n} C_{m_{i}} :
m_{0} \in \omega \text{ and } \ m_{i+1} \in A_{m_{i}} \hbox{ for } i \leq
2^n\right\} .$$

We have to check that ${\mathcal S}_n$ satisfies conditions of
Definition~\ref{gooddef}.

$(1)$ Let $U$ be a dense open subset of $^{\omega}2$. 
Note
that $A_{k} \cap \left\{l \in \omega: U_n \cap C_{l} \subseteq
  U\right\} \neq
\emptyset$ for every $k \in \omega$, by  the density of $U$. 

Now define by induction a sequence
$\left\{m_{i}: i \leq 2^n\right\}$ such that $C_{m_{i}} \subseteq U$
and $m_{i+1} \in A_{m_{i}}$ for $i < 2^n$. Clearly
$U \supseteq \bigcup_{i \leq 2^n}C_{m_{i}} \in {\mathcal S}_n .$

\vspace{0.1in}

$(2)$ Suppose that $V_{1},  V_{2},  \ldots,  V_{2^n} \in
{\mathcal S}_n$. For any $j \leq 2^n$, $V_{j}=\bigcup_{i \leq
2^n} C_{m_{i}^{j}}$, where $m^{j}_{i} \in A_{m^{j}_{i}}$
for $i,j \leq 2^n$. Order the sets $V_{j}$ in such a way that
$m^{i}_{i} \leq m^{j}_{i}$ for $i \leq j \leq 2^n$.

It is easy to show by induction that $\bigcap_{j \leq 2^n}
V_{j} \supseteq \bigcap_{j \leq 2^n} C_{m^{j}_{j}} \neq
\emptyset$.
    \end{proof}

Let ${\mathcal S}=\bigcup_n {\mathcal S}_n=\{S^n_m:n,m\in \omega\}$. 
For $ f \in \omom$ define $ \varphi_-(f) = f \in {\mathcal C}$.
For $S \in
 {\mathcal C}$ let $ \varphi_+(S)=f \in \omom $ be such that 
$$(M)_f \supseteq  \twoom  \setmin  \bigcap_{m \in \omega}
\bigcup_{n>m} \bigcap_{i \in S(n)} S^{n}_{i} .$$
Since $|S(n)| \leq 2^n$ for all but finitely many $n$, by
Lemma~\ref{2.10}, 
$$  \emptyset \neq U_n \cap \bigcap_{i \in S(n)} S^{n}_{i}.$$ 
Now suppose that $ \varphi_-(f) \subseteq^{\star} S$.
This assumption means that there exists $n_{0} \in \omega$ such that 
$ f(m)
\in S(m)$ for $m \geq n_{0}$.
It follows that
$$(M)_{\varphi_+(S)}
 \supseteq  \twoom \setmin  \bigcap_{m \in \omega}
\bigcup_{n>m} \bigcap_{i \in S(n)} S^{n}_{i} \supseteq \twoom
\setmin  \bigcap_{m\in \omega} \bigcup_{n>m} S^{n}_{f(n)}.$$
  \end{proof}
Theorem~\ref{main} follows immediately; compose the morphisms
constructed in Lemma~\ref{firstlemma} and Lemma~\ref{secondlemma}.
\end{proof}

\begin{theorem}
  $(\reals,\N,\in) \tles (\M,\reals,\not\ni)$; there are Borel
  functions $ \varphi_-,\varphi_+: \reals \longrightarrow
  \omom $ such that for $x,y \in \reals$,
$$x \in (N)_{\varphi_+(y)} \hbox{ whenever } y \not \in
(M)_{\varphi_-(x)}.$$ 
In particular, $\COV(\N) \subseteq \UNIF(\M)$ and $\COV(\M) \subseteq
\UNIF(\N)$, $\cov(\N) \leq \unif(\M)$ and $\cov(\M) \leq \unif(\N)$.
\end{theorem}
\begin{proof}
  Let $B$ be a $G_\delta $  null set whose complement is meager. Use
  Theorem \ref{classic}(3) and Theorem \ref{classic1} to find
Borel functions $ \varphi_-,\varphi_+: \reals \longrightarrow
  \omom $ such that
$$ \forall x \ \left(B+x \subseteq (N)_{\varphi_-(x)}\right) \hbox { and } \forall
y \ \left(\twoom \setmin (B+y) \subseteq (M)_{\varphi_+(y)}\right).$$
If $y \not \in (N)_{\varphi_-(x)}$ then $y \not\in B+x$. It follows
that $x \in \twoom \setmin (B+y) \subseteq (M)_{\varphi_+(y)}$.
\end{proof}

\begin{theorem}\label{lastlem}
  $(\M,\M, \not \supseteq)\tles (\omom, \omom,
  {}^\star\not\geq)$;
there are Borel mappings $ \varphi_-,\varphi_+: \omom
  \longrightarrow \omom $ such that for $f, g \in
  \omom $
$$(M)_f \not \supseteq (M)_{\varphi_+(g)}, \hbox{ whenever }
\varphi_-(f) {}^\star\not\geq g.$$
In particular, $\D \subseteq \COF(\M)$ and $\ADD(\M) \subseteq
\B$, ${\mathfrak d} \leq \cof(\M)$ and $\add(\M) \leq {\mathfrak b}$.
\end{theorem}
\begin{proof}
Let ${\mathcal S}_n$ be the family of clopen sets $C$ such that
there exists $k>n$ and $s \in {^{[n,k)}}2 $ such that 
$$C=\{x \in \twoom : x \rest [n,k)=s\}.$$
Note that the family ${\mathcal S}=\bigcup_n {\mathcal S}_n$ is good
(given the appropriate choice of the sequence $\{U_n: n \in \omega\}$).

For $f \in \omom $ let $ \varphi_-(f)(n)=k$ if and only if
$\dom(S^n_{f(n)})=[n,k)$.

For a strictly increasing function $g \in \omom $ define
$\varphi_+(f)=h\in \omom $ such that 
$$(M)_h=\{x \in \twoom : \forall^\infty n \ \exists i \in \lft1[n,f(n)\rgt1) \
\lft1(x(i) \neq 0\rgt1)\}.$$
Note that the image of $ \omom $ under $ \varphi_+$ is rather
small, $ {\varphi_+}{``}(\omom)$ is not even cofinal in $\M$.

To finish the proof it is enough to show that if $
\varphi_-(f)(n)<g(n)$ for 
infinitely many $n$, then
\begin{multline*}
\left\{x \in \twoom : \forall^\infty n \ \exists i \in \lft1[n,g(n)\rgt1) \
x(i) \neq 0\right\}\not \subseteq \\ 
\left\{x \in \twoom :
\forall^\infty n \ x 
\rest \lft1[n,\varphi_-(f)(n)\rgt1) \neq S^n_{f(n)} \right\}.
\end{multline*}

Find a sequence $\{n_k: k\in \omega \}$ such that  for all $k$,
$$n_k< \varphi_-(f)(n_k)<g(n_k)<n_{k+1}.$$ 
Construct a real $z$ such that 
$z \rest \lft1[n_k,\varphi_-(f)(n_k)\rgt1)=S^{n_k}_{f(n_k)}$. Thus $z  \in (M)_f$
but $z \rest \lft1[n,g(n)\rgt1)\not \equiv 0$ for all $n$, so  $z
\in \{x \in \twoom : \forall^\infty n \ \exists i \in \lft1[n,f(n)\rgt1) \ 
(x(i) \neq 0)\}$.
\end{proof}
\begin{theorem}
  $(\omom,\omom,\not\geq) \tles (\reals,\M,\not\ni)$;
  there are mappings $ \varphi_-: \omom \longrightarrow
  \omom $ and $ \varphi_+: \reals \longrightarrow
  \omom $ such that for $f \in \omom $ and $y \in
  \reals$,
$$ f \not \geq^\star \varphi_+(y) \hbox{ whenever } y \not\in
(M)_{\varphi_-(f)}.$$
In particular, $\COV(\M) \subseteq \D$ and $\B \subseteq
\UNIF(\M)$, $\cov(\M) \leq {\mathfrak d}$ and ${\mathfrak b} \leq \unif(\M)$.
\end{theorem}
\begin{proof}
  Identify $\reals \setmin \rationals$ with $ \omom $ and
  define 
$ \varphi_-(f)=h$ such that $$(M)_h=
\left\{z \in \omom : \forall^\infty n \ \lft1(z(n) \leq
  f(n)\rgt1)\right\}$$ and $ \varphi_+(y)=y$.
\end{proof}

\begin{theorem}
  $(\N,\N, \not \supseteq) \tles (\N,\reals, \in)$; there are Borel
  functions $ \varphi_-: \omom \longrightarrow \reals$ and $
  \varphi_+: \omom \longrightarrow \omom $ such that
  for $f, g \in \omom $,
$$(N)_f \not \supseteq (N)_{\varphi_+(g)}, \hbox{ whenever }
\varphi_-(f) \in (N)_g.$$ 
\end{theorem}
\begin{proof}
  Let $ \varphi_-:\omom \longrightarrow \reals $ be any Borel
  function such that for $f \in \omom $, $ \varphi_-(f)
  \not\in (N)_f$ (see Theorem~\ref{classic}(4)) and let $
  \varphi_+(g)=g$ for $g \in \omom $. 
Verification that both functions have the required properties is
  straightforward. 
\end{proof}

\bigskip

We conclude this section with some remarks concerning Luzin sets. 
\begin{definition}
  Given ${\mathbf A}=(A,A_-,A_+)$ and two cardinals $ \kappa \leq
  \lambda $ we call a set $X \subseteq A_-$ a {\it Luzin set} if $|X| \geq
  \lambda $ and for every $Y \subseteq X$, $|Y|=\kappa $, $Y \in
  \gb({\mathbf A})$. 
\end{definition}

When ${\mathbf A}=(\reals,\M,\in)$, $ \kappa =\mathbf\aleph_1 $
and $\lambda = 2^{\mathbf\aleph_0} $ then we get the original
Luzin set.
The set given by $(\reals,\N,\in)$, $ \kappa =\mathbf\aleph_1 $
and $\lambda = 2^{\mathbf\aleph_0} $ is usually called Sierpinski
set.

\begin{lemma}
Suppose that $X$ is a Luzin set determined by ${\mathbf A}$ and $
\kappa \leq \lambda $.
Then $\nor{{\mathbf A}} \geq \lambda $ and $\nor{{\mathbf A}^\perp}
\leq \kappa $.
\end{lemma}
\begin{proof}
  Since every set $Y \subseteq X$, $|Y|=\kappa $ belongs to
  $\gb({\mathbf A})=\gd({\mathbf A}^\perp)$, we get the second
  inequality.

For the first inequality note that if $y \in A_+$ then
$\{x \in X \cap A_-:A(x,y)\}$ has size $< \kappa \leq |X|$. Thus any
family that dominates $X$ has to have a size at least $|X|\geq \lambda 
$. 
\end{proof}
  
Morphisms preserve Luzin sets.

\begin{lemma}\label{cichonspaper}
  Suppose that ${\mathbf A} \leq {\mathbf B}$ and $X$ is a $(\kappa,
  \lambda) $ Luzin set in ${\mathbf A}$. Then $\varphi_-{``}(X)$ is a
  $\kappa, \lambda $  Luzin set in ${\mathbf B}$.
\end{lemma}
\begin{proof}
  Clearly every subset of size $\kappa $ of $\varphi_-{``}(X)$ is
  unbounded. Moreover, for every $B \in B_-$, $\varphi^{-1}(B) \cap X$ 
  has size $< \kappa $. Thus $|\varphi_-{``}(X)| \geq \lambda $.
\end{proof}
{\bf Historical remarks} 
Families  of small sets as defined here appeared in various
contexts. Rec{\l}aw  \cite{RecCh}  suggested considering
small 
sets rather than cardinal characteristics.

Many people contributed to the proof of the Theorem~\ref{8}.
In the last diagram:
\begin{itemize}
\item Rothberger \cite{Rothb38} showed that $\cov(\M) \leq \unif(\N)$
  and $\cov(\N) \leq \unif(\M)$.
\item Miller \cite{Mil81Som} and Truss \cite{Tru77Set} showed that
$\add(\M)=\min\{{\mathfrak b}, \cov(\M)\}$ and Fremlin showed that 
$\cof(\M)=\max\{{\mathfrak d}, \unif(\M)\}$.
\item Bartoszynski \cite{Bar84Add} and Raisonnier and Stern in their
  \cite{RaiSte85Str} showed that $\add(\N) \leq \add(\M)$ and
  $\cof(\M) \leq \cof(\N)$.  
Different proofs of these inequalities have been found --- forcing proof
  by Judah and Repick\'y \cite{JuReAmRe} and a very general
  combinatorial argument, Theorem
  \ref{stevo1} of this paper.
\end{itemize}
Fremlin \cite{fre:cichon} first realized that Tukey embeddings are
responsible for the 
inequalities in the Cicho\'n diagram.
Pawlikowski \cite{Paw85Leb} proved Lemma~\ref{secondlemma}, which
was the crucial step in the proof of $\M \tles \N$.

The first  diagram of Theorem~\ref{8}:
\begin{itemize}
\item Vojt\'a\v s \cite{Vojtas} proved it with arbitrary morphisms,
\item Rec{\l}aw \cite{RecCh} proved a version with Borel  morphisms
  (which gives the second diagram),
\item Pawlikowski and Rec{\l}aw in their \cite{PawRecPar} proved the
  existence of 
  continuous morphisms. 

 \end{itemize}
Lemma \ref{cichonspaper} was proved in \cite{Cich89Car}.

\section{Combinatorial characterizations}

This section is devoted to the combinatorics associated with the
cardinal invariants of the Cicho\'n diagram.
We will find the combinatorial equivalents of most of the invariants
as well as characterize membership in the corresponding classes of
small sets.
We conclude the section with a characterization 
of the ideal $(\N, \subseteq)$ as maximal in the sense of
Tukey connections among a large class of partial orderings.

\begin{theorem}\label{main1}
  The following conditions are equivalent:
  \begin{enumerate}
  \item $X \in \COV(\M)$,
  \item for every Borel function $x \rightsquigarrow f^x \in \omom$
    there exists a function $g \in 
  \omom$ such that 
$$\forall  x \in X \ \exists^\infty n \
\lft1(f^x(n)=g(n)\rgt1).$$
  \end{enumerate}
\end{theorem}
\begin{proof}
 $(1) \rightarrow (2)$. Suppose that $x \rightsquigarrow f^x \in \omom$
is a Borel mapping.
Let $H=\{\<x,h\> \in \twoom  \times \omom : \forall^\infty n
\ \lft1(h(n) \neq f^x(n)\rgt1)\}$. Clearly $H$ is a Borel set with all $(H)_x$
meager and if $g \not \in
\bigcup_{x \in X} (H)_x$ then $g$ has required properties.

\vspace{0.1in}

$(2) \rightarrow (1)$.  
We will need several lemmas. To avoid repetitions let us define:
\begin{definition}
  Suppose that $X \subseteq \twoom$. $X$ is {\it nice} if for every  Borel function
  $x \rightsquigarrow f^x \in \omom$ there exists a function $g \in
  \omom$ such that 
$$\forall  x \in X \ \exists^\infty n \
\lft1(f^x(n)=g(n)\rgt1).$$
\end{definition}

\begin{lemma}\label{lem1}
  Suppose that $X$ is nice.
Then for every Borel function 
$x \rightsquigarrow \<Y^x,f^x\> \in {}^\omega[\omega] \times \omom$
there exists $g \in  
\omom$ such that  
$$\forall x \in X \ \exists^\infty n \in Y^x \ (f^x(n)=g(n)).$$
\end{lemma}
\begin{proof}
Suppose that a Borel mapping $x \rightsquigarrow \<Y^x,f^x\>$ is given.
Let $y^x_n$ denote the
$n$-th element 
of $Y^x$ for $x \in \twoom$.
For every $x \in \twoom $ define a function $h^x$ as
follows:
$$h^x(n) = f^x\rest
\left\{y^x_0,\ y^x_1, \ldots,
y^x_n\right\} \hbox{ for } n \in \omega .$$
Since the mapping $x \rightsquigarrow h^x$ is Borel and functions $h^x$ can be
coded as elements of 
$^{\omega}\omega $ there is a function $h$ such that
$$\forall x \in X \ \exists^{\infty}n \ \lft1(h^x(n) =
h(n)\rgt1) .$$
Without loss of generality we can assume that $h(n)$ is a function
from an $n+1$-element subset of $\omega$ into $\omega$. 

Define $g \in {}^{\omega}\omega $ in the following way.
Recursively choose 
$$z_{n} \in  \dom\lft1(h(n)\rgt1)\setmin \left\{z_{0},
z_{1},\ldots, 
z_{n-1}\right\} \hbox{ for } n \in \omega.$$
Then let $g$ be any function such that
$g(z_{n}) = h(n)(z_{n})$ for $n \in \omega$.

We show that the function $g$ has the required properties.
Suppose that $x \in X$. Notice that the equality $h^x(n)= h(n)$ 
implies that
$$f^x(z_{n}) = g(z_{n}) \hbox{ and } z_{n} \in Y^x .$$
That finishes the proof since $h^x(n) = h(n)$ for infinitely
many $n \in \omega$.
\end{proof}

\begin{lemma}\label{lem1a}
  Suppose that $X$ is nice. Then for every Borel mapping $x \rightsquigarrow
  f^x \in \omom$ there exists an increasing sequence $\<n_k :
  k \in \omega\>$ such that 
$$\forall x \in X \ \exists^\infty k \ \lft1(f^x(n_k)<n_{k+1}\rgt1).$$
\end{lemma}
\begin{proof}
 Suppose that the lemma is not true and
let $x \rightsquigarrow f^x$ be a counterexample. Without loss of generality we can
assume that $f^x$ is  increasing for all $x \in X$. To get a
contradiction we will define a Borel mapping $x \rightsquigarrow g^x \in
\omom$ such that $\{g^x: x \in X\}$ is a dominating family.
That will contradict the assumption that $X$ is nice.

Define
$$g^x(n)=\max\{\underbrace{f^x\circ f^x\circ \dots \circ f^x}_{j+1 \hbox{
    times }}(i) : i,j \leq
  n\} \quad\hbox{ for } n \in \omega.$$
Suppose that $g \in \omom$ is an increasing function. By the
assumption there exist $x \in X$ and $k_0$ such that 
$$\forall k\geq k_0 \ \lft1(f^x\lft1(g(k)\rgt1) \geq g(k+1)\rgt1).$$
In particular,
$$\forall k \geq g(k_0) \ \lft1(g(k) \leq g^x(k)\rgt1)$$
which finishes the proof.
\end{proof}

%\begin{lemma}\label{lem2}
%  Suppose that $X$ is nice. Then for every Borel mapping $x \rightsquigarrow
%  Y^x \in [\omega]^\omega$ there exists a set $Y=\{u_n: n \in
%  \omega\}$ such that 
%  \begin{enumerate}
%    \item $u_{n+1} \geq u_n+2$ for all $n$,
%    \item $\forall x \in X \ |Y \cap Y^x |=\alef_0$.
%  \end{enumerate}
%\end{lemma}
%\begin{proof}
%By applying \ref{lem2}, we can find an increasing function $f \in
%\omom$ such that 
%$$\forall x \in X \ |Y^x \setmin \{f(n): n \in \omega\}| =
%\alef_0.$$
%Let $A_0 = \{2k : k \in \omega\}$ and $A_1 = \{2k+1 : k \in \omega\}$.
%Define Borel mapping $x \rightsquigarrow \<Z^x,g^x\> \in [\omega]^\omega
%\times \twoom$ as follows: $Z^x = \dom(g^x)$ and for $n \in \omega$,
%$$g^x(n)=\left\{\begin{array}{ll}
% 0 & \hbox{if } Y^x \cap \lft1(f(n),f(n+1)\rgt1) \cap A_0 \neq
%  \emptyset,\\
%1  & \hbox{if } Y^x \cap \lft1(f(n),f(n+1)\rgt1) \cap A_1 \neq
%  \emptyset\\
%& \hbox{undefined otherwise.} 
%\end{array}\right. . $$

%Note that the first two conditions of this definition are not
%exclusive. We use either value when that happens.

%By Lemma~\ref{lem2}, there exists a function $h \in \twoom$ such that 
%$$\forall x \in X \ \exists^\infty n \in Z^x \ h(n)=g^x(n).$$
%Define
%$$Y = \bigcup_{n \in \omega} \lft1((f(n),f(n+1)\rgt1) \cap A_{h(n)}.$$
%It is clear that $Y$ has required properties.
%\end{proof}
 
We now return to  the proof that $(2)$ implies $(1)$ for \ref{main1}.
Let $x \rightsquigarrow f_x \in \omom $ be a Borel mapping. We
want to 
show that $\bigcup_{x \in X} (M)_{f_x} \neq \twoom $.

Without loss of generality we can assume that $M$ is the set built
using the family from the proof of Lemma~\ref{lastlem}.
For each $x$ let $g_x \in \omom $ and $\{s^x_n: n \in
\omega\}$ be such that 
$S^n_{f_x(n)}=\{x \in \twoom : x \rest \lft1[n, g_x(n)\rgt1)=s^x_n\}$.

By Lemma~\ref{lem1a},
  there exists a sequence $\<n_k: k \in \omega\>$ such that
  \begin{enumerate}
    \item $n_{k+1} > \sum_{i=0}^k n_i$, for all $k$,
    \item $\forall x \in X \ \exists^\infty n \ \lft1(g_x(n_k)<n_{k+1}\rgt1)$.
  \end{enumerate}

For $x \in X$ let $Z^x = \{k : g_x(n_k)<n_{k+1}\}$. By Lemma~\ref{lem1},
there exists a sequence 
$\<s_{k} : k \in \omega \>$ 
such that 
$$\forall x \in X \ \exists^{\infty} k \in Z^x \ \lft1(s^{x}_{n_k} =
s_{k}\rgt1). $$
Without loss of generality we can assume that $s_k:[n_k,m_k)
\longrightarrow 2$, where $m_k<n_{k+1}$.
Choose $z \in \twoom $ such that $s_k \subseteq z$ for all $k$.
%Define mapping $x \rightsquigarrow Y^x$ by
%$$Y^x = \{k\in Z^x : s_k=s^x_{n_k}\}.$$
%Let $Y$ be a set obtained by applying Lemma~\ref{lem2} to this family.
%Define
%$$z = s_{l_0}\!^\frown s_{l_1}\!^\frown s_{l_2}\!^\frown \dots ,
%\hbox{ where } l_0<l_1<l_2 \dots \hbox{ is the increasing enumeration of
%  $Y$.}$$
%Note that if $l_{k+1} \in Y \cap Y^x$ then
%$$|s_{l_0}\!^\frown s_{l_1}\!^\frown \dots^\frown s_{l_k}| < \sum_{j
%  \leq l_k} n_{l_j+1} <n_{l_k+2} \leq n_{l_{k+1}}$$
%and $$s_{l_{k+1}}=s^x_{n_{l_{k+1}}}.$$
It follows that $z \not\in (M)_{f_x}$ for every $x \in X$.
\end{proof}

As a corollary we have:
\begin{theorem}\label{charcovm}
The following  are equivalent:
\begin{enumerate}
\item $\cov(\M) > \kappa$,
\item $\forall F \subseteq [\omom]^\kappa \ \exists g \in
  \omom \ \forall f \in F \ \exists^\infty n \ \lft1(f(n)=g(n)\rgt1).$
  \end{enumerate}
\end{theorem}

The above proof can be dualized to give:
\begin{theorem}\label{dual}
    The following conditions are equivalent:
  \begin{enumerate}
  \item $X\times X \in \UNIF(\M)$,
  \item for every Borel function $x \rightsquigarrow f^x \in \omom$
    there exists a function $g \in 
  \omom$ such that 
$$\forall  x \in X \ \forall^\infty n \
f^x(n)\neq g(n).$$
  \end{enumerate}
\end{theorem}
We only explain  why we have $X \times X$ in (1) rather than $X$.
If we analyze the proof of Theorem~\ref{main1}, we see that in order to
produce  real $z$ such that $z \not\in \bigcup_{x \in X} (M)_{f_x}$
we had to diagonalize (find an infinitely often equal real) twice.

Similar situation arises here; each element of $X$ produces two
functions, and a real that avoids a given meager set is constructed from
two such functions, each coming from a different point of $X$.

As a corollary we get:
\begin{theorem}\label{charnonm}
$\unif({\mathcal M})$ is the size of the smallest family
$F \subseteq {}^{\omega}\omega $ such that
$$\forall g \in {}^{\omega}\omega  \ \exists f \in F \ \exists^{\infty}n
\ \lft1(f(n)=g(n)\rgt1) .$$
\end{theorem}

\begin{theorem}\label{bandcov}
  $\ADD(\M)=\B \cap \COV(\M)$. In particular, 

$\add(\M)=\min\{{\mathfrak
  b}, \cov(\M)\}$.
\end{theorem}
\begin{proof} 
The inclusion $ \subseteq $ follows immediately from Theorem~\ref{8}.

\bigskip

Suppose that $X \in \B \cap \COV(\M)$. 
Let $x \rightsquigarrow f_x \in \omom $ be a Borel mapping. Since $X
\in \COV(\M)$ there is a real $z$ such that $z \not\in \bigcup_{x \in
  X} (M)_{f_x}$.
For $x \in X$ define for $ n \in \omega $,
$$g_x(n)=\min\left\{l: \forall t \in {^n 2} \ \left(\lft1[t^\frown z \rest
  [n,l)\rgt1] \subseteq \bigcup_{m>n} S^m_{f_x(m)} \right)\right\}.$$
The mapping $ x \rightsquigarrow g_x$ is also Borel. Since $X \in \B$, it follows
  that there is an increasing function $h \in \omom $ such
  that 
$$ \forall x \in X \ \forall^\infty n \ \lft1(g_x(n) \leq h(n)\rgt1).$$
Consider the set
$$G= \bigcap_n \bigcup_{m>n} \bigcup \left\{\lft2[t^\frown z \rest
\lft1[m,h(m)\rgt1)\rgt2]: t \in {^m 2}\right\}.$$
Clearly $G$ is a dense $G_\delta $ set. 
Moreover, for every $x \in X$ there is $n$ such that 
$$\bigcup_{m>n} \bigcup \left\{\lft2[t^\frown z \rest
\lft1[m,h(m)\rgt1)\rgt2]: t \in {^m 2}\right\} \subseteq \bigcup_{m>n}
S^m_{f_x(m)}.$$
It follows that
$$\bigcup_{x \in X} (M)_{f_x} \subseteq \twoom \setmin G,$$
which finishes the proof.
\end{proof}

From the Theorem~\ref{8} it follows that $\D \cup \UNIF(\M) \subseteq
\COF(\M)$. The other inclusion does not hold. We only have the
following result dual to  Theorem~\ref{bandcov}.
\begin{theorem}
  If $X \not \in \D$ and $Y \not\in \UNIF(\M)$ then $X \times Y
  \not\in \COF(\M)$. In particular, $\cof(\M)=\max\{\unif(\M), {\mathfrak
  d}\}$. 
\end{theorem}
\begin{definition}
Let $\reals_+=\{x \in \reals : x \geq 0\}$ and define
$$\ell^1=\left\{f \in {}^\omega \reals_+ : \sum_{n=1}^\infty f(n)<
  \infty\right\}.$$
For $f, g \in \ell^1$, $ f \leq^\star g$ if $f(n) \leq g(n)$ holds for
  all but finitely many $n$.
\end{definition}
\begin{theorem}\label{312}
    The following  are equivalent:
    \begin{enumerate}
    \item $X \in \ADD(\N)$,
    \item for every Borel function $x \rightsquigarrow S^x \in {\mathcal C}$
    there exists $S \in 
  {\mathcal C}$ such that 
$$\forall  x \in X \ \forall^\infty n \
\lft1(S^x(n) \subseteq  S(n)\rgt1).$$
    \item for every Borel function $x \rightsquigarrow f^x \in \ell^1$
    there exists a function $f \in 
  \ell^1$ such that 
$$\forall  x \in X \ \forall^\infty n \
\lft1(f^x(n)\leq f(n)\rgt1).$$
    \end{enumerate}
In particular,   the following conditions are equivalent:
  \begin{enumerate}
  \item[a.] $\add(\N) > \kappa $,
  \item[b.] for every family $F \subseteq \omom $ of size $
    \kappa 
    $ there exists $S \in {\mathcal C}$ such that 
$$ \forall f \in F \ \forall^\infty n \ \lft1(f(n) \in S(n)\rgt1),$$
\item[c.] for every family $F \subseteq \ell^1 $ of size $ \kappa
    $ there exists $g \in\ell^1$ such that 
$$ \forall f \in F \ \forall^\infty n \ \lft1(f(n) \leq g(n)\rgt1).$$
  \end{enumerate}
\end{theorem}
\begin{proof}
We will establish the equivalence of (1) and (2). 
Suppose that $X \in \ADD(\N)$ and $x \rightsquigarrow S_x$ is a Borel
mapping. Consider the morphism $( \varphi_-, \varphi_+)$ witnessing
that $ {\mathcal C} \tles \N$. 
Let $f $ be such that $ \bigcup_{x \in X} (N)_{\varphi_-(S_x)}
\subseteq (N)_f$. Then $ \varphi_+(f) \in {\mathcal C}$ is the object
we are looking for. 

Suppose that $X \not \in \ADD(\N)$. Let $F : X \longrightarrow
\omom $ be a Borel mapping such that $ \bigcup_{x \in X}
(N)_{F(x)} \not \subseteq (N)_f$ for $f\in \omom $.
Consider the morphism $( \varphi_-, \varphi_+)$ witnessing
that $ \N \tles {\mathcal C}$.
It follows that there is no $S \in {\mathcal C}$ such that 
$$ \forall x \in X \ \forall^\infty n \ 
\lft1(\varphi_-\lft1(F(x)\rgt1)(n) \subseteq^\star S(n)\rgt1).$$

Equivalence of (2) and (3) follows from:
\begin{lemma}
  ${\mathcal C} \equiv \ell^1$.
\end{lemma}
\begin{proof}
To show that $\ell^1 \tles {\mathcal C}$ 
  define $ \varphi_-: \ell^1 \longrightarrow {\mathcal C}$ as
$$ \varphi_-(f)(n)=\left\{k: 2^{-n}>f(k)\geq 2^{-n-1}\right\}.$$
Similarly, define $ \varphi_+: {\mathcal C} \longrightarrow \ell^1$ is defined
  by:
$ \varphi_+(S)(n)=\max\{2^{-k}: n \in S(k)\}.$
It is easy to see that these mappings have the required properties.

\bigskip

To show that ${\mathcal C} \tles \ell^1$ identify $ \omega \times \omega
$ with $ \omega $ via functions $L,K \in \omom $.
For $S \in {\mathcal C}$ let
$$ \varphi_-(S)(n)=\left\{
  \begin{array}{ll}
2^{-n} & \hbox{if } K(n)\in S\lft1(L(n)\rgt1)\\
0 & \hbox{otherwise}
  \end{array} \right. .$$
For $f \in \ell^1$ let
$$ \varphi_+(f)(n)=\left\{k: \frac{1}{2^{n-1}}> f(k) \geq
    \frac{1}{2^n}\right\}.$$ 
\end{proof}
The second part of \ref{312}  follows readily from the first.
\end{proof}

The dual version yields:
\begin{theorem}\label{314}
    The following  are equivalent:
    \begin{enumerate}
    \item $X \in \COF(\N)$,
    \item for every Borel function $x \rightsquigarrow S^x \in {\mathcal C}$
    there exists $S \in 
  {\mathcal C}$ such that 
$$\forall  x \in X \ \exists^\infty n \
\lft1(S(n) \not\subseteq  S^x(n)\rgt1).$$
    \item for every Borel function $x \rightsquigarrow f^x \in \ell^1$
    there exists a function $f \in 
  \ell^1$ such that 
$$\forall  x \in X \ \exists^\infty n \
\lft1(f^x(n)\leq f(n)\rgt1).$$
    \end{enumerate}
In particular, the following  are equivalent:
  \begin{enumerate}
  \item[a.] $\cof(\N) < \kappa $,
  \item[b.] for every family $F \subseteq \omom $ of size $ \kappa
    $ there exists $S \in {\mathcal C}$ such that 
$$ \forall f \in F \ \exists^\infty n \ \lft1(f(n) \not\in S(n)\rgt1),$$
\item[c.] for every family $F \subseteq \ell^1 $ of size $ \kappa
    $ there exists $g \in\ell^1$ such that 
$$ \forall f \in F \ \exists^\infty n \ \lft1(f(n) \leq g(n)\rgt1).$$
  \end{enumerate}
\end{theorem}

Additivity of measure, $\add(\N)$, has a special place among cardinal
invariants of the continuum as being provably smaller than a large
number
 of them. It has been conjectured (wrongly in \cite{BagJud}) that
this is because additivity of measure 
is equivalent to the Martin Axiom for a
large class of forcing notions (Suslin ccc). Only very recently
this phenomenon has been explained as 
being directly related to the combinatorial complexity of the measure
ideal. 
 
 \begin{definition}
We say that an ideal $ {\mathcal J}  \subseteq \Power(\omega)$ is a
{\it p-ideal} if for every family $\{X_n: n \in \omega\} \subseteq {\mathcal
  J} $ there is $X \in {\mathcal J} $ such that $X_n \subseteq^\star
X$ for $ n \in \omega $.

Define 
$$\add^\star({\mathcal J})=\min\{|{\mathcal A}| : {\mathcal A}
\subseteq {\mathcal J}\ \& \ \neg \exists Y \in {\mathcal J} \ \forall
X \in {\mathcal A} \ X \subseteq^\star Y\}.$$
\end{definition}

It is easy to see that the $\cof^\star$ defined analogously is
equivalent to the true cofinality.

Many ideals of
Borel subsets of $\reals $ are Tukey equivalent to analytic (${\mathbf
  \Sigma}^1_1$) ideals of 
subsets of $ \omega $.
For example:
\begin{itemize}
\item Ideal of null sets $\N$. 
Work with ${\mathcal C}$ instead of $\N$. For $S \in {\mathcal C}$ let
$A_S=\{(n,k) \in \omega\times \omega : k \in S(n)\} $. The family
$\{A_S: S \in {\mathcal C}\}$ generates an analytic p-ideal on $
\omega \times \omega \simeq \omega$.
\item $(\omom, \leq^\star)$.
This is very similar. For $f \in \omom $ let
$A_f = \{(n,k): k \leq f(n)\}$.
\item The ideal of meager sets $\M$.
Let $\{C_n: n \in \omega\}$ be an enumeration of open basic subsets of 
$\omega $. Consider the ideal generated by sets $X \subseteq \omega $
such that for every $k$, $\bigcup_{n \not\in X\cup k} C_n$ is open dense.
This is an analytic p-ideal on $\omega $ which is equvalent to $\M$.
\item $\ell^1$ is Tukey equivalent to the p-ideal of summable sets
$$\left\{X \subseteq \omega : \sum_{n \in X}
  \frac{1}{n}<\infty\right\}.$$    
\end{itemize}

Moreover, in all these cases the additivity of the ideal is equal to
the $^\star$additivity of the associated ideal on $ \omega $.
For example,  $\add(\N)=\add^\star({\ell^1})=\add^\star({\mathcal
  C})$, etc. 
In the remainder of this section we will show 
for a (nontrivial) analytic p-ideal $ {\mathcal J} $ on $ \omega $ we
have
$$ \omom \tles {\mathcal J} \tles \N.$$
In particular, by the above remarks, 
 $ {\mathfrak b} \geq \add(\M) \geq
\add(\N)$. 

We need a few general facts
about analytic p-ideals.
  To simplify the notation let us identify $ \twoom $ with
  $\Power(\omega)$ via characteristic functions. 

Let $\K(\twoom)$ be the collection of compact subsets of $\twoom $
with Hausdorff metric $d_H$ defined as follows.
For two nonempty compact sets $K,L \subseteq \twoom $ let
$$d_H(K,L)=\max\lft1(\rho(K,L),\rho(L,K)\rgt1),$$
where $\rho(K,L)=\max_{x \in L} d(x, K)$ ($d$ is the usual metric in
$\twoom $).

Let $\Mon \subseteq \K(\twoom )$ be the collection of compact subsets
of $\twoom $ which are downward closed.
We will use the following well known facts:
\begin{lemma}
  \begin{enumerate}
  \item $\K(\twoom)$ is a compact Polish space,
  \item $\Mon$ is a closed subspace of $\K(\twoom)$.
  \end{enumerate}
\end{lemma}
\begin{proof}
  See \cite{Kechris95} 4.F.
\end{proof}
Let 
$$\FF =\{K \in \Mon: \forall X \in {\mathcal J} \ \exists n \ (X
\setmin n \in K)\}.$$
It is clear that $\FF$ is a filter.

\begin{lemma}\label{onestep}
  Suppose that $H \subseteq \FF$ is a closed set. There exists a
  relatively clopen set $U \subseteq H$ such that 
$$ \bigcap_{K \in U} K \in \FF.$$
In particular, $H=\bigcup_{n \in \omega}H_n$, where for each $n$,
$ \bigcap_{K \in H_n} K \in \FF.$
\end{lemma}
\begin{proof}
Let $\<U_n: n \in \omega\>$ be an enumeration of clopen subsets of
$\K(\twoom)$.

  For $X \in {\mathcal J} $ and $n$ define
$H_n(X)=\{K\in H: X \setmin n \in K\}$.
The sets $H_n(X)$ are closed and $H = \bigcup_{n \in \omega} H_n(X)$ for
every $X \in {\mathcal J} $.
By the Baire Category theorem for each $X$ there is a pair
$\lft1(n(X),m(X)\rgt1)\in \omega \times \omega$ such that 
$$H_{n(X)}(X) \cap U_{m(X)}=H \cap U_{m(X)}.$$
Since $ {\mathcal J} $ is a p-ideal, we can find $(n,m)$ such that 
$$\{X: n(X)=n \ \&\  m(X)=m\} \hbox{ is cofinal in } {\mathcal J} .$$
It follows that $ \bigcap_{K \in U_m \cap H} K \in \FF$.
That finishes the proof of the first part.

\bigskip

To prove the second part let $H_1=H \cap U_m$. Next, apply the above
construction to $H \setmin H_1$ (which is closed) to get $H_2$, etc.
\end{proof}
\begin{lemma}\label{fsig}
  $\FF$ is $F_\sigma$ in $\Mon$.
\end{lemma}
\begin{proof}
  Consider $\G=\Mon \setmin \FF$.
Note that for $K \in \Mon$ we have
$$K \in \G \iff \exists X \in {\mathcal J} \ \forall n \ (X \setmin n
\not \in K).$$
It follows that $\G$ is an analytic ideal.
Moreover, $\G$ is a $ \sigma $-ideal; if $\{K_n: n \in \omega\}
\subseteq \G$ and $K \subseteq \bigcup_n K_n$ then $K \in \G$. 
To see this let $X_n $ witness that $K_n \in \G$. Find $X \in
{\mathcal J} $ such that $X_n \subseteq^\star X$ for all $n$. 
Clearly, $X \setmin n \not \in K$ for all $n$.
Now the lemma follows immediately from the following:
\begin{theorem}\label{chrstr}
  Let $ {\mathcal I}$ be an analytic $ \sigma $-ideal of compact sets
  in a compact metrizable space $E$. Then ${\mathcal I}$ is actually
  $G_\delta $.
\end{theorem}
\begin{proof}
  See \cite{Chr}, \cite{str1} or \cite{kelouwoo:sigmai}.
\end{proof}
\end{proof}
\begin{lemma}\label{coungen}
  $\FF$ is countably generated.
\end{lemma}
\begin{proof}
Using Lemma~\ref{fsig}
  represent $\FF=\bigcup_n H_n$, where each $H_n$ is closed. 
Apply Lemma~\ref{onestep} to write for $ n \in \omega $,
 $H_n=\bigcup_{m \in \omega } H^n_m$, where $G^n_m=\bigcap_{K \in
  H^n_m} K 
  \in \FF$. It is clear that $\{G^n_m: n, m \in \omega \}$ generates
  $\FF$.  
\end{proof}

Let $\<G_n : n \in \omega\>$ be a descending 
sequence generating $\FF$. The following lemma gives a simple
($F_{\sigma \delta }$)
description of $ {\mathcal J} $ in terms of $\<G_n: n \in
\omega\>$. 
\begin{lemma}\label{char}
  $X \in {\mathcal J} \iff \forall n \ \exists m \ (X \setmin m \in
  G_n).$ 
\end{lemma}
\begin{proof}
  Implication $( \rightarrow )$ is obvious.

\bigskip

$( \leftarrow )$ We will use the following result.
\begin{theorem}\label{meafil}
  Suppose that $ {\mathcal I} \subseteq \Power(\omega)$ is an ideal
  containing all finite sets.
The following conditions are equivalent:
\begin{enumerate}
\item ${\mathcal I}$ has the Baire property,
\item ${\mathcal I}$ is meager,
\item there exists a partition $\{I_n: n \in \omega\}$ of $ \omega $
  into disjoint intervals such that 
$$ \forall X \in {\mathcal I} \ \forall^\infty n \ (I_n \not \subseteq
X).$$ 
\end{enumerate}
\end{theorem}
\begin{proof}
See  \cite{Tal80Com} or  \cite{BJbook}.
\end{proof}
Suppose that $X \not \in {\mathcal J} $. 
The ideal $ {\mathcal J} \rest X = \{Y \cap X: Y \in {\mathcal
  J} \} \subseteq \Power(X)$ is  analytic and hence has the Baire
property. 
By  Theorem~\ref{meafil}(3) there exists a partition $\{I_n: n \in \omega\}$ of $X$
  into finite sets such that 
$$ \forall Z \in {\mathcal J} \ \forall^\infty n \ (I_n \not \subseteq
Z)  .$$
Choose $n_0 \in \omega $ such that 
the set
$$K=\{Y: \forall n>n_0 \ (I_n \not \subseteq Y)\}\in \FF.$$
 Let $n$ be such that $G_n \subseteq K$. It follows that for every $m
 \in \omega $,
$$ X \setmin n =^\star \bigcup_{n \in \omega} I_n  \not\in G_n,$$
which finishes the proof of \ref{char}.
\end{proof}

For $K, L \in \K(\twoom)$ define 
$K \oplus L =\{X \cup Y: X \in K, \ Y \in L\}$.
($\cup$ is in $\Power(\omega)$  the same as coordinate-wise 
maximum in $\twoom$)

Let $\<G_n : n \in \omega\>$ continue to be a descending 
sequence generating $\FF$.
\begin{lemma}\label{sol}
  For every $K \in \FF$ there exists $m$ such that $G_m \oplus G_m
  \subseteq K$. 
\end{lemma}
\begin{proof}
  Fix $X \in {\mathcal J} $ and using the fact that $ {\mathcal J} $
  is a p-ideal find $k$ such that $\{Y \in {\mathcal J} : X \setmin
  m \subseteq Y\}$ is cofinal in $ {\mathcal J} $.
The set
$$H_X=\{Y: (X \setmin k) \cup Y \in K\} \in \FF.$$
Let $n(X)$ be such that $G_{n(X)} \subseteq H_X$.
We have
$$\{X \setmin n(X)\} \oplus G_{n(X)} \subseteq K.$$
Choose $n$ such that $\{X: n(X)=n\}$ is cofinal in $ {\mathcal J} $.
The set $L = \{X: \{X \setmin n\} \oplus G_n \subseteq K\} \in \FF.$
Let $m \geq n$ be such that $ G_m \subseteq L$.
It follows that $G_m \oplus G_m \subseteq K$.
\end{proof}

We are ready to formulate the first result.

\begin{theorem}\label{stevo1}
Suppose that $ {\mathcal J} $ is an analytic p-ideal on $ \omega
$. Then
$ {\mathcal J} \tles \ell^1.$  In particular, $\add^\star({\mathcal
  J})\geq 
\add(\N)$ and $\cof({\mathcal J}) \leq \cof(\N)$.
\end{theorem}
\begin{proof}
Use Lemma~\ref{sol} to find a descending sequence
$\<G_n : n \in \omega\>$ 
generating $\FF$ such that for each $n$,
$$\underbrace{G_{n+1} \oplus \cdots \oplus G_{n+1}}_{2^{n+1} \hbox{
    times}} \subseteq G_n.$$
For $X \in {\mathcal J} $ let $\<k_n(X): n \in \omega \>$ be an
    increasing sequence such that 
$$ \forall n \ (X \setmin k_n(X) \in G_{n+2}).$$

Identify $ \omega $ with $[\omega]^{<\omega}$ and define
$ \varphi_-: {\mathcal J}
\longrightarrow {\mathcal C}$ and $ \varphi_+:
{\mathcal C} \longrightarrow {\mathcal J} $ such that 
$$X \subseteq^\star \varphi_+(S), \hbox{ whenever } \varphi_-(X)
\subseteq^\star S.$$
Since ${\mathcal C} \equiv \ell^1 \equiv \N$ this will finish the
proof. 
For $X \in {\mathcal J} $ and $ n \in \omega $ define
$$ \varphi_-(X)(n)=X \cap k_n(X) \in 
[\omega]^{<\omega} \simeq \omega .$$
Mapping $ \varphi_+$ will be defined as follows.
Suppose that $S \in {\mathcal C}$ is given (with $S(n) \subseteq
[\omega]^{<\omega}$).  For $ n \in \omega $ let
$$Z_n=\left\{(t,s)\in S(n+1)\times S(n): s \subseteq t \ \&\ t
  \setmin \max(s) \in G_{n+2}\right\}.$$
Now define 
$$v_n=\bigcup_{(t,s) \in Z_n} t \setmin \max(s).$$
Note that
$v_n$ is a sum of at most $2^{n+1}$ terms, each belonging to
$G_{n+2}$. Thus, $v_n \in G_{n+1}$ for all $n$.

The motivation for this definition is following: if 
$\varphi_-(X)(n)=X \cap k_n(X) \in S(n)$ and
$\varphi_-(X)(n+1)=X \cap k_{n+1}(X) \in S(n+1)$,  then
$$X \cap k_{n+1}(X) \setmin \max\lft1(X \cap k_n(X)\rgt1) = X \cap \lft1[k_n(X),
k_{n+1}(X)\rgt1) \subseteq v_n.$$
  The requirements of the
definition describe this situation and filter out ``background noise''
coming with $S$.

Finally define 
$$ \varphi_+(S) = Y=\bigcup_n v_n.$$
By the remarks above it is clear that if $X \in {\mathcal J} $ and $S
\in {\mathcal C}$ then from the fact that
$$ \forall^\infty n \ \lft1(\varphi_-(X)(n) \in S(n)\rgt1)$$
it follows that $X \subseteq^\star Y=\varphi_+(S)$.
To finish the proof it remains to show that the range of $ \varphi^+$
is contained in $ {\mathcal J} $.

Let $ \varphi_+(S)=Y = \bigcup_n v_n$ be defined as above.
For $j \in \omega $, let $Y_j = \bigcup_{n \geq j} v_n$.
Since $Y \setmin Y_j$ is finite for every $j$, 
by the lemma above, in order to show that $Y \in {\mathcal J} $ it
would suffice to  show that $Y_j \in G_j$.
\begin{lemma}
  For each $l\in \omega $,
$$v_n \cup v_{n+1} \cup \dots \cup v_{n+l} \in G_n.$$
\end{lemma}
\begin{proof}
We prove this by induction on $l$. For each $n$, $v_n \in G_{n+1}$ so
the lemma is true for $l=0$.  Suppose it holds for $l$ and all $n$.
We have 
$$v_n \cup v_{n+1} \cup v_{n+l+1} = v_n \cup (v_{n+1} \cup \dots v_{n+1
  + l}) \in G_{n+1} \oplus G_{n+1} \subseteq G_n,$$
which finishes the proof.  
\end{proof}
Since sets $G_n$ are closed, we
  conclude that $Y_j \in G_j$. In particular, by Lemma~\ref{char}, $Y=Y_0
  \in {\mathcal J}$.
\end{proof}

The last theorem gave us an lower bound for $\add^\star({\mathcal
  J})$. The next theorem  gives us an upper bound.

Suppose that $ {\mathcal J} \subseteq \Power(\omega)$ is an ideal.
We say that $ {\mathcal J} $ is {\it atomic} if there is $Z \in {\mathcal J}
$ such that $ {\mathcal J} = \{X \subseteq \omega : X \subseteq^\star
Z\}$. It is clear that $\add^\star({\mathcal J})$ is undefined (or
equal to $\infty$) for an
atomic ideal. 

\begin{theorem}\label{stevo2}
  Suppose that $ {\mathcal J} $ is an analytic p-ideal which is not
  atomic. 
Then $\omom \tles {\mathcal J}$. In particular,
$ \add^\star({\mathcal J}) \leq {\mathfrak b}$ and $\cof({\mathcal J})
  \geq {\mathfrak d}$.
\end{theorem} 
\begin{proof}
For $X \subseteq \omega $ let $\overline{X} \in \prod_{n \in \omega} n$ be
defined as $\overline{X}(n)=|X \cap n|$.
Let $ {\mathcal X} = \left\{\overline{X}: X \in \Power(\omega)\right\}$. 
It is easy to see that $ {\mathcal X} $ is a compact space.
For $\overline{X}, \overline{Y} \in {\mathcal X} $ define
$ \overline{X} \leq^\star \overline{Y}$ if there is a finite set $Z \subseteq
\omega $ such that 
$$ \forall n \ \overline{X}(n) \leq \overline{Y \cup Z}(n)$$
or alternatively
 $$ \forall n \ \overline{X \setmin Z}(n) \leq \overline{Y}(n).$$

  \begin{lemma}
    $( {\mathcal J}, {\mathcal J}, \not \supseteq^\star) \tles
({\mathcal X}, {\mathcal X}, \not \geq^\star)$. 
  \end{lemma}
  \begin{proof}
    Define $ \varphi_-: {\mathcal J} \longrightarrow \omom $ and 
$ \varphi_+: \omom \longrightarrow {\mathcal J} $ as
$ \varphi_-(X)=\overline{X}$ and $ \varphi_+(\overline{X})=X$.
Suppose that $ \varphi_-(X) \not \geq^\star \overline{Y}$. That means that
for any finite set $Z$, 
$$ \exists^\infty n \ \lft1(\overline{X \cup Z}(n) < \overline{Y}(n)\rgt1).$$
It follows that $Y \not \subseteq X \cup Z$, hence $ Y \not
\subseteq^\star X$.
  \end{proof}

Let $\Mon \subseteq \K({\mathcal X})$ be the collection of compact
subsets 
of $ {\mathcal X}  $ which are downward closed (with respect to
$\leq$).
Let 
$$\FF=\left\{K \in \Mon: \forall X \in {\mathcal J} \ \exists n
\ \lft1(\overline{X \setmin n} \in K\rgt1)\right\}.$$
As before, by Lemma~\ref{coungen}, $\FF$ is countably generated. 
Let $\<G_n: n \in \omega\>$ be any sequence generating $\FF$.

We will show that
$( {\mathcal X}, {\mathcal X}, \not \geq^\star) \tles
(\omom, \omom, \not \geq^\star)$. We need functions
$ \varphi_-: {\mathcal X} \longrightarrow \omom $ and 
$ \varphi_+: \omom \longrightarrow {\mathcal X} $ such that
for $\overline{X} \in {\mathcal X} $ and $f \in \omom $,
$$\overline{X} \not \leq^\star \varphi_+(f), \hbox{ whenever }
\varphi_-(\overline{X}) 
\not \geq^\star f.$$
Clearly, the dual morphism witnesses
 that  $\omom \tles {\mathcal
  X}\tles {\mathcal J}$.

Each set $G_n$ is a set of branches of some tree. By taking the
rightmost branch (towards the larger values) at every node of the
tree,  we produce a countable
family $\left\{\overline{Z}^n_m: m \in \omega \right\}$ such that 
$$ \forall X \in {\mathcal J} \ \exists m \ \lft1(\overline{X} \leq^\star
\overline{Z}^n_m\rgt1).$$
Since $ {\mathcal J} $ is a p-ideal, for each $n$ there is $m$ such
that 
$$\left\{X\in {\mathcal J} : \overline{X} \leq^\star \overline{Z}^n_m\right\} \hbox{ is
  cofinal in } {\mathcal J} .$$
Denote the first such $\overline{Z}^n_m$ by $\overline{U}_n$.
Without loss of generality we can assume $\overline{U}_n \geq
  \overline{U}_{n+1}$ for all $n$.

We have the following two cases:\\
{\sc Case 1}. \hspace{0.1in} There exists $\overline{Z} \in {\mathcal X} $ such that 
\begin{enumerate}
\item $\overline{X} \leq^\star \overline{Z}$ for every $X \in {\mathcal J} $.
\item $\overline{Z} \leq^\star \overline{U}_n$ for every $n$.
\end{enumerate}
In this case $ {\mathcal J} $ is atomic. 
Note that condition (2) implies that
$$ \forall n \ \exists m \ \lft1(\overline{Z \setmin m} \leq \overline{U}_n
\in G_n\rgt1).$$ 
Thus by Lemma~\ref{char}, $Z \in {\mathcal J}$.
Condition (1), together with the fact that $ {\mathcal J} $ is an
ideal, implies that $X \subseteq^\star Z$ for every $X \in {\mathcal
  J}$.

\bigskip

{\sc Case 2}.  \hspace{0.1in} Suppose that there is no $\overline{Z}$ as in {\sc Case 1}

For $X \in {\mathcal J} $ and $ n \in \omega $ define
$$ \varphi_-(\overline{X})(n) = \max\left\{\max\left\{j: \overline{X}(j) >
  \overline{U}_i(j)\}: i \leq n\right\}\right\}.$$ 
For $f \in \omom $ let 
$$ \varphi_+(f) \hbox{ be any } X \in {\mathcal
  J} \hbox{ such that } f \leq^\star \varphi_-(\overline{X}).$$
It is clear that these mappings have the required properties provided
that they are correctly defined.
Thus, the following lemma will complete the proof:
\begin{lemma}
  For every $f \in \omom $ there exists an $X \in {\mathcal J} $
  such that $ f \leq^\star \varphi_-(\overline{X})$. 
\end{lemma}
\begin{proof}
  Suppose not and let $f \in \omom $ be a strictly increasing
  function such that for every set $X \in {\mathcal J} $ the set 
$$Z_X=\left\{n: \varphi_-(\overline{X})(n) > f(n)\right\}$$
is coinfinite. 
It follows that the family $\{Z_X: X \in {\mathcal J} \}$ generates a
proper analytic ideal $ {\mathcal H}$.
As before, $ {\mathcal H} $ has the Baire property, hence by
Theorem 
\ref{meafil}(3), there exists a partition $\{I_n: n \in \omega\}$ such
that 
$$ \forall X \in {\mathcal J} \ \forall^\infty n \ \lft1(I_n \not \subseteq
Z_X\rgt1).$$ 
Let $h(n)=f\lft1(\max(I_{n+1})\rgt1)$ for $ n \in \omega $ and
consider the function  
$$\overline{U} = \bigcup_n \overline{U}_n \rest
\lft1[h(n), h(n+1)\rgt1).$$  
Clearly, $\overline{U}_n \geq^\star \overline{U}$ for all $n$. 
We will show that $\overline{X} \leq^\star \overline{U}$ for  $X \in
{\mathcal J} $,  which will give the
contradiction.

Fix $X \in {\mathcal J} $. 
Suppose that $n \in I_m$ and $k \in I_{m+1} \setmin Z_X$. Clearly, 
$m \leq n$ and $k \leq \max\lft1(I_{m+1}\rgt1)$.
We have
$$\max\left\{j: \overline{X}(j) > \overline{U}_n(j)\right\} \leq \varphi_-(X)(k) \leq f(k)
\leq 
f\lft1(\max(I_{m+1})\rgt1) \leq h(n).$$
In particular,
$$ \forall j > h(n) \overline{X}(j) \leq \overline{U}_n(j).$$
It follows that 
$\overline{X} \leq^\star \overline{U}$.
\end{proof}
\end{proof}

\bigskip

{\bf Historical remarks}
Theorem~\ref{main1} was proved in \cite{PawRecPar} and
\cite{BarJud93Cov}. 
Theorem~\ref{dual} is due to Pawlikowski and Rec{\l}aw in their
\cite{PawRecPar}. Theorems \ref{charcovm} and 
 \ref{charnonm} were proved in
\cite{Bar87Com}. Theorem~\ref{bandcov} was proved in \cite{PawRecPar}.
The second part is due to Miller \cite{Mil81Som}. The first part of
Theorem~\ref{312} was proved in \cite{PawRecPar} and the second in
\cite{Bar84Add}. 
Todorcevic proved  Theorem
\ref{stevo2} 
\cite{Tod96Gaps}.
I learned Theorem~\ref{stevo1} from Todorcevic
\cite{TodTuk96}. The result can  also be attributed to
Louveau and Velickovic (see their \cite{LoVe}, Theorem 5).
 Methods used in the proof, in particular Lemmas \ref{fsig} and
 \ref{sol} are 
 due to Solecki  
\cite{Solecki} and  \cite{Solecki1}.
Similar ideas were already present in \cite{Tod96Gaps} and earlier in 
\cite{JustIF}.  
Theorem \ref{chrstr} is due to Christensen and Saint Raymond. It was
generalized in \cite{kelouwoo:sigmai}.
Theorem~\ref{meafil} was proved by Talagrand.

\section{Cofinality of $\cov({\mathcal J})$ and $\COV({\mathcal J})$}
It is clear that cardinal invariants $\add$, $\unif$ and $\cof$ have
uncountable cofinality and families $\ADD$, $\UNIF$ and $\COF$ are 
$ \sigma $-ideals.
It this section we investigate $\cov$ and $\COV$ for both ideals $\M$
and $\N$.

\begin{theorem}\label{mill}
  $\COV(\M)$ is a $ \sigma $-ideal. In particular, $
  \cf\lft1(\cov(\M)\rgt1) >
  \alef_0$.
\end{theorem}
\begin{proof}
  Suppose that $\{X_n: n \in \omega\} \subseteq \COV(\M)$. 
Let $x \rightsquigarrow f_x \in \omom $ be a Borel mapping.
It is enough to  find $g \in \omom $ such that 
$$ \forall n \ \forall x \in X_n \ \exists^\infty m \ \lft1(g(m)=f_x(m)\rgt1).$$
Let $\{A_k: k\in \omega\}$ be a partition of $ \omega $ into
infinitely many infinite pieces.
For each $n$ consider the mapping $x \rightsquigarrow f_x \rest A_n$
and find
$g_n \in {}^{A_n}\omega$ such that 
$$ \forall  x \in X_n \ \exists^\infty k\in A_n\ \lft1(f_x(k)=g_n(k)\rgt1).$$
Then $g =\bigcup_n g_n$ is as required.
\end{proof}
In the presence of many dominating reals we have a similar result for
the measure ideal.
\begin{theorem}\label{covandb}
  If $\cov(\N) \leq {\mathfrak b}$ then $
  \cf\lft1(\cov(\N)\rgt1) >\alef_0$. 
\end{theorem}
\begin{proof}
  See \cite{Bar88Cov} of \cite{BJbook}.
\end{proof}
The following surprising result of Shelah shows that 
without any additional assumptions it is not possible to show that 
$\cov(\N)$ has uncountable cofinality.

\begin{theorem}\label{sh592}
  It is consistent with $\ZFCa$ that $\COV(\N)$ is not a $ \sigma
  $-ideal and $ \cf\lft1(\cov(\N)\rgt1)=\alef_0 $.
\end{theorem}
The proof of this theorem will occupy the rest of this section.
  The model will be obtained by a two-step finite support iteration. 
We start with a suitably chosen model $\eL$ satisfying $2^{\alef_0}=\alef_1$ and 
add $\alef_\omega $ Cohen reals followed by a finite support iteration 
of subalgebras of the random algebra $ \mathbf B $. 
We start by developing various tools needed for the construction.

\bigskip 

{\bf He random real algebra.} Recall that the
  random real algebra can be represented as
$$ {\mathbf B} =\{P \subseteq \twoom : \mu(P)>0 \hbox{ and } P
\hbox{ is closed}\}.$$
For $P_1, P_2 \in \mathbf B $, $P_1 \geq P_2$ if $P_1 \subseteq P_2$.
Elements of $ \mathbf B $ can be coded by reals in the following way.
Let $\widetilde{P} \in \eL $ be a universal closed set,
i.e. $\widetilde{P} \subseteq \twoom  \times \twoom $ is Borel 
and for every
closed set $P \subseteq \twoom $ there is $x$ such that
$P=(\widetilde{P})_x$. 
Let $H = \left\{x: \mu\lft1((\widetilde{P})_x\rgt1)>0\right\}$. By
Theorem 
\ref{classic}(1), $H$ is a Borel set.
Define $\widetilde{B}=(H \times \twoom) \cap \widetilde{P}$.
If $M$ is a model of $\ZFCa$ then we define
$$ {\mathbf B}^M=\left\{P \in {\mathbf B}  : \exists x \in M \cap
\twoom \ \lft1(P = (\widetilde{B})_x\rgt1)\right\}.$$

\bigskip

{\bf $\Delta$-systems.}
The following concepts will be crucial for the construction of the model.
\begin{definition}
Let ${\mathcal R}\in \eL$ be a forcing notion
Suppose that $\bar{p}=\<p_n: n \in \omega\>$ is a sequence of
  conditions in $ {\mathcal R}$. Let $
  \dot{X}_{\bar{p}}$ be the $ {\mathcal R}$-name
  for the set
$\{n : p_n \in \dot{G}_{{\mathcal R}}\}$.
In other words, for every $n$,  $p_n=\bv{ n \in \dot{X}_{\bar{p}}}$.
\end{definition}
At the moment we will be concerned with the case when ${\mathcal R}={\mathbf
  C}_{\alef_{\omega+1}}$ is the forcing notion adding
$\alef_{\omega+1}$ Cohen reals.
\begin{definition}\label{delta}
Let $\Delta \subseteq \left[{\mathbf
  C}_{\alef_{\omega+1}}\right]^\omega $ be the collection of all
sequences $\bar{p}=\{p_n: n \in \omega \}$ such that 
 there exists $k,l \in
\omega $ and $g \in {}^{l\times \omega }\omega $, $s \in {}^k \omega $
such that 
\begin{enumerate}
\item $\dom(p_n)=\{\beta_1, \dots, \beta_k\} \dot{\cup} \{\alpha^n_1,
  \dots , \alpha^n_l\}$, with $\beta_1<\cdots<\beta_k$ and
  $\alpha^n_1< \cdots < \alpha^n_l$ for $n \in \omega$ (so the $\dom(p_n)$'s form a
  $\Delta$-system with root $\{\beta_1,\dots,\beta_k\}$),
\item $\alpha^n_i<\alpha^m_i$ for $n<m$, 
\item $p_n(\alpha^n_i)=g(i,n)$ for every $i \leq l, n \in \omega $,
\item $p_n(\beta_i)=s(i)$ for $i \leq k, n \in \omega $.
\end{enumerate}
Let $p_{\bar{p}}=p_0\rest \{\beta_1,\dots, \beta_k\}$.
\end{definition}
 
Note that  if $\bar{p} \in \Delta$ then $f_{\bar{p}} = \bigcup_{n \in
\omega} p_n$ is a function. Moreover, $p_{\bar{p}}= f_{\bar{p}} \rest
\{\beta_1, \dots, \beta_k\}$ and
$p_{\bar{p}} \forces_{{\mathbf
  C}_{\alef_{\omega+1}}} X_{\bar{p}} \text{ is infinite}$.

\begin{definition}
  A subset $\Delta' \subseteq \Delta$ is filter-like if for any 
$\bar{p}^1, \dots , \bar{p}^n \in \Delta'$ there exists $q $ such
that 
$$q \forces_{{\mathbf
  C}_{\alef_{\omega+1}}} \bigcap_{i \leq n} X_{\bar{p}^i} \text{ is
infinite}.$$ 
\end{definition}

 \begin{theorem}\label{delta1}
   Suppose that $\V \thinks 2^{\alef_0}=\alef_1 \ \&\
   2^{\alef_1}=\alef_{\omega+1}$. 
Then
   $\Delta$  is the union of $\alef_1$ filter-like sets.
 \end{theorem}
 \begin{proof}
   Let $T$ be the collection of $\<k,l,v,\{f_{i,n},g_j: i\leq
   l,j\leq k, n \in
   \omega \},g,s\>$ such that 
   \begin{enumerate}
   \item $ k,l \in \omega $,
   \item $v \in [\alef_1]^{\leq \alef_0}$,
   \item $g_j,f_{i,n} \in {}^v \omega$ are pairwise different for
     $i\leq l,j\leq k, n \in \omega $,
   \item $g \in {}^{l \times \omega }\omega $,
   \item $s \in {}^k \omega $.
   \end{enumerate}

From the assumption about the cardinal arithmetic in $\V$ it follows
that
$\V \thinks \alef_n^{\alef_0}=\alef_n$ for $n \geq 1$. In particular
$\V\thinks |T|=\alef_1$. Moreover, since $\V \thinks
2^{\alef_1}=\alef_{\omega+1}$ we can find in $\V$ an enumeration
$\<h_\alpha : \alpha < \alef_{\omega+1}\>$ of $2^{\alef_1}$.

Given $t = \<k,l,v,\{f_{i,n},g_j: i\leq l,j\leq k, n \in
   \omega \},g,s\> \in T$ define $\Delta_t \subseteq
   \Delta$ to be collection of 
   all $\bar{p}=\<p_n: n \in \omega \>$ such that 
   \begin{enumerate}
   \item $\dom(p_n)=\{\beta_1<\dots< \beta_k\} \dot{\cup} \{\alpha^n_1<
     \dots < \alpha^n_l\}$,
   \item $p_n(\alpha^n_i)=g(i,n)$, 
   \item $p_n(\beta_i)=s(i)$,
   \item $\forall i\leq l \ \lft1(h_{\alpha^n_i}\rest v=f_{i,n}\rgt1)$,
   \item $\forall j \leq k \ \lft1(h_{\beta_j} \rest v=g_j\rgt1)$.
   \end{enumerate}

\begin{lemma}
$\Delta_t$ is filter-like for every $t \in T$.
\end{lemma}
\begin{proof}
Suppose that $\bar{p}^1, \bar{p}^2 \in \Delta_t$.
First we show that 
$f_{\bar{p}^1}\cup f_{\bar{p}^2}$ is a function.
Suppose that  $\alpha \in \dom(f_{\bar{p}^1}) \cap
\dom(f_{\bar{p}^2})$. Consider the function $h_\alpha$ and note that 
exactly one of the following possibilities happens:
\begin{enumerate}
\item there exists exactly one pair $(n,i)$ such that 
$h_\alpha \rest v=f_{i,n}$. In this case $f_{\bar{p}^1}, f_{\bar{p}^2}$
agree on $\alpha $ with the value $g(i,n)$,
\item there exists exactly one $j \leq k$ such that $h_\alpha \rest
  v=g_j$ (so $f_{\bar{p}^1}(\alpha)= f_{\bar{p}^2}(\alpha)=s(j)).$
\end{enumerate}

Now, put
$q= p_{\bar{p}_1} \cup p_{\bar{p}_2}$ 
and note that
$q$ has the required property.
 \end{proof}

To finish the proof of \ref{delta1} note that
  $\Delta=\bigcup_{t\in T} \Delta_t$.
  Suppose that $\bar{p}=\<p_n : n \in \omega \> \in \Delta$.
Let $k,l,g$ and $s$ be as in \ref{delta}, and
put $v$ to be a countable set such that 
$h_{\alpha_i^n} \rest v$ and $h_{\beta_j} \rest v$  are pairwise different.
\end{proof}

\bigskip

{\bf Finitely additive measures on $ \omega $.}

\begin{definition}
A set $ {\mathcal A}  \subseteq \Power(\omega)$ is an {\it algebra} if 
\begin{enumerate}
\item $X \cup Y \in {\mathcal A}$ whenever $X,Y \in {\mathcal A} $,
\item $\omega \setmin X \in {\mathcal A} $ whenever $X \in {\mathcal
    A} $,
\item $\emptyset, \omega \in {\mathcal A} $, $\{n\} \in {\mathcal A} $ for $ n \in \omega $.
\end{enumerate}
Given an algebra $ {\mathcal A} $,
  a function $m : {\mathcal A}  \longrightarrow [0,1]$ is a {\it finitely
  additive measure} if 
  \begin{enumerate}
  \item $m(\omega)=1$ and $m(\emptyset)=m(\{n\})=0$ for every $n$,
  \item if $X,Y \subseteq \omega $ are disjoint, then $m(X \cup Y)=m(X)+m(Y)$.
  \end{enumerate}
\end{definition}
Any non-principal filter on $\omega $ corresponds to a finitely
additive measure and any ultrafilter is a maximal such measure.

\begin{definition}
For a
real valued function $f
: \omega \longrightarrow [0,1]$ define
$$\int_\omega f \ dm = \lim_{n \rightarrow \infty} \sum_{k=0}^{2^n}
\frac{k}{2^n} \cdot m(A_k),$$ 
where 
$$A_k=\left\{n:\frac{k}{2^n} \leq f(n) <\frac{k+1}{2^n}\right\}.$$ 
\end{definition}
We leave it to the reader to verify that integration with respect to
$m$ has its usual properties.

The following is the special case of the Hahn--Banach theorem.
\begin{theorem}[Hahn--Banach]
  Suppose that $m$ is a finitely additive measure on an algebra
  ${\mathcal A} $, and $X \not \in {\mathcal A} $. Let $a \in [0,1]$
  be such that 
$$\sup\{m(A): A \subseteq X, A \in {\mathcal A} \} \leq a \leq
\inf\{m(B): X \subseteq B, B \in {\mathcal A} \}.$$
Then there exists a measure $\bar{m}$ on $\Power(\omega)$ extending
$m$ such that $\bar{m}(X)=a$.
\end{theorem}

We will need several results concerning the existence of measures in
forcing extensions.
\begin{lemma}\label{joint}
Let $m_0\in \V$ be a finitely additive measure on $ \omega $.
For $i=1,2$ let ${\mathcal R}_i$ be a forcing notion and let $
\dot{m}_i$ be a ${\mathcal R}_i$-name for a finitely additive measure
on $\V^{{\mathcal R}_i} \cap \Power(\omega) $ extending $m_0$.
Then there exists a $ {\mathcal R}_1 \times {\mathcal R}_2$-name for a
measure $ \dot{m}_2$ extending both $ \dot{m}_1$ and $ \dot{m}_2$.
\end{lemma}
\begin{proof}
We extend the measures using the Hahn-Banach theorem and we only
  need to 
  check that the requirements are consistent. 
Suppose that we have ${\mathcal R}_1$-name
  $\dot{X}$ and ${\mathcal R_2}$-name $\dot{Y}$ such that 
$\forces_{{\mathcal R}_1 \times {\mathcal R}_2} \dot{X}
  \subseteq^\star \dot{Y}$. 
A necessary and sufficient condition for both measures to have a
  common extension is that in such a case $m_1(\dot{X}) \leq
  m_2(\dot{Y})$. 
Let $(\bar{p},\bar{q}) \in {\mathcal R}_1 \times {\mathcal R}_2$ and
  $\bar{n}$ be such that 
$$(\bar{p},\bar{q}) \forces_{{\mathcal R}_1 \times {\mathcal R}_2} 
\dot{X} \setmin \bar{n} \subseteq \dot{Y}.$$
Let  
$$Z=\left\{n>\bar{n} :\exists p \in {\mathcal R}_1 \ \lft2(p \geq \bar{p} \ \&\ p \forces_{{\mathcal R}_1} n \in
  \dot{X}\rgt2)\right\}.$$
Set $Z$ belongs to $\V$ and $\bar{p} \forces_{{\mathcal R}_1} \dot{X}
  \setmin \bar{n} \subseteq Z$. Similarly $\bar{q} \forces_{{{\mathcal
  R}_2}} Z \subseteq \dot{Y}$. 
In particular, 
$$(\bar{p},\bar{q})\forces_{{{\mathcal R}_1 \times{\mathcal R}_2}}
\dot{m}_1(\dot{X}) \leq \dot{m}_1(Z)=m_0(Z)=\dot{m}_2(Z)  \leq
\dot{m}_2(\dot{Y}).$$ 
\end{proof}

We will need the following theorem.
\begin{theorem}\label{fromjech}
  Suppose that $m \in \V$ is a finitely additive atomless measure on $
  \omega $ and $ v \in \mathbf B $. For
  a $ \mathbf B $-name 
  $ \dot{X}$ for an element of $[\omega]^\omega $ define 
$$ \dot{m}^v_{\mathbf B}(\dot{X})=\sup\left\{\inf\left\{\int_\omega  \frac{\mu\lft1(q \cap \bv{ n \in
      \dot{X}}\rgt1)}{\mu(q)} \ dm : q \geq p\right\}:{p \geq v, \ p\in \dot{G}_{{\mathbf
      B}}}\right\}.$$
The name $ \dot{m}^v_{\mathbf B}$ has the following properties:
\begin{enumerate}
\item $v\forces_{\mathbf B} \dot{m}^v_{\mathbf B}: \Power(\omega)
  \longrightarrow [0,1]=\{x \in \reals: 0\leq x\leq 1\}$,
\item $v\forces_{\mathbf B} \dot{m}^v_{\mathbf B}$ is a finitely additive
  atomless measure,
\item for $X \in \V \cap \Power(\omega) \ v\forces_{\mathbf B}
  \dot{m}^v_{\mathbf B}(X)=m(X)$, 
\item if $ \dot{X}$ is a $ \mathbf B $-name for a subset of $ \omega $
  and $\mu\left(\bv{n \in \dot{X}}_{\mathbf B}\cap v\right)/\mu(v)=a>0$ for all
  $n$, then there is a condition $p \in \mathbf B, \ p\geq v $ such that 
$p \forces_{\mathbf B} \dot{m}^v_{\mathbf B}(\dot{X}) \geq a$.
\end{enumerate}
\end{theorem}
\begin{proof}
Without loss of generality we can assume that $v=2^\omega $ and
therefore we will drop the superscript $v$ altogether.

  (1) is clear.

\bigskip

(2) For a $ \mathbf B $-name $ \dot{X}$ for a subset of $ \omega $ and
$ p \in \mathbf B $ let
$$m_p(\dot{X}) = \int_\omega  \frac{\mu\lft1(p \cap \bv{ n \in
      \dot{X}}\rgt1)}{\mu(p)} \ dm$$
and
$$m^\star_p(\dot{X})=\inf\{m_q(\dot{X}): q \geq p\}.$$ 
Clearly, $ \dot{m}_{\mathbf B}(\dot{X}) = \sup_{p \in \dot{G}} \inf_{q
  \geq p} m_q(\dot{X})= \sup_{p\in \dot{G}} m^\star_p(\dot{X}).$
Note that if $p \forces_{\mathbf B} \dot{X} \subseteq \dot{Y}$ then 
$p \cap \bv{n \in \dot{X}} \subseteq p \cap \bv{n \in \dot{Y}}$ for
every $n$. It follows that
$m_p(\dot{X}) \leq m_p(\dot{Y})$ and $m^\star_p(\dot{X}) \leq
m^\star_p(\dot{Y})$. 

Similarly, if $p\forces_{\mathbf B} \dot{X} \cap \dot{Y} = \emptyset$
and $ \dot{Z}$ is a name for $ \dot{X} \cup \dot{Y}$ then 
$m_p(\dot{X})+m_p(\dot{Y})=m_p(\dot{Z})$ and
$m_p^\star(\dot{X})+m_p^\star(\dot{Y})\leq m_p^\star(\dot{Z})$.

\begin{lemma}\label{512}
  $p \forces_{\mathbf B} \dot{m}_{\mathbf B}(\dot{X}) \geq r \iff
  m^\star_p(\dot{X}) \geq r$.
\end{lemma}
\begin{proof}
  $ \leftarrow$ is obvious.

\bigskip

$ \rightarrow $. Suppose that $ p \forces_{\mathbf B} \dot{m}_{\mathbf
  B}(\dot{X}) \geq r$.
Fix any real $t < r$ and $p' \geq p$. It follows that 
$$D=\{q \geq p': m_q^\star(\dot{X})\geq t\}$$
is dense below $p'$.
Let $\{q_n: n \in \omega\}$ be a maximal antichain in $D$.
We have $m_{q_n}(\dot{X}) \geq m^\star_{q_n}(\dot{X}) \geq t$.
It follows that
$\mu(q_n \cap \bv{n \in \dot{X}}) \geq t \cdot \mu(q_n).$
Since $p' = \bigcup_n q_n$ and $q_n$ are pairwise disjoint, we get 
$$\mu(p' \cap \bv{n \in \dot{X}}) \geq t \cdot \mu(p').$$
We conclude that $m^\star_p(\dot{X}) \geq t$ and since $t$ was
arbitrary,  $m^\star_p(\dot{X}) \geq r$.
\end{proof}

Now we show that $ \dot{m}_{\mathbf B}$ is a finitely additive measure.

Suppose that $\forces_{\mathbf B} \dot{X} \subseteq \dot{Y}$.
Suppose that $p \forces_{\mathbf B} \dot{m}_{\mathbf B}(\dot{X}) > 
\dot{m}_{\mathbf B}(\dot{Y})$. Let $q \geq p$ and $r$ be such that 
$q \forces_{\mathbf B} \dot{m}_{\mathbf B}(\dot{X}) > r \geq 
\dot{m}_{\mathbf B}(\dot{Y})$. Then $m^\star_q(\dot{X}) \geq r$ and
$m^\star_q(\dot{Y}) < r$ -- contradiction.

\bigskip

Suppose that $\forces_{\mathbf B} \dot{X} \cap \dot{Y} = \emptyset$
and let $ \dot{Z}$ be  a name for $ \dot{X} \cup \dot{Y}$.
Let $p$, $r_1, r_2$ be such that 
$p \forces_{\mathbf B} \dot{m}_{\mathbf B}(\dot{X}) \geq r_1$ and 
$p \forces_{\mathbf B} \dot{m}_{\mathbf B}(\dot{Y}) \geq r_1$.
It follows that $m^\star_p(\dot{X})\geq r_1$ and
$m^\star_p(\dot{Y})\geq r_2$. Thus $m^\star_p(\dot{Z})\geq r_1+r_2$,
so 
$p \forces_{\mathbf B} \dot{m}_{\mathbf B}(\dot{Z}) \geq
\dot{m}_{\mathbf B}(\dot{X}) +
\dot{m}_{\mathbf B}(\dot{Y}).$

Suppose that
$p \forces_{\mathbf B} \dot{m}_{\mathbf B}(\dot{Z}) >
\dot{m}_{\mathbf B}(\dot{X}) +
\dot{m}_{\mathbf B}(\dot{Y}).$
There are reals $r_1, r_2$ and $q \geq p$ such that 
$q \forces_{\mathbf B} \dot{m}_{\mathbf B}(\dot{X}) < r_1$,
$q \forces_{\mathbf B} \dot{m}_{\mathbf B}(\dot{Y}) < r_2$  and 
$q \forces_{\mathbf B} \dot{m}_{\mathbf B}(\dot{Z}) > r_1+r_2$.
Use \ref{512}, to find $q' \geq q$ such that 
$m_{q'}(\dot{X}) < r_1$ and $m_{q'}(\dot{Y})< r_2$. By \ref{512},
$m_{q'}(\dot{Z}) \geq m^\star_{q'}(\dot{Z}) \geq r_1+r_2$. On the
other hand, since $m_{q'}$ is additive, $m_{q'}(\dot{Z})< r_1+r_2$ --
contradiction.

\bigskip

(3) Suppose that $ \dot{X}$ is a $ \mathbf B $-name and for some $p\in
\mathbf B $ and $ X \in \V \cap \Power(\omega)$, $p \forces_{\mathbf
  B} \dot{X}=X$. 
That means that for every $q \geq p$,
$$\frac{\mu(q \cap \bv{n \in \dot{X}})}{\mu(q)} = \left\{
  \begin{array}{ll}
1 & \hbox{if } n \in X\\
0 & \hbox{if } n \not\in X
  \end{array}\right. .$$
It follows that $p \forces_{\mathbf B} \dot{m}_{\mathbf B}(X)
\geq m(X).$
Since $ \dot{m}_{\mathbf B}$ is a measure, by looking at the
complements we get,
$p \forces_{\mathbf B} 1-\dot{m}_{\mathbf B}(X)
\geq 1-m(X)$, hence $p \forces_{\mathbf B} \dot{m}_{\mathbf B}(X)
=m(X).$

\bigskip

(4) Suppose that $\mu\left(\bv{n \in \dot{X}}_{\mathbf
    B}\right)=a>0$  for $ n \in \omega
$. 
Let 
$$D =\{p: \exists \varepsilon>0 \ m_p(\dot{X}) \leq
(1-\varepsilon)\cdot a\}.$$
If $D$ is not dense in $ \mathbf B $, then the condition witnessing
that has the required property.

So suppose that $D$ is dense and work towards a contradiction. Let
$\{q_n : n \in \omega\}$ be a maximal antichain in $D$. 
Clearly $\sum_{n=0}^\infty \mu(q_n)=1$.
Let $ \varepsilon_0 >0 $ be such that 
$m_{q_0}(\dot{X}) \leq
  (1-\varepsilon_0)\cdot a$, which means that 
$$\int_\omega \mu\left(q_0 \cap \bv{n \in \dot{X}}\right) \ dm \leq
(1-\varepsilon_0)\cdot a \cdot \mu(q_0).$$
Similarly for $n>0$,
$$\int_\omega \mu\left(q_n \cap \bv{n \in \dot{X}}\right) \ dm \leq
a \cdot \mu(q_n).$$
Let $q=\bigcup_{i\leq n} q_n$.
We have
$$\int_\omega \mu\left(q \cap \bv{n \in \dot{X}}\right) \ dm \leq
(1-\varepsilon_0)\cdot a \cdot \mu(q_0)+\sum_{j=1}^n a\cdot \mu(q_j)=
a\cdot \mu(q) - \varepsilon_0\cdot a \cdot \mu(q_0).$$
This is a contradiction since 
$$\lim_{\mu(q) \rightarrow 1}  \int_\omega \mu\left(q \cap \bv{n \in
  \dot{X}}\right) \ dm = a.$$
\end{proof}

\bigskip

{\bf The iteration}

Let $\eL$ be a model satisfying $2^{\alef_0}=\alef_1$ and
$2^{\alef_1}=2^{\alef_2}=\dots=\alef_{\omega+1}$.
In $\eL$ we will define the following objects:
\begin{enumerate}
\item A finite support iteration $\<{\mathcal P}_\alpha, \dot{{\mathcal 
      Q}}_\alpha: \alpha < \alef_{\omega+1}\>$. 
\item A sequence $\<A_\alpha : \alef_\omega \leq \alpha <
  \alef_{\omega+1}\>$.
\item A sequence $\<\dot{m}^\xi_\alpha: \alef_\omega \leq \alpha <
  \alef_{\omega+1}, \xi<\alef_1\>$ such that 
  \begin{enumerate}
  \item $ \dot{m}^\xi_\alpha $ is a ${\mathcal P}_\alpha $-name for a
  finitely additive measure on $\omega $,
\item $\dot{m}^\xi_\alpha $ extends $\bigcup_{\beta<\alpha}
  \dot{m}^\xi_\beta $,
\item if $ \cf(\gamma)>\alef_0$ then $\dot{m}^\xi_\gamma
  =\bigcup_{\beta<\gamma }  \dot{m}^\xi_\beta $.
  \end{enumerate}
\end{enumerate}
The definition is inductive. Formally, given $ {\mathcal P}_\alpha $,
$\left\{\dot{m}^\xi_\alpha: \xi<\alef_1\right\}$ and $A_\alpha $ we define
$\left\{\dot{m}^\xi_{\alpha+1}: \xi<\alef_1\right\}$ followed by $A_{\alpha+1}$
and then ${\mathcal P}_{\alpha+1}= {\mathcal P}_\alpha \star
\dot{{\mathcal Q}}_{\alpha}$. 

For limit $\alpha $, ${\mathcal P}_\alpha $ and $\{\dot{m}^\xi_\alpha:
\xi<\alef_1\}$ will be defined by the previous values and
$A_\alpha=\emptyset$.
Since the definition of $ \dot{m}^\xi_\alpha $ is most complicated it
is more natural to proceed in the reverse order by making commitments
about the defined objects as we go along. 

We will use the following notation: suppose that $\<{\mathcal
  P}_\alpha, \dot{{\mathcal Q}}_\alpha: \alpha < \delta\>$ is a finite 
support iteration and $A \subseteq \delta, A \in \eL $.
Let ${\mathcal P}(A)$ be the subalgebra generated by $\dot{G} \rest
A$ and let $\eL[\dot{G} \rest A]$ denote model $\eL[\dot{G} \cap {\mathcal P}(A)]$.  
Note that if $|A|=\alef_n$, $n>0$ then $\eL[\dot{G} \rest A] \thinks 2^{\alef_0}=\alef_n$.

To define the iteration we require that:

\begin{enumerate}
\item[A0.] $A_\alpha \subseteq \alpha $ for $\alpha<\alef_{\omega+1}$.
 
\end{enumerate}

Let $\<{\mathcal P}_\alpha, {\mathcal Q}_\alpha : \alpha <
\alef_{\omega+1}\>$ be  a finite support iteration such that 
$$\forces_\alpha {\mathcal Q}_\alpha = \left\{
  \begin{array}{ll}
{\mathbf C} & \hbox{if } \alpha< \alef_\omega \\
{\mathbf B}^{\eL[\dot{G}\rest A_\alpha]} & \hbox{if }
\alpha\geq  \alef_\omega 
  \end{array}\right. . $$

\begin{lemma}\label{extra}
  Suppose that $G$ is ${\mathcal P}_\alpha $-generic over $\eL$ and $x 
  \in \eL[G] \cap \Power(\omega)$.  Then $x$ that can 
be computed from countably many generic reals with indices in $A$.
In other words,
there exists a countable set $A
  \subseteq \alpha $, $A \in \eL$ and a 
Borel
  function $f \in {^\omega (\twoom)} \longrightarrow \twoom $, $f \in
  \eL $ and a set $\{\alpha_n: n \in \omega\} \in \eL \cap [A]^\omega$ 
  such that 
$x=f\lft1(\dot{G}(\alpha_1), \dots, \dot{G}(\alpha_n), \dots\rgt1)$.
\end{lemma}
\begin{proof}
  Induction on $ \alpha $. 

{\sc Case 1} $\alpha = \beta+1$.
Let $G \subseteq {\mathcal P}_\alpha $ be a generic filter and let $x
\in \eL[G]$.
Work in the model $\eL[G \rest \alpha]$. Since $ {\mathcal P}_\alpha = 
{\mathcal P}_\beta \star {\mathbf B}^{\eL[\dot{G}\rest A_\beta]}$
there exists a Borel function $\tilde{f} \in \eL[G \cap {\mathcal
  P}_\beta]$  such that 
$$\eL[G \cap {\mathcal P}_\beta] \thinks
\tilde{f}\lft1(G(\alpha)\rgt1)=x.$$
Since $\tilde{f}$ is coded by a real, there exists a set $A=\{\alpha_n: 
n \in \omega\}
\subseteq \beta $ and a function $f \in \eL$ such that 
$$\tilde{f} = f\lft1(G(\alpha_1), \dots, G(\alpha_n), \dots\rgt1).$$
Function $f$ and the set $A \cup \{\beta\}$ are the objects we are
looking for.

\bigskip

{\sc Case 2} $ \cf(\alpha)=\alef_0$.
Fix an increasing sequence $\<\alpha_n : n \in \omega\>$ such that
$\sup_n \alpha_n = \alpha $ and suppose that $x$ is a ${\mathcal
  P}_\alpha $-name for a real number (i.e. a set of countably many
antichains.
Let $x_n$ be a ${\mathcal P}_{\alpha_n}$-name for a real obtained by
restriction conditions in these antichains to $ \alpha_n$. 
Note that
$\forces_{{\mathcal P}_\alpha} \lim_n x_n=x$.
Apply the induction hypothesis to $x_n$'s to get Borel functions $f_n$ 
and countable sets $A_n$. Let $A = \bigcup_n A_n$ and let 
$f : {}^{\omega\times \omega }(\twoom) \longrightarrow \twoom $ be
defined as
$$f(\dots,x^n_m, \dots)=\lim_n f_n(\dots, x^n_m, \dots).$$

{\sc Case 3} $ \cf(\alpha)> \alef_0$.
Since no reals are added at the step $ \alpha $ there is nothing to prove.
\end{proof}

Furthermore, we will require that
\begin{enumerate}
\item[A1.] $|A_\alpha|< \alef_\omega $, for any $
  \alef_\omega\leq \alpha<\alef_{\omega+1}$.
  \item[A2.] For every set $A \in [\alef_{\omega+1}]^{<
    \alef_\omega}$ there are cofinally
    many $ \alpha $ with $A \subseteq A_\alpha$.
\end{enumerate}

To state the next requirement we will need the following notation: 
suppose that $A \subseteq \alef_{\omega+1}$. 
Let ${\mathcal P} 
\rest A= \{p \in {\mathcal P}  : \supp(p) \subseteq
  A\}$. Suppose that $ \dot{f} \subseteq \twoom \times \twoom$ is
  a name for an arbitrary function from 
  $\twoom $ to $\twoom $ (not necessarily Borel). Then $ \dot{f} 
\rest A =\{(\dot{x},\dot{y}) \in \dot{f}: \dot{x}, \dot{y}  \hbox{ are
  } {\mathcal P}  \rest A
\hbox{-names}\}$.

\begin{enumerate}
\item[A3.] 
$\dom(\dot{m}_{\alpha}^\xi \rest A_{\beta}
  )=\Power(\omega) \cap \eL[\dot{G} \rest A_{\beta }]$ for every
  $\xi<\alef_1$ and $\alef_\omega \leq \beta\leq \alpha<\alef_{\omega+1}$. In
  other words, $ \dot{m}_{\alpha}^\xi \rest A_{\beta } $ is a name for finitely
  additive measure on $\Power(\omega) \cap   \eL[\dot{G} \rest
  A_{\beta}]$.   
\end{enumerate}

Suppose that $\{\dot{m}^\xi_{\delta}: \xi<\alef_1\}$ is given.

Assume that $ \delta=\alpha+1$ and
that 
in order to meet the requirement A2 we have to cover a certain set $A$ of
size $\alef_n$. 
Define a sequence $\<A^\gamma_{\alpha+1}: \gamma< \omega_1\>$ such that 
\begin{enumerate}
\item $A^0_{\alpha+1} =A$,
\item $A^\beta_{\alpha+1} \subseteq A^\delta_{\alpha+1} $ for $ \beta \leq
  \delta $,
\item $A^\delta_{\alpha+1} = \bigcup_{\beta<\delta} A^\beta_{\alpha+1} $ for
  limit $ \delta $,
\item for every set $ X \in \eL[\dot{G} \rest
  A^\beta_{\alpha+1}]$ and $\xi<\alef_1$, $ \dot{m}_{\alpha+1}^\xi(X) \in \eL[\dot{G} \rest
  A^{\beta+1}_{\alpha+1}]$,
\item $|A^\gamma_{\alpha+1}|=\alef_n+\alef_1$ for all $ \gamma $.
\end{enumerate}

Note that since $ \eL[\dot{G} \rest A^\beta_{\alpha+1}] \thinks
2^{\alef_0} = \alef_n$, in order to produce
$A^{\beta+1}_{\alpha+1} $ we have to add to $A^\beta_{\alpha+1} $ at most $
\alef_n+\alef_1$ countable sets.
Finally let $A_{\alpha+1}  = \bigcup_{\gamma<\omega_1} A^\gamma_{\alpha+1} $.
It is clear that $A_{\alpha+1} $ is as required.

If $ \delta $ is limit then we put $ A_\delta=\emptyset$.
Note that in both cases
condition A3 is satisfied by the induction hypothesis and the fact that 
$\dot{m}^\xi_{\delta}$ extends 
$\bigcup_{\alpha<\delta} \dot{m}^\xi_{\alpha}.$

In order to finish the construction we have to define measures
$\{\dot{m}^\xi_\alpha: \alef_\omega\leq \alpha< \alef_{\omega+1}\}$.

We start with the definition of a certain dense subset of $ {\mathcal
  P} $ and from now on use only conditions belonging to this subset.
Let $D \subseteq {\mathcal P} $ be a subset such
that $p \in D$ if
\begin{enumerate}
\item $\dom(p) \in [\alef_{\omega+1}]^{<\omega}$,
\item $p(\alpha) \in
  {}^{< \omega}\omega \simeq \mathbf C$,  for $\alpha \in \dom(p) \cap
  \alef_\omega$,
\item for each $ \alpha \in \dom(p) \setmin \alef_\omega $,
  \begin{enumerate}
  \item   $\forces_{\alpha} p(\alpha) \in { \mathbf B}^{{
  \eL}[\dot{G}\rest A_\alpha]}$, 
\item there is a clopen set $C_\alpha  \subseteq {}^\omega
  2$ such that 
$$\forces_{\alpha} \frac{\mu\lft1(C_\alpha  \cap
  p(\alpha)\rgt1)}{\mu(C_\alpha)} \geq 1-\frac{1}{2^{n-j+5}},$$
 where $n=|\dom(p) \setmin \alef_\omega |$ and $j=|\alpha \cap
 (\dom(p) \setmin \alef_\omega) |$. 
  \end{enumerate}
\end{enumerate}

\begin{lemma}
  $D$ is dense in $ {\mathcal P} $.
\end{lemma}
\begin{proof}
  Induction on $\max(\supp(p))$.
\end{proof}

Let $\clopen$ be the collection of clopen subsets of $2^\omega $.
Represent $ {\mathbf  C}_{\alef_{\omega+1}}$  as the colection of
functions $q$ such that $\dom(q) \in [\alef_{\omega+1}]^{<\omega}$ and
$q(\alpha) \in \mathbf C $ for $ \alpha < \alef_\omega $ and
$q(\alpha) \in \clopen$ for $\alpha \geq
\alef_\omega$.

Note that there is a natural projection $\pi$ from $D$ to $ {\mathbf
  C}_{\alef_{\omega+1}}$ defined as 
$$\pi(p)(\alpha)=\left\{
  \begin{array}{ll}
p(\alpha) & \text{if } \alpha < \alef_\omega \\
C_\alpha  & \text{if } \alpha \geq \alef_\omega 
  \end{array}\right. .$$
For a sequence $\bar{p}=\<p_n: n \in \omega\>$ let
$\pi(\bar{p})=\<\pi(p_n): n \in \omega\>$. 
Suppose that $\bar{p}$ is such that $\pi(\bar{p}) \in \Delta$, as
defined in \ref{delta}.
We will define a condition $p_{\bar{p}}$ in the following way;
$\dom(p_{\bar{p}})=\widetilde{\Delta}$, where $\widetilde{\Delta}$ is
the root of the $\Delta$-system $\{\dom(p_n): n \in \omega\}$.

{\sc Case 1} $\alpha \in \widetilde{\Delta} \cap \alef_\omega $.

Let $p_{\bar{p}}(\alpha)$ be the common value of $p_n(\alpha)$ for $ n
\in \omega $.

\bigskip

{\sc Case 2} $\alpha \in \widetilde{\Delta} \setmin \alef_\omega $.

Work in the  model $\V=\eL[\dot{G} \rest A_\alpha ]$ and let
$C=\pi(p_n(\alpha))$. Clearly $\V \thinks C \in \mathbf B $.
It follows that for some $k \in \omega $ and every $ n \in \omega $,
$$\V \thinks \frac{\mu\lft1(C  \cap
  p_n(\alpha)\rgt1)}{\mu(C)} \geq 1-\frac{1}{2^{k}}.$$
Let $\dot{X}$ be a $\mathbf B $-name such that $\bv{n \in \dot{X}}=C
\cap p_n(\alpha)$.
Apply, \ref{fromjech} in $\V$, to find a condition $r \in \mathbf B,
r \geq C$ such that 
$$r \forces_{\mathbf B} \dot{m}^C_{\mathbf B}(\dot{X}_{\bar{p}}) \geq 1-\frac{1}{2^{k}}.$$
Let $p_{\bar{p}}(\alpha)=r$.

\bigskip

Now we turn our attention to the sequence $\< \dot{m}_\alpha^\xi:
\alef_\omega \leq \alpha <
\alef_{\omega+1}\>$.
Let $\{t_\xi: \xi<\alef_1\}$ be an enumeration of the set $T$ and let
 $\Delta=\bigcup_{\xi<\alef_1} \Delta_{t_\xi}$ be the decomposition
from \ref{delta1}.
For $\xi< \alef_1$ let
$$\Delta^\xi=\{\bar{p} \in [{\mathcal P}]^\omega : \pi(\bar{p}) \in
\Delta_{t_\xi}\}.$$
The measure $\dot{m}^\xi_\alpha $ will be first defined on the set
$$\left\{\dot{X}_{\bar{p}}: \bar{p} \in \Delta^\xi \cap [{\mathcal
    P}_\alpha]^\omega\right\}.$$ 
We will do it in such a way that for 
$\bar{p} \in \Delta^\xi \cap [{\mathcal
    P}_\alpha]^\omega$ 
$$p_{\bar{p}} \forces_\alpha \dot{m}^\xi_\alpha(\dot{X}_{\bar{p}}) >
0,$$
where $p_{\bar{p}}$ is the condition defined above.
Next $\dot{m}^\xi_\alpha$ will be extended arbitrarily to the set 
$\Power(\omega) \cap \eL^{{\mathcal P}_\alpha}$.

Fix $\xi< \alef_1$ and define $\dot{m}^\xi_\alpha $ as follows:

\bigskip

{\sc Case 1}. \hspace{0.1in} 
$ \alpha = \alef_\omega$.
Consider the family 
$$ \dot{{\mathcal H}}_\xi = \{ \dot{X}_{\bar{p}}: \bar{p} \in
\Delta^\xi \cap [{\mathcal P}_{\alef_\omega}]^\omega, \
p_{\bar{p}} \in \dot{G}_{{\mathcal P}}\}.$$
It is easy to see that $ \dot{{\mathcal H}}_\xi $ is a $ {\mathcal P}_{\alef_\omega}
$-name for a filter base. Let $ \dot{\F}_\xi$ be any $ {\mathcal P}
$-name for an 
ultrafilter extending $ 
\dot{{\mathcal H}}_\xi $ and let $ \dot{m}^\xi_{\alef_\omega }$ be the
corresponding measure. In other words, for $\dot{X} \in \dot{{\mathcal H}}_\xi $,
$$\forces_{\alef_\omega} \dot{m}^\xi_{\alef_\omega }(\dot{X})=1.$$

\bigskip

{\sc Case 2}. \hspace{0.1in} 
$ \alpha > \alef_\omega$ and $ 
\cf(\alpha)=\alef_0$. 

Since $ \dot{m}^\xi_\alpha $ extends $\bigcup_{\beta<\alpha}
\dot{m}^\xi_\beta $, we have to define  $ \dot{m}^\xi_\alpha $ on the set
$$ \left\{\dot{X}_{\bar{p}}: \bar{p} \in \Delta^\xi \cap \left([{\mathcal
    P}_{\alpha}]^\omega \setmin \bigcup_{\beta<\alpha} [{\mathcal
    P}_{\beta}]^\omega\right)\right\}.$$
Put $ {\mathcal A} = \Delta^\xi \cap \left([{\mathcal
    P}_{\alpha}]^\omega \setmin \bigcup_{\beta<\alpha} [{\mathcal
    P}_{\beta}]^\omega\right)$ and 
for $\bar{p} \in  {\mathcal A} $ let 
$j=j_{\bar{p}} \in \omega $ be the
such that 
$$ \beta = \sup_{n \in \omega} \alpha^n_{j-1} < \sup_{n \in \omega}
\alpha^n_j = \alpha ,$$
where $\alpha^n_i$ is the $i$'th element of $\dom(p_n)$.
Consider sequences $\bar{p}^- = \<p_n \rest \alpha^n_j : n \in
\omega\>$ and  $\bar{p}^+ = \<p_n \rest [\alpha^n_j,\alpha) : n \in
\omega\>$.
Let $\dot{{\mathcal H}}_\xi$ be a ${\mathcal P}_\alpha $-name for the
family
$\{\dot{X}_{\bar{p}^+}: \bar{p} \in {\mathcal A}\}$. 
Note that 
\begin{enumerate}
\item $\forces_\alpha \dot{{\mathcal H}}_\xi$ is a filter base,
\item $ \forall \dot{X} \in \dot{{\mathcal H}}_\xi \ \forall \beta<
  \alpha \ \forall \dot{Y} \in [\omega]^\omega  \cap \eL^{{\mathcal
      P}_\beta} \ \forces_\alpha \dot{X} \cap \dot{Y} $ is infinite.
\end{enumerate}
Suppose that $\bar{p} \in {\mathcal A} $ and note that 
$$p_{\bar{p}^-} \forces_\beta \ \dot{m}^\xi_\beta(\dot{X}_{\bar{p}^-})=a>0.$$
By the remarks made above, we can set
$\dot{m}_\alpha(\dot{X}_{\bar{p}^+})=1$ and
$\dot{m}_\alpha(\dot{X}_{\bar{p}})=a$. Finally note that the value $a$
is forced by
$p_{\bar{p}}$.

\bigskip

{\sc Case 3}. \hspace{0.1in} 
$ \alpha $ is limit and $ \cf(\alpha )> \alef_0$.

Let $ \dot{m}_\alpha^\xi = \bigcup_{\beta<\alpha } \dot{m}_\beta^\xi $.
This definition is correct since no subsets  of $ \omega $ are added
at the step $ \alpha $.

\bigskip

{\sc Case 4}. \hspace{0.1in} $ \alpha=\delta+1$. 

As before we have to define $\dot{m}^\xi_\alpha $ on 
$$\left\{\dot{X}_{\bar{p}}: \bar{p} \in \Delta^\xi \cap \left([{\mathcal 
    P}_{\alpha}]^\omega \setmin [{\mathcal
    P}_{\delta}]^\omega\right)\right\}.$$
Set $ {\mathcal A} = \Delta^\xi \cap \left([{\mathcal 
    P}_{\alpha}]^\omega \setmin [{\mathcal
    P}_{\delta}]^\omega\right)$ and note that if $ \bar{p} \in
{\mathcal A} $ then $\delta \in \bigcap_{n \in \omega}
\dom(p_n)$. Thus, let $C$ be a clopen set such that 
$\pi(p_n(\delta))=C$ for $n \in \omega $.
Let $\V= {
  \eL}[\dot{G}\rest A_\delta]$.
Find a forcing notion
${\mathcal R}$ such that 
$ {\mathcal P}_\delta   = ({\mathcal P}_\delta  \rest
  A_\delta) \star {\mathcal R} $. It follows that
$ \eL^{{\mathcal P}_\alpha }= \eL^{{\mathcal
  P}_{\delta+1}} ={\V}^{{\mathcal R}\times {\mathbf B}}$.
By the induction hypothesis $m= \dot{m}_\delta^\xi \rest A_\delta$ is a finitely additive 
  measure. In other words $m \in \V$ is a
  finitely additive measure defined on $\Power(\omega) \cap \V$. 
Clearly $ \dot{m}_\delta^\xi $ is an extension of $m$ to $\V^{\mathcal R}
  \cap \Power(\omega)$. On the other hand let $\dot{m}^C_{\mathbf B}$
  be an extension of $m$ to $\V^{\mathbf B}$ as given by
  \ref{fromjech}.  
Let $\dot{m}^\xi_\alpha=\dot{m}^\xi_{\delta+1}$ be the common extension of $
\dot{m}_\delta^\xi $ and $\dot{m}^C_{\mathbf B}$ guaranteed by \ref{joint}.
It is clear that $\dot{m}^\xi_\alpha$ has the required properties.

Finally let $\dot{m}^\xi = \bigcup_{\alef_\omega \leq \alpha<
  \alef_{\omega+1}} \dot{m}^\xi_\alpha$. Note that each $\dot{m}^\xi$
is a ${\mathcal P} $-name for a finitely additive measure on
$\Power(\omega) \cap \eL^{{\mathcal P}}$.

\bigskip

{\bf Proof of the Theorem \ref{sh592}}

We are ready now for the proof of the main theorem.
The following lemma gives the lower bound for $\cov(\N)$.
\begin{lemma}
$ \eL^{{\mathcal P} } \thinks \cov(\N) \geq 
\alef_\omega$. In particular, $[\reals]^{<\alef_\omega} \subseteq
\COV(\N).$ 
\end{lemma}
\begin{proof}
 Suppose that $\{H_\alpha : \alpha<\kappa <
 \alef_\omega\}$ is a family of measure zero sets in $
 \eL^{{\mathcal P} }$.
Let $N$ be a master set for $\N$ defined earlier. Without loss of
 generality we can assume that for some $f_\alpha \in \omom $,
$H_\alpha =(N)_{f_\alpha}$, and let $ \dot{f}_\alpha $ be a $
 {\mathcal P} $-name for $f_\alpha $. As in \ref{extra},
let $K_\alpha \in
 [\alef_{\omega+1}]^{\alef_0}\cap \eL$ be the set such that $f_\alpha
 \in \eL[\dot{G} \rest K_\alpha]$. 
Find $ \beta$ such that $\bigcup_{\alpha<\kappa} K_\alpha \subseteq
 A_\beta $.
The random real added by $ {\mathbf B}^{\eL[\dot{G} \rest
 A_\beta]}$ avoids all null sets coded in $ \eL[\dot{G} \rest
 A_\beta]$, in particular, all $H_\alpha$'s.  
\end{proof}

It remains to be checked that $\cov(\N) \leq \alef_\omega $ in the extension.

Let $X = \{f_\alpha : \alpha < \alef_\omega \}=\dot{G}\rest
\alef_\omega$ be the
sequence of first $\alef_\omega $ Cohen reals added by $ {\mathcal P} $. 
Our intention is to show that $X \not \in \COV(\N)$.
In fact we will show that
$$\bigcup_{\alpha< \alef_\omega} (N)_{f_\alpha} = \twoom,$$
where $N$ is the master set defined in the previous section. 
That will finish the proof since $X$ is a countable union of sets of
smaller size (so they are all 
in $\COV(\N)$) and thus $X$ witnesses that $\COV(\N)$ is not a $
\sigma $-ideal and that $\cov(\N) \leq
\alef_\omega$.

Suppose the opposite  and let $z$ be such that 
$$ \eL^{{\mathcal P} } \thinks z \not\in
  \bigcup_{\alpha< 
  \alef_\omega} (N)_{f_\alpha}.$$

\begin{lemma}\label{4.3}
  There exists a $ {\mathcal P}   $-name $ \dot{Y}$ for a subset of $
  \alef_\omega $ and $\bar{n} \in \omega $ such that 
  \begin{enumerate}
  \item $\forces_{{\mathcal P} } \dot{Y} \in
    [\alef_\omega]^{\alef_1}$,
  \item $\forces_{{\mathcal P} } \twoom \setmin \bigcup_{\alpha \in
      \dot{Y}} \bigcup_{n>\bar{n}} C^n_{f_\alpha(n)}$ is uncountable.
  \end{enumerate}
\end{lemma}
\begin{proof}
Denote by $ \dot{z}$ a $ {\mathcal P}   $-name for
$z$ and let $ \delta< \alef_{\omega+1}$ be the least
ordinal such that $ \dot{z} $ is a $ {\mathcal P}_\delta $-name. 
We have the following two cases:

\vspace{0.1in}

%{\sc Case 1} \hspace{0.1in} $ \delta< \alef_\omega $. 

%Since $ \eL^{{\mathcal P}}$ is obtained by adding $
%\alef_\omega $ Cohen reals we have

%\bigskip

{\sc Case 1}. \hspace{0.1in} $ \delta = \lambda +1$ is a
successor ordinal.

Suppose first that $ \delta> \alef_\omega$. Work in $\V={
  \eL}^{{\mathcal P}_\lambda}$ and let $ {\mathbf B}_\lambda 
= {\mathbf B}^{\eL[\dot{G}\rest A_\lambda]}$.  
For each $ \alpha <
\alef_\omega $ choose $q_\alpha  \in {\mathbf B}_\lambda $ and
$n_\alpha\in
\omega $ such that $\V \thinks q_\alpha \forces_{{\mathbf B}_\lambda}
\dot{z} \not \in \bigcup_{n>n_\alpha } 
C^n_{f_\alpha(n)} $. Since 
$ {\mathbf B}_\lambda $ has a dense subset of size $<
\alef_\omega $, we can find $q \in {\mathbf B}_\lambda $ and
$\bar{n} \in \omega $ such that the set 
$$Y=\{\alpha : q_\alpha = q \ \&\ n_\alpha = \bar{n}\}$$ 
is uncountable.
Consider the  set
$C=\twoom \setmin \bigcup_{\alpha\in Y} \bigcup_{n>\bar{n}}
C^n_{f_\alpha(n)}$ in $\V$.
Observe that $C$ is a closed set and if it was countable then all its
elements would be in $\V$. However, ${\V}^{{\mathbf B}_\lambda}
 \thinks z \in C$ and $z \not\in \V$. 

\bigskip

If $ \delta < \alef_\omega$ then the argument is identical except
that we use $ {\mathbf C} $ instead of $ {\mathbf B}_\lambda $.
In fact one can show that
$$\eL^{\mathcal P} \cap \twoom \subseteq
\bigcup_{\alpha<\omega_1} (N)_{f_\alpha} \subseteq
\bigcup_{\alpha<\alef_\omega} (N)_{f_\alpha}.$$

\vspace{0.1in}

{\sc Case 2}. \hspace{0.1in} $ \delta$ is limit and $
\cf(\delta)=\alef_0$.

In $\eL^{{\mathcal P}_\delta}$ we can find
$\bar{n} \in \omega $ and an uncountable set $Z \subseteq
\alef_\omega $
such that  
$$\eL^{{\mathcal P}_\delta} \thinks z \not\in \bigcup_{\alpha\in Z}
\bigcup_{n>\bar{n}} C^n_{f_\alpha(n)}.$$
Let $ \dot{Z}$ be a $ {\mathcal P}_\delta $-name for $ Z$.
Suppose that $G \subseteq {\mathcal P}_\delta $ is a generic filter
over $\eL$.
For each $ \alpha < \omega_1$ choose $p_\alpha \in {\mathcal P}_\delta
\cap G$ and $\eta_\alpha$
such that $p_\alpha \forces_{{\mathcal P}_\delta }
\dot{Z}(\alpha)=\eta_\alpha$, where $ \dot{Z}(\alpha)$ is a $
{\mathcal P}  $-name for the  $ \alpha $-th element of $Z$.

There is an uncountable  set $I \subseteq \omega_1 $, 
and $
\lambda < \delta $ such that $p_\alpha \in {\mathcal P}_\lambda \cap
G$ for $ \alpha \in I$.  Let 
$Y=\{\eta_\alpha : \alpha \in I\}$ and let $ \dot{Y}$ be a $
{\mathcal P}_\lambda $-name for $Y$.  As in the previous
case, consider the set $C=\twoom \setmin \bigcup_{\alpha\in Y}
\bigcup_{n>\bar{n}} C^n_{f_{\eta_\alpha}(n)}$ in $\eL^{{\mathcal
    P}_\lambda }$. We see 
that $C$ is uncountable because it contains an element
which does not belong to $\eL^{{\mathcal P}_\lambda
  }$.
\end{proof}

Find different 
ordinals $\{\eta_\alpha : \alpha <\omega_1\}$ and conditions
$\{p_\alpha :
\alpha < \omega_1\} \subseteq {\mathcal P}  $ such that 
$p_\alpha  \forces_{{\mathcal P}}
\eta_\alpha  \in \dot{Y}.$ 
Using the $\Delta$-lemma we can assume that there are $\widetilde{k},
\widetilde{l} \in \omega $, $s \in {}^{\widetilde{k}} \omega$ and clopen 
sets $\{C_j: j \leq \widetilde{l}\}$ such that  
\begin{enumerate}
\item $\supp(p_\alpha)$ form a $\Delta$-system,
\item $\supp(p_\alpha)=
\{\gamma^\alpha_1< \cdots < 
  \gamma^\alpha_{\tilde {k}}< \alef_\omega \leq
  \delta^\alpha_1<\cdots< 
  \delta^\alpha_{\tilde{l}}\}$, 
\item $\forall \alpha \ \forall j \leq  \widetilde{k}\
 \lft1( p_\alpha(\gamma^\alpha_j)=s(j)\rgt1)$,
\item for all $j \leq \widetilde{l}$ 
$$\forces_{\alpha_j} \frac{\mu\lft1(C_j \cap p_\alpha(\delta^\alpha_j)\rgt1)}{\mu(C_j)} \geq
1-\frac{1}{2^{\tilde{l}-j+5}}.$$
\end{enumerate}
Without loss of generality we can assume that $\eta_\alpha \in
\dom(p_\alpha)$. Furthermore we can assume that for some $j_0 \leq
\widetilde{k}$, $\eta_\alpha = \gamma^\alpha_{j_0}$ and that
$s(j_0)=s^\star$ with $|s^\star|=n^\star$.

Consider the first $ \omega $ conditions $\bar{p}=\{p_n: n \in \omega \}$ 
Our next step is to extend $p_n$'s slightly to get a new sequence
$\bar{p}^\star$.
We will need the following definition.
\begin{definition}
For a clopen set $C \subseteq \twoom$ define 
$\suppp(C)$ to be the smallest set $F \subseteq \omega $ such that 
$C=(C \cap {{}^F 2}) \times {^{\omega \setmin F}}2$. In other words,
support of $C$ is the set of  coordinates that carry information about
$C$. 
  \end{definition}

Let $K_n =\left\{m: \suppp\left(C^{{n^\star}}_m\right) \subseteq
  n\right\}$ and let $\{J_n: n \in \omega\}$ be a partition of $
\omega $ such that $|J_n|=|K_n|$ for each $n$. Fix a function $o \in
\omom$ such that $o{``}(J_n)=K_n$ for every $n$.
Define
$$p^\star_n=\left\{\begin{array}{ll}
p_n(\alpha) & \hbox{if } \alpha \neq \eta_n\\
{s^\star}{}^\frown \lft1(n^\star, o(n)\rgt1)& \hbox{if } \alpha=\eta_n
\end{array}\right. .$$

Observe that there is $\xi< \alef_1$ such that
$\bar{p}^\star=\{p_n^\star: n \in \omega\} \in \Delta^\xi$. 
This is being witnessed by the $\widetilde{k},
\widetilde{l}$, $s \in {}^{\widetilde{k}} \omega$, clopen 
sets $\{C_j: j \leq \widetilde{l}\}$ and function $g$ defined as
$$g(i,n)=\left\{
  \begin{array}{ll}
s(i)& \text{if } i \leq \widetilde{k}, \ i \neq j_0\\
{s^\star}{}^\frown \lft1(n^\star, o(n)\rgt1)& \text{if } i=j_0
  \end{array}\right. .$$

Our goal is to show:
\begin{theorem}\label{crucial}
There exists a condition $p^{\star\star}$ and $\varepsilon>0$ such that 
$$p^{\star\star} \forces_{{\mathcal P}}
  \exists^\infty n \ \frac{\left|\left\{m\in J_n: 
  p^\star_m \in \dot{G}_{{\mathcal P}}\right\}\right|}{\left|J_n\right|}
  \geq     
  \varepsilon .$$ 
\end{theorem}
Before we prove this theorem let us see that the theorem follows readily
from it.
Recall that in Lemma~\ref{4.3} we showed that
$\forces_{{\mathcal P}} \twoom \setmin \bigcup_{\alpha \in
      \dot{Y}} \bigcup_{n>\bar{n}} C^n_{f_\alpha(n)}$ is uncountable.
Since this set is closed,  there is a $ {{\mathcal P}} $-name for a tree $ \dot{T}$ 
 such that 
$\forces_{{\mathcal P}} \bigcup_{\alpha \in
      \dot{Y}} \bigcup_{n>\bar{n}} C^n_{f_\alpha(n)} \cap
    [\dot{T}]=\emptyset$.
Let $\dot{Z}_n=\left\{m\in J_n: 
  p^\star_m \in \dot{G}_{{\mathcal P}}\right\}$ for $ n \in \omega $.
It follows that for every $n$,
$$p^{\star\star} \forces_{{\mathcal P}} \left(\bigcup_{k\in \dot{Z}_n}
  C^{{n^\star}}_k \right) \rest n \cap \dot{T} \rest n = \emptyset.$$
This is because for a clopen set $C$ and a tree $T$,
if $C \cap [T]=\emptyset$ then $\lft1(C \rest \suppp(C)\rgt1) \cap \lft1(T\rest
  \suppp(C)\rgt1)=\emptyset$. 
 Fix $n \in \omega $ and suppose that $ |\dot{T}\rest n| = m$.
The size of the set $J_n$ is equal to $\left(\begin{array}{c} {2^n}\\
{2^{n-{n^\star}}}\end{array}\right)$. On the other hand the number of sets
  $C^{{{n^\star}}}_k$ which are disjoint with $ \dot{T} \rest n$ is at
  most $\left(\begin{array}{c} {2^n-m}\\
{2^{n-{n^\star}}}\end{array}\right)$. Put $ 2^{-{n^\star}}=\epsilon$.
It follows, (after some calculations), that for some constant $a\geq 1$:
$$\frac{|\dot{Z}_n|}{|J_n|} \leq  \frac{\left(\begin{array}{c}{2^n-m}\\
{2^{n-{n^\star}}}\end{array}\right)}{\left(\begin{array}{c}{2^n}\\
{2^{n-{n^\star}}}\end{array}\right)} = \prod_{j=1}^m
  \left(1-\frac{2^{n-{n^\star}}}{2^n-m+j}\right) \leq a\cdot e^{-\epsilon
  \cdot m}.$$
Thus
$$\frac{|\dot{Z}_n|}{|J_n|} \leq a\cdot e^{-\epsilon \cdot
  |\dot{T}\rest n|}.$$
Since $p^{\star\star} \forces_{{\mathcal P}}  \limsup_n \displaystyle
\frac{|\dot{Z}_n|}{|J_n|}\geq
  \varepsilon $ we get that
$p^{\star\star} \forces_{{\mathcal P}} \lim_n
  |\dot{T} \rest n| < \infty$ (the size of $T\rest n$ increases with $n$). 
In particular, 
$$p^{\star\star} \forces_{{\mathcal P}} \dot{T} \hbox{ is not
  perfect,}$$
 which gives a contradiction.

 \bigskip

{\bf Proof of the Theorem \ref{sh592}: conclusion}

In order to finish the proof of \ref{sh592} we have to prove \ref{crucial}.
We will need one more modification of the sequence $\bar{p}^\star$ and 
we will require the construction described below.

\begin{lemma}\label{lastone}
Let $\widetilde{\Delta}$ be a finite subset of $\alef_{\omega+1}
\setmin \alef_\omega $. 
Suppose that $\{q_i: i \leq N\}$ is a sequence of conditions in $
{\mathcal P} $ such that 
\begin{enumerate}
\item $\dom(q_i)=\widetilde{\Delta}$,
\item $\forall \alpha \in \widetilde{\Delta} \ \exists a_\alpha \
  \forall i \leq N \ \forces_\alpha
  \mu\lft1(q_i(\alpha)\rgt1)=a_\alpha>3/4$.
\end{enumerate} 

There exists a condition $ q^\star$ such that 
\begin{enumerate}
\item $\dom(q^\star)=\widetilde{\Delta}$,
\item $q^\star \in {\mathcal P} $,
\item $ \forall \alpha \in \widetilde{\Delta} \ \forces_\alpha
  \mu\lft1(q^\star(\alpha)\rgt1) \geq 2a_\alpha -1,$
\item  $q^\star \forces_{{\mathcal P}} \{k \leq N: \forall \alpha \in
  \widetilde{\Delta}\ 
  q^\star \rest \alpha \forces_{\alpha} 
  \dot{x}_\alpha  \in q_k(\alpha)\}$ has at least
  $2^{-|\widetilde{\Delta}|}\cdot N\cdot \prod_{\alpha\in
    \widetilde{\Delta}} a_\alpha $ elements, where $\dot{x}_\alpha$ is 
  the generic real added by $\dot{G}(\alpha)$.
\end{enumerate}
\end{lemma}
\begin{proof}
If $\widetilde{\Delta}=\emptyset$, then there is nothing to prove.

Suppose that 
$\widetilde{\Delta}\neq \emptyset$ and let $\beta=\max(\widetilde{\Delta})$.
Let $q'_k = q_k \rest \beta $ for $ k\leq N$.
Apply the induction hypothesis to get a condition $q'$ such that 
\begin{enumerate}
\item $\dom(q')=\widetilde{\Delta}\setmin \{\beta \}$,
\item $q' \in {\mathcal P} $,
\item $ \forall \alpha \in \widetilde{\Delta}\setmin \{\beta\} \ \forces_\alpha
  \mu\lft1(q'(\alpha)\rgt1) \geq 2a_\alpha -1,$
\item  $q' \forces_{{\mathcal P}} \{k \leq N: \forall \alpha \in
  \widetilde{\Delta}\setmin \{\beta\}\ 
  q' \rest \alpha \forces_{\alpha} 
  \dot{x}_\alpha  \in q_k(\alpha)\}$ has at least
  $2^{-|\widetilde{\Delta}|}\cdot N\cdot \prod_{\alpha\in
    \widetilde{\Delta}\setmin \{\beta \}} a_\alpha $ elements.
\end{enumerate}

Let $ \dot{W}$ be a $ {\mathcal P} $-name for the set 
$$\{k \leq N: \forall \alpha \in
  \widetilde{\Delta}\setmin \{\beta\}\ 
  q' \rest \alpha \forces_{\alpha} 
  \dot{x}_\alpha  \in q_k(\alpha)\}.$$
Let $\{W^i: i \leq \ell\}$ be the list of all subsets of $N$ of size
at least $2^{-|\widetilde{\Delta}|}\cdot N\cdot \prod_{\alpha\in
    \widetilde{\Delta}\setmin \{\beta \}} a_\alpha $.
Find a maximal antichain $\{q^i: i \leq \ell\}$ below $q'$ such that 
$q^i \forces_{{\mathcal P}} \dot{W}=W^i$ for $i \leq \ell$.

We will need the following easy observation.
\begin{lemma}\label{triv} 
  Suppose that $\{A_n: n<N\}$ is a family of subsets of $\twoom $ of
  measure $a>0$. Let
$$B=\left\{x \in \twoom : x \hbox{ belongs to at least
    $\displaystyle\frac{N\cdot a}{2}$ sets $A_i$}\right\}.$$
Then $\mu(B) \geq \max\{a/2, 2a-1\}$.
\end{lemma}
\begin{proof} Let $\chi_{A_{i}}$ be the characteristic function of the set $A_{i}$ for
$i \leq N$.
It follows that $\int \sum_{i \leq N} \chi_{A_{i}} = N
\cdot a$. On the other hand, estimation of this integral yields,
$$ N \cdot \mu(B) + \frac{N\cdot a}{2}\lft1(1-\mu(B)\rgt1) \geq N
\cdot a$$
and after simple computations we get
$\mu(B) \geq \displaystyle\frac{a/2}{1-a/2}.$ It follows that we get
the following estimates:
$$\mu(B) \geq \left\{
  \begin{array}{ll}
a/2 & \hbox{if } a \hbox{ is close to } 0\\
2a-1 & \hbox{if } a \hbox{ is close to } 1
  \end{array}\right. .$$
 \end{proof}
 
Work in $\V^{{\mathcal P}_\beta}$ and for each $i \leq \ell$ apply
\ref{triv} to the 
family $\{q_k(\beta): k \in W^i\}$ and obtain a condition $r^i \in
\mathbf B^{\eL[\dot{G} \rest A_\beta]}$ such
that 
$$r^i \forces \{k \in W^i:   \dot{x}_\beta   \in q_k(\beta)\} \text{
  has at least } \frac{|W^i|}{2}\cdot a_\beta \text{ elements},$$
and $\forces_\beta \mu(r^i)\geq 2a_\beta -1$.

Finally, define
$q^\star$ to be a $ {\mathcal P} $-name such that for $i \leq \ell$,
$q^i \forces q^\star(\beta)=r^i$.
It is easy to see that $q^\star$ is as required.
\end{proof}

Let $q_k = p^\star_k \rest \widetilde{\Delta}$, where
$\widetilde{\Delta}=\{\alpha_1<\dots <\alpha_\ell\}$ is the root of
the $\Delta$-system $\left\{\dom\lft1(p_k^\star \rest
  [\alef_\omega,\alef_{\omega+1}\rgt1): k\in \omega 
  \right\}$. 

For each $n$ apply  \ref{lastone} to the family $\{q_k: k \in
J_n\}$ to get a condition $q_n^\star$ such that 
\begin{enumerate}
\item $\dom(q_n^\star)=\widetilde{\Delta}$,
\item $ \forall i\leq \ell\ \forces_{\alpha_i}
  \displaystyle\frac{\mu(q^\star(\alpha_i) \cap C_{\alpha_i})}{\mu(C_{\alpha_i})} 
  \geq 2\left(1-\frac{1}{2^{\ell-i+5}}\right)-1=\frac{1}{2^{\ell-i+4}},$
\item  $q^\star_n \forces_{{\mathcal P}} \left|\left\{k \leq N: \forall \alpha \in
  \widetilde{\Delta}\ 
  \left(q^\star_n \rest \alpha \forces_{\alpha} 
  \dot{x}_\alpha  \in q_k(\alpha)\right)\right\}\right|\geq \displaystyle\frac{|J_n|}{2^{\ell+1}}$. 
\end{enumerate}

Define for $ k \in \omega $,
$$p_k^{\star\star}(\alpha)=\left\{
  \begin{array}{ll}
p^\star_k(\alpha) \cap q_n^\star& \text{if } \alpha\in \widetilde{\Delta}, \ k \in
J_n\\
p^\star_k(\alpha)& \text{otherwise}
  \end{array}\right. . $$
Let $\bar{p}^{\star\star}=\{p^{\star\star}_n: n \in \omega\}$.
Find $\xi<\alef_1$ such that $\bar{p}^{\star\star} \in \Delta^\xi$.
According to our definitions, 
$$p^{\star\star}=p_{\bar{p}^{\star\star}} \forces_{{\mathcal P}}
\dot{m}^\xi(\dot{X}_{\bar{p}^{\star\star}})>0.$$
In particular,
$$p^{\star\star} \forces_{{\mathcal P}}
\dot{X}_{\bar{p}^{\star\star}}=\{n: p^{\star\star}_n \in
\dot{G}_{\mathcal P}\}\text{ is infinite.}$$ 
Let $ \varepsilon  = 2^{-\ell-1}$ and note that 
$$p^{\star\star}_n \forces_{{\mathcal P}} \frac{\left|\{k \in J_n:
    p^\star_k \in \dot{G}_{\mathcal P}\}\right|}{|J_n|} \geq
\varepsilon .$$
It follows that $p^{\star\star}$ is the condition required in \ref{crucial}.

\bigskip

{\bf Historical remarks}
Theorem~\ref{mill} was proved by Miller \cite{Mil82Bai}.  
Better estimates are true (see \cite{BarJud93Cov} and
\cite{BarIho89Cof} or \cite{BJbook}).
Theorem
\ref{covandb} was proved in \cite{Bar88Cov} (see \cite{BJbook}).
Theorem~\ref{sh592} is due to Shelah. His paper \cite{Sh592} contains a
more a general construction, where in addition ${\mathbf {MA}}_{\alef_1}$
holds.

\section{Consistency results and counterexamples}
This section is devoted to the consistency results involving cardinal
invariants of the Cicho\'n diagram and non-inclusion between the
corresponding classes of small sets.
We will describe several such constructions in  detail.

Suppose that $ {\mathcal P} $ is a forcing notion.
  Let ${\mathcal D}({\mathcal P})$ denote the family of all dense
  subsets of elements of ${\mathcal P}$ 
and ${\mathcal G}({\mathcal P})$ the family of all filters on 
 $ {\mathcal P} $.
With $ {\mathcal P} $ we can associate 
the following cardinal invariants:
\begin{enumerate}
\item$\m({\mathcal P})=\min\{|{\mathcal A}|: {\mathcal A} \subseteq
    {\mathcal D}({\mathcal P}) \ \&\ \neg \exists G \in {\mathcal
    G}({\mathcal P}) \ \forall D \in {\mathcal A} \ (G \cap D \neq
    \emptyset)\} $,
\item $\w({\mathcal P})=
\min\{|\mathcal G|:{\mathcal G}\subseteq{\mathcal G}({\mathcal P})\
\&$ for every countable sequence $\{D_n: n \in \omega \} \subseteq
{\mathcal D}({\mathcal P}) \ \exists G\in{\mathcal G}({\mathcal P}) \
\forall n \ (G\cap 
D_n\neq\emptyset)\}. $
\end{enumerate}

In other words, $\m({\mathcal P})$ is the size of the smallest family of
dense subsets of ${\mathcal P}$ for which there is no filter intersecting
all of them and $\w({\mathcal P})$ is the size of the smallest family of
filters such that for every countable family of dense  subsets of
${\mathcal P}$ there is a filter in the family that intersects all of
them.
 
Consider the forcing notions:
\begin{itemize}
\item Amoeba forcing ${\mathbf A}=\{U \subseteq \twoom : U \hbox{ is
    open and } \mu(U)>1/2\}$. For $U, V \in {\mathbf A}$, $U \geq V$
    if $U \supseteq V$,
  \item Random real forcing $ {\mathbf B} =\{P \subseteq \twoom : P$
    is a closed set of positive measure $\}$.
  \item Cohen forcing $ {\mathbf C} $.
  \item Dominating real forcing ${\mathbf D} = \{\<n,f\> : n \in \omega, f \in \omom \}.$
 For
$\<n,f\>, \<m, g\> \in {\mathbf D}$ let 
$\<n, f\> \geq \<m,g\>$ if $n
\geq m \ \&\ f \rest m = g \rest m \ \& \ \forall k \ f(k) \geq g(k).$
 \end{itemize}

We have the following result (see \cite{BJbook} for the proof):
\begin{theorem}\label{4chars}
  \begin{enumerate}
  \item $\add(\N)=\m({\mathbf A})$ and $\cof(\N)=\w({\mathbf A}).$
  \item $\cov(\N)=\m({\mathbf B})$ and $\unif(\N)=\w({\mathbf B}).$
  \item $\cov(\M)=\m({\mathbf C})$ and $\unif(\M)=\w({\mathbf C}).$
  \item $\add(\M)=\m({\mathbf D})$ and $\cof(\M)=\w({\mathbf D}).$
  \end{enumerate}
\end{theorem}

This description is particularly well suited to use with the finite
support iteration. If $ {\mathcal P} $ is a ccc forcing notion having
``nice'' definition and $ {\mathcal P}_\kappa $ is a finite support
iteration of $ {\mathcal P} $ of length $ \kappa $ then 
\begin{enumerate}
\item If $\V \thinks 2^{\alef_0}=\alef_1 $
  then $\V^{{\mathcal P}_{\omega_2}} \thinks \m({\mathcal
  P})=\alef_2$,
\item If $\V \thinks 2^{\alef_0} =\alef_2 $
  then $\V^{{\mathcal P}_{\omega_1}} \thinks \w({\mathcal
  P})=\alef_1$.
\end{enumerate}
This example motivates the following definition: a pair of models $\V$
and $\V'$ is dual if
$$\V \thinks \m({\mathcal P})=2^{\alef_0} \iff \V' \thinks
\w({\mathcal P})< 2^{\alef_0} .$$
For our purpose we restrict our attention
 to the coefficients of the Cicho\'n
diagram and define that $\V$ is dual to $\V'$ if all of the following
hold: 
\begin{enumerate}
\item $\V \thinks \cov=2^{\alef_0} \iff \V' \thinks \unif<
  2^{\alef_0} $,
\item $\V \thinks \add=2^{\alef_0} \iff \V' \thinks \cof<
  2^{\alef_0} $,
\item $\V \thinks \unif=2^{\alef_0} \iff \V' \thinks \cov<
  2^{\alef_0} $,
\item $\V \thinks \cof=2^{\alef_0} \iff \V' \thinks \add<
  2^{\alef_0} $,
\item $\V \thinks {\mathfrak b}=2^{\alef_0} \iff \V' \thinks
  {\mathfrak d}<
  2^{\alef_0} $,
\item $\V \thinks {\mathfrak d}=2^{\alef_0} \iff \V' \thinks
  {\mathfrak b}<
  2^{\alef_0} $.
\end{enumerate}

To illustrate this consider the following theories:

$$\ZFCa + \add(\M)=\cov(\N)= \alef_2 +
\add(\N)=\alef_1 $$
and 
$$\ZFCa + \cof(\N)=\alef_2 + \cof(\M)=\unif(\N)=
\alef_1. $$
The first of these models can be obtained by a finite support
iteration of $ {\mathbf B} \star {\mathbf D}$ of length $
\alef_2 $ over a model for $\CH$ and the second by 
iteration of $ {\mathbf B} \star {\mathbf D}$ of length $
\alef_1 $ over a model for $\neg \CH$. 
It is clear that $\add(\M)$, $\cov(\N)$ and $\cof(\M)$ and $\unif(\N)$
have the required values. What is less obvious is that $\add(\N)=
\alef_1 $ in the first and $\cof(\N)= \alef_2
$ in the second case. 
To check that we need a preservation result which ensures that the
iteration which we use does not change the value of these invariants. 
Such theorems were proved in \cite{JudShel90Kun}, \cite{RepProp} and
\cite{BJbook}. 

We will not study these examples any further because this method has
one fundamental weakness: it can give us only some of the models we
need. This is because the finite support iteration adds Cohen reals.
 We will use however the notion of duality outlined above.
From now on we will focus on obtaining the models using countable
support iteration. 
To this end we will associate with every cardinal invariant of the
Cicho\'n diagram a proper forcing notion and a preservation theorem as
follows:
\begin{itemize}
\item $\add(\N)$ $\leftrightsquigarrow$ Amoeba forcing $\mathbf A$,
  preservation of ``not adding amoeba reals''
\item $\cov(\N)$ $\leftrightsquigarrow$ random real forcing ${\mathbf
    B}$,   preservation of ``not adding random reals''
\item $\cov(\M)$ $\leftrightsquigarrow$ Cohen forcing $\mathbf C$,   preservation of ``not adding Cohen reals''
\item $\unif(\M)$ $\leftrightsquigarrow$ forcing ${\mathbf {PT}}_{f,g}$,
  preservation of non-meager sets,
\item $\mathfrak b$ $\leftrightsquigarrow$ Laver forcing ${\mathbf {LT}}$, preservation of ``not adding unbounded reals''
\item $\mathfrak d$ $\leftrightsquigarrow$ rational perfect set forcing
  $\mathbf {PT}$, preservation of ``not adding dominating reals''
\item $2^{\alef_0}$ $\leftrightsquigarrow$ Sacks forcing
  $\mathbf S$, preservation of Sacks property,
\item $\cof(\N)$ $\leftrightsquigarrow$ forcing ${\mathbf {S}}_2$,
\item $\unif(\N)$ $\leftrightsquigarrow$ forcing ${\mathbf
  {S}}_{g,g^\star}$, preservation of positive outer measure.
\end{itemize}

We do not assign anything to $\add(\M)$ and $\cof(\M)$ because they
are expressible using the remaining invariants.
We refer the reader to \cite{BJbook} for the definitions of all these
forcing notions and the formulation of the preservation theorems. 
We will illustrate the problems with the following  examples.

\begin{example}
Dominating number $\mathfrak d$. 
Rational perfect set forcing ${\mathbf {PT}}$ associated with ${\mathfrak
  d}$ is one of the forcing notions that increase ${\mathfrak d}$ without
affecting other characteristics. 
The preservation theorem can be stated as follows.
We say that a proper forcing notion $ {\mathcal P} $ is $
\omom $-bounding if
$$ \forall f \in \V^{{\mathcal P}} \cap \omom \ \exists g \in
\V \cap \omom \ \forall n \ \lft1(f(n) \leq g(n)\rgt1).$$
It is clear that $ {\mathcal P} $ is $ \omom $-bounding if
and only if $ {\mathcal P} $ preserves dominating families.

\begin{theorem}\label{boundCor}
  The countable support iteration of proper $\omom$-bounding
  forcing notions is $\omom$-bounding.
\end{theorem}

This is the ideal situation -- no matter what forcing notion we assign to
$\cov(\N)$, 
$\unif(\M)$, $\cof(\N)$ and $\unif(\N)$ it has to be $ \omom
$-bounding and this property is preserved under countable support
iteration iteration. 
\end{example}

\begin{example}
Covering numbers $\cov(\M)$ and $\cov(\N)$. The choice of forcing
notions that we assign to these invariants is determined by the
Theorem~\ref{4chars}, it has to be equivalent to Cohen and random real
forcing respectively.

The preservation theorem could be stated as follows (see
\cite{JuReIteration} or \cite{BJbook}).
\begin{theorem}
Suppose that $ {\mathcal P}_\delta = \lim_{\alpha<\delta} {\mathcal
  P}_\alpha $ ($ \delta $ limit) is a countable support iteration of
  proper forcing notions such that for every $ \alpha < \delta $, $
  {\mathcal P}_\alpha $ does not add random reals. Then $ {\mathcal
  P}_\delta $ does not add random reals.
\end{theorem}

The question whether this theorem remain true if we replace words
``random'' by ``Cohen'' is open. 
However, even if the preservation theorem for not adding Cohen reals
is true, both result cover only limit stages of the
iteration. For the successor steps we do not have an analog of Theorem
\ref{boundCor}, and indeed 
we can find two ccc forcing notions $ {\mathcal P} $
and $ {\mathcal Q} $ such that neither $ {\mathcal P} $ nor $
{\mathcal Q} $ adds random reals but $ {\mathcal P} \star {\mathcal Q}
$ adds random reals. Similarly for Cohen reals. 

These facts impose the following requirements:
\begin{itemize}
\item any iteration of finite length of forcing notions assigned to
${\mathfrak b}$, $ \unif(\M)$, ${\mathfrak d}$, $\cov(\M)$, $\unif(\N)$ and
$\cof(\N)$ does not add random reals,
\item iteration of any length of forcing notions assigned to
${\mathfrak b}$, $ \cov(\N)$, $\unif(\M)$, $\unif(\N)$ and
$\cof(\N)$ does not add Cohen reals.
\end{itemize}
It is easy to verify that each of the forcing notions chosen for these
invariants have
the required properties. However, the reasons why they, for example,
do not add Cohen reals are different in each case. Thus, the
preservation theorems are often difficult, technical and at the same
time not very general.  
\end{example}

The full proof that the construction outlined above is possible can be
found in \cite{BJbook}. A preservation theorem for not adding Cohen
reals that covers the cases we are interested in can be found in
\cite{RoSh470}. 

We will take all these constructions for granted 
and present some applications.

Let us consider the following examples:
 \begin{theorem}\label{model2}
  It is consistent with $\ZFCa$ that
$$ {\mathfrak b}=\alef_2 +
\cov(\N)=\unif(\N)=\alef_1 .$$
\end{theorem}
\begin{proof}
Recall that for any tree $T$,
$\stem(T)$ is an element of $T$ such that for all $t \in T$, $t
\subseteq \stem(T)$ or $\stem(T)\subseteq t$ and 
for $s \in T$, $\suc_T(s)=\{t: s \subseteq t \ \& \ |t|=|s|+1\}$.

 The Laver forcing ${\mathbf {LT}}$ is the following forcing notion:
$$T \in {\mathbf {LT}} \iff T \subseteq {}^{<\omega}\omega \text{ is a
  tree } \& \ 
\forall s \in T \ \lft1(|s| \geq \stem(T) \rightarrow
|\suc_T(s)|=\alef_{0}\rgt1).$$

For $T, T' \in {\mathbf {LT}}$, $T \geq T'$ if $T \subseteq T'$.

\begin{lemma}\label{laver}
\begin{enumerate}
  \item $\V^{\mathbf {LT}} \thinks \hbox{$\V \cap \omom$ is
      bounded in $\omom$}$.
  \item $\V^{\mathbf {LT}} \thinks \V
    \cap \twoom \not \in \N$.
  \item ${\mathbf {LT}}$ does not add random reals.
\end{enumerate}
Moreover (2) and (3) hold for the countable support iteration of Laver
forcing as well.
\end{lemma}
\begin{proof}
  See \cite{BJbook}.
\end{proof}
Let ${\mathcal P}_{\omega_2}$ be a countable support iteration of
length $ \alef_2 $ the
Laver 
forcing.
It follows from Lemma~\ref{laver} that ${\mathfrak  b}=\alef_2$ in 
$\V^{{\mathcal P}_{\omega_2}}$, while both $\cov(\N)$ and $\unif(\N)$
are equal to $\alef_1$.
\end{proof}

\begin{theorem}\label{model1}
It is consistent with $\ZFCa$ that
$$ {\mathfrak d}=\alef_1 +
\cov(\N)=\unif(\N)=\alef_2 .$$
  \end{theorem}
  \begin{proof}
We will use forcing notion $\mathbf {EE}$ defined below rather than
$\gS$, it has a much simpler definition and has the required
properties (the difficulties appear when unbounded reals are added).

The infinitely equal forcing notion $\mathbf {EE}$ 
 is defined as follows:
$p \in \mathbf {EE}$ if the following conditions are satisfied:
\begin{enumerate}
 \item $p:\dom(p) \longrightarrow \twolom$,
 \item $p(n) \in {^n 2}$ for all $n \in \dom(p)$, and
  \item $\dom(p) \subseteq \omega, \ |\omega \setmin
    \dom(p)|=\alef_0$.
\end{enumerate}
For $p,q \in \mathbf {EE}$  we define $p \geq q$ if $ p \supseteq q$.

\begin{lemma}\label{miller}
Forcing ${\mathbf {EE}}$ has the following properties:
\begin{enumerate}
  \item $\V^{\mathcal P} \thinks \V \cap \twoom \in \N$. In fact,
$$ \forall x \in \V \cap \twoom \ \exists^\infty n \ (x \rest n
=f_G(n)), $$
where $f_G$ is a generic real.
  \item ${\mathcal
      P}$ does not add random reals,
  \item ${\mathcal P}$ is $\omom$-bounding.
    \end{enumerate}
\end{lemma}
\begin{proof}
  See \cite{BJbook}.
\end{proof}
Let $\left\{{\mathcal P}_\alpha, \dot{{\mathcal Q}}_\alpha :
\alpha<\omega_2\right\}$ be 
a countable support iteration such that for every $\alpha < \omega_2$,
\begin{enumerate}
\item $\forces_\alpha \dot{{\mathcal Q}}_\alpha \simeq {\mathbf {EE}}$ if
$\alpha$ is even, and
\item $\forces_\alpha \dot{{\mathcal Q}}_\alpha \simeq {\mathbf B}$ if
$\alpha$ is odd.
\end{enumerate}

Let $G$ be a ${\mathcal P}_{\omega_2}$-generic filter over $\V \models
\CH$.

It is clear that $\V[G] \models
\unif(\N)=\cov(\N)=\alef_2$. To see that ${\mathfrak 
  d}=\alef_1$ in the extension note that both forcing
notions $\mathbf {B}$ and ${\mathbf {EE}}$ are
$\omom$-bounding and 
use \ref{boundCor}.
  \end{proof}

Now consider the corresponding problem concerning the families of 
small sets. 
The question is whether the models constructed for the Cicho\'n
diagram correspond to  the sets witnessing the strict inclusion
between the corresponding classes of sets. 

It is clear that we cannot show that in $\ZFCa$ alone. For example, it
is consistent that $\ADD(\N)=\ADD(\M)=\COV(\M)=[\reals]^{\leq
  \alef_0 }$ (a model for Borel Conjecture, see \cite{BJbook}). 

However, the theory $\ZFCa +\CH$ provides a sufficiently rich universe in which
$<$-results about invariants $\add$, $\cov$, etc in a natural way
yield $ \subsetneq$ results about $\ADD$, $\COV$, etc.

We will describe here several such constructions in detail. First
consider those that involve only forcing notions satisfying ccc.

\begin{theorem}\label{exam0}

  ($\ZFCa+\CH$) There is a set $X \subseteq \reals $ such that
$X \in \D$ and 
$X \not\in \UNIF(\N) \cup \UNIF(\N)$.
\end{theorem}
\begin{proof}
The construction is canonical. 
Set the cardinal invariants corresponding to the families that $X$
belongs to to $\mathbf\aleph_2$ and the other ones to
$\mathbf\aleph_1 $. 
In our case $\gd=\mathbf\aleph_2 $ and
$\unif(\N)=\unif(\M)=\mathbf\aleph_1 $. Now consider the forcing
notion that produces the  model for the dual setup,
i.e. $\gb=\mathbf\aleph_1 $ and
$\cov(\N)=\cov(\M)=\mathbf\aleph_2 $.
According to our table it is the iteration of Cohen and random
forcings, ${ \mathbf C} \star {\mathbf B} $.
Let $\{M_\alpha : \alpha < \mathbf\aleph_1\}$ be an increasing
sequence of contable submodels of $\HH(\lambda)$ such that 
\begin{enumerate}
\item $\twoom \subseteq \bigcup_{\alpha<\omega_1} M_\alpha $.
\item for ever $ \alpha < \omega_1$, $M_{\alpha+1} \thinks M_\alpha $
  is countable,
\item $\{M_\beta: \beta\leq \alpha\} \in M_{\alpha+1}$.
\end{enumerate}

For each $\alpha $ choose a pair $(c_\alpha, r_\alpha) \in
M_{\alpha+1}$ such that $(c_\alpha, r_\alpha)$ is ${\mathbf C}\star
{\mathbf B} $ over $M_\alpha $. Note that such a pair will also be
generic over $M_\beta $ for $ \beta < \alpha $.
Let $z_\alpha $ encode $(c_\alpha, r_\alpha )$ as 
$$z_\alpha(n)=\left\{
  \begin{array}{ll}
c_{\alpha}(k) & \text{if }n=2k\\
r_\alpha(k) & \text{if }n=2k+1
  \end{array}\right. .$$
Let $X=\{z_\alpha: \alpha <\omega_1\}$.
We will show that $X$ has the required properties.

To show that $X \in \D$ fix a Borel function $F: \reals
\longrightarrow \omom$ and find $\alpha_0$ such that $F$ is coded in
$M_{\alpha_0}$.
Let $f$ be any function which dominates $M_{\alpha_0} \cap \omom$.
For any $\alpha < \omega_1$, $F(z_\alpha)\in M_{\alpha_0}^{{\mathbf
  C}\star {\mathbf B}}$. Since ${\mathbf C}\star {\mathbf B} $ does
not add dominating reals it follows that for every $\alpha $ there is
a function $g \in M_{\alpha_0}\cap \omom$ such that $g \not\leq^\star
F(z_\alpha)$. Since $g $ is dominated by $f$ we conclude that
$f \not \leq^\star F(z_\alpha)$ for every $\alpha < \omega_1$.

To see that $X \not \in \UNIF(\M)\cup \UNIF(\N)$ let $Y=\{c_\alpha:
\alpha<\omega_1\}$. Observe that $Y$ is a continuous image of $X$. 
Moreover, if $F \in M_{\alpha_0}$ is a meager set then $c_\alpha
\not\in F$ for $\alpha>\alpha_0$ since $c_\alpha $  is a Cohen real
over $M_{\alpha_0}$.
The argument that $X \not\in \UNIF(\N)$ is identical.
\end{proof}

Observe that the crucial point of the above construction is that the
real $z_\alpha $ efined at the step $ \alpha $ is generic not only
over model $M_\alpha $ but also over models $M_\beta $ for
$\beta<\alpha $.
To illustrate this point  suppose that ${\mathcal P} $ is a forcing
notion, $M \subseteq N$ are two submodels of $\HH(\lambda)$ and
${\mathcal P} \in M$.
Let ${\mathcal A} \in M$ be a maximal antichain in ${\mathcal P} $. If 
${\mathcal P} $ satisfies ccc then $ {\mathcal A} \subseteq M$, as a
range of a function on $\omega $. If $
{\mathcal P} $ is absolutely ccc then $N \thinks {\mathcal A} $ is an
maximal antichain, so a ${\mathcal P} $-generic real over $N$ is also
${\mathcal P} $-generic over $M$.
If $ {\mathcal P} $ is not absolutely ccc then we no longer know if
${\mathcal A} $ is a maximal antichain in $N$. In fact, we do not know 
if $ {\mathcal A} $ is an antichain at all. However, if both $M$ and
$N$ are elementary submodels of $\HH(\lambda)$, then $N \thinks {\mathcal A} $ is a
maximal antichain. Finally, if ${\mathcal P} $ does not satisfy ccc,
then it is no longer true that $ {\mathcal A} \subseteq M$, so a
${\mathcal P} $-generic real over $\V$ may not be generic over $M$.
Recall that a condition $p \in {\mathcal P} $ is $(M,{\mathcal
  P})$-generic if $p$ forces that the above situation does not
happen. If for every countable $M \prec \HH(\lambda)$ the collection of
$(M, {\mathcal P})$-generic conditions is dense in $ {\mathcal P} $,
then ${\mathcal P} $ is proper.

The following strengthening of properness will allow us to carry out
the construction from the proof of \ref{exam0} for non-ccc posets.
\begin{definition}
  Suppose that ${\mathcal P} $ is a forcing notion and
  $\alpha<\omega_1$ is an ordinal. We say that ${\mathcal P} $ is
  {\it $\alpha $-proper} if for every sequence $\{M_\beta : \beta\leq
  \alpha\}$ such that 
  \begin{enumerate}
  \item for every $\beta $, $M_\beta $ is a countable elementary
    submodel of $\HH(\lambda)$,
  \item $\{M_\gamma:\gamma\leq \beta \}\in M_{\beta+1}$,
  \item $M_{\beta+1} \thinks M_\beta $ is countable,
  \item $M_\lambda = \bigcup{\beta<\lambda} M_\beta $ for limit
    $\lambda $,
  \item ${\mathcal P} \in M_0$,
  \end{enumerate}
and for every $p \in {\mathcal P} \cap M_0$, there exists $q \geq p$
which is $(M_\beta, {\mathcal P} )$-generic for $\beta \leq \alpha $. 
\end{definition}

%Note that if ${\mathcal P} $ is $ \alpha  $-proper for each $\alpha <
%\omega_1$ then we can carry out the construction as in the proof of
%\ref{exam0} with $ {\mathcal P} $ instead of ${\mathbf C}\star
%{\mathbf B}$.

\begin{definition}\label{axiomA}
A forcing notion ${\mathcal P}$ satisfies {\it axiom ${\mathsf A}$} if there
exists a sequence $\left\{\leq_n : n \in \omega\right\}$ of orderings
on ${\mathcal P}$ (not necessarily transitive) such that  
\begin{enumerate}
\item if $p \geq_{n+1} q$, then $p \geq_n q$ and $p \geq q$ for $p,q
  \in {\mathcal P}$, 
\item if $\<p_n:n \in \omega\>$ is a sequence of conditions such that
  $p_{n+1} \geq_{n} p_n$,  then there exists
$p \in {\mathcal P}$ such that $p \geq_n p_n$ for all $n$, and 
\item if ${\mathcal A} \subseteq {\mathcal P}$ is an antichain, then for every
$p \in {\mathcal P}$ and $n \in \omega$ there exists $q \geq_n p$ such
that $\left\{r \in {\mathcal A} : q \hbox{ is compatible with } r\right\}$ is countable.
\end{enumerate}
\end{definition}

All forcing notions assigned to the cardinal invariants from Cicho\'n
diagram  satisfy axiom ${\mathsf A}$.

\begin{lemma}
  If $ {\mathcal P} $ satisfies Axiom ${\mathsf A}$ then ${\mathcal P} $ is
  $\alpha $-proper for every $ \alpha< \omega_1$.
\end{lemma}
\begin{proof}
  Induction on $ \alpha $. Let $\<M_\beta: \beta \leq \alpha\>$ be a
  sequence of models having the required properties. Fix  $p \in
  {\mathcal P} \cap M_0$ and $n \in \omega $.
We will find $q \geq_n p$ which is $M_\beta $-generic for $\beta \leq
\alpha $.
If $\alpha =0$ then it is the usual proof that Axiom ${\mathsf A}$ implies
properness.
If $ \alpha = \gamma+1$ then first find $q' \geq_n p$ which is
$M_\delta $-generic for $\delta \leq \gamma $ and then use properness
of ${\mathcal P} $ to get $q \geq_n q'$ which is $M_\alpha $-generic.
If $\alpha $ is limit then fix an increasing sequence $\<\alpha_n: n
\in \omega\>$ such that $\sup_n \alpha_n=\alpha $.
Use the induction hypothesis to find conditions $\{p_k: k\in \omega \}$ 
such that 
\begin{enumerate}
\item $p_{k+1} \in M_{\alpha_k+1}$,
\item $p_k$ is $M_\gamma $-generic for $\gamma < \alpha_k$,
\item $p_{k+1} \geq_{n+k} p_k$ for each $k$.
\end{enumerate}
Let $q $ be such that $q \geq_{n+k} p_k$ for each $k$, it is the
condition we are looking for.
\end{proof}

\begin{theorem}\label{exam1}
  ($\ZFCa+\CH$) There is a set $X \subseteq \reals $ such that
$X \in \B$ and 
$X \not\in \COV(\N) \cup \UNIF(\N)$.
\end{theorem}
\begin{proof}
In terms of
cardinal invariants the statement of the theorem corresponds to the
dual to  the model for  
${\mathfrak b}=\alef_2 $ and
$\cov(\N)=\unif(\N)=\alef_1$, 
that is the one where ${\mathfrak
  d}=\alef_1 $ and $\cov(\N)=\unif(\N)=\alef_2$. 
The set we are looking for is defined using the forcing notion used to
construct that model (cf. Theorem~\ref{model1}).
 
Let $\{f_\alpha: \alpha<\omega\}$ be an enumeration of $\reals$.
  Let $\<M_\alpha : \alpha<\omega_1\>$ be a sequence of countable
  elementary submodels of $\HH(\lambda)$ such that 
  \begin{enumerate}
  \item $f_\alpha \in M_\alpha $,
  \item $\<M_\beta: \beta \leq \alpha\> \in M_{\alpha+1}$, and
    $M_{\alpha+1} \thinks M_\alpha $ is countable,
  \item $M_\gamma =\bigcup_{\alpha<\gamma} M_\alpha $ for limit $
    \gamma $.
  \end{enumerate}
Note that from (2) it follows that for every $ \beta < \alpha $,
$M_\alpha \thinks $ ``$M_\beta $ is countable.''
Let $\<e_\alpha, r_\alpha : \alpha<\omega_1\>$ be a sequence of reals
such that 
\begin{enumerate}
\item $e_\alpha, r_\alpha \in M_{\alpha+1}$,
\item $e_\alpha $ is ${\mathbf {EE}}$-generic over $M_\beta $ for $
  \beta \leq \alpha $,
\item $r_\alpha $ is ${\mathbf {B}}$-generic over $M_\beta[e_\alpha] $ for $
  \beta \leq \alpha $.
\end{enumerate}
For $ \alpha < \omega_1$ define
$$z_\alpha(n)=\left\{
  \begin{array}{ll}
e_\alpha(k) & \hbox{if } n=2k\\
r_\alpha(k) & \hbox{if } n=2k+1
  \end{array}\right. .$$
Let $Z=\{z_\alpha : \alpha<\omega_1\}$.

\bigskip

$Z \not\in \UNIF(\N)$. The set $X=\{r_\alpha: \alpha \in \omega_1\}$ is
a Borel image of $Z$. Given $f \in \omom $ find $ \alpha $
such that $f=f_\alpha $. Notice that $r_\beta \not \in (N)_f$ for $
\beta>\alpha $. In particular, no uncountable subset of $Z$ is in
$\UNIF(\N)$.

\bigskip

$Z \not\in \COV(\N)$. Consider the set $Y=\{e_\alpha: \alpha <
\omega_1\}$ which is a Borel image of $X$.
Let $\overline{P}=\left\{f \in {}^\omega \left([\omega]^{<\omega}\right) :
\forall n \ \lft1(f(n) \in {^n 2}\rgt1)\right\}$.
Let 
$$\widetilde{H}=\{(f,x): f \in \overline{P}, x \in \twoom \ \&\
\exists^\infty n \ x\rest n = f(n)\}.$$
It is easy to see that $\widetilde{H}$ is a Borel set in
$\overline{P}\times \twoom $ and $(\widetilde{H})_f \in \N$ for
every $f$.
Suppose that $x \in \twoom $. Find $ \alpha $ such that $x \in
M_\alpha $ and note that for $ \beta>\alpha $, $x \in
(\widetilde{H})_{e_\beta}$. It follows that no uncountable subset of $Z$ is in
$\COV(\N)$.
 
\bigskip

$X \in \D$ 
Let $F: X \longrightarrow \omom  $ be a Borel
mapping. Find $ \alpha $ such that $F$ is coded in $M_\alpha $.
 Let $f \in \omom $ be such that for every $g \in M_\alpha
 \cap \omom$, $ g \leq^\star f$. Since $M_\alpha $ is
 countable, such an $f$ exists.
Since both $ \mathbf B $ and $ \mathbf {EE}$ are $ \omom
$-bounding (so is $ \mathbf B \star {\mathbf {EE}}$) for every  $
\beta > \alpha $, there exists $g \in M_\alpha 
$ such that $F(z_\alpha) \leq^\star g \leq^\star f$.
\end{proof}
\begin{theorem}\label{exam3}
  ($\ZFCa+\CH$) There is a set $X \subseteq \reals $ such that 
$X \in \COV(\N) \cap \UNIF(\N)$ and $X \not \in \D$.
\end{theorem}
\begin{proof}
    Let $\<M_\alpha : \alpha<\omega_1\>$ be a sequence of countable
  elementary submodels of $\HH(\lambda)$ as in the previous proof.

In this case we use the Laver forcing from Theorem~\ref{model2}.
The only difference is that in order to ensure that the constructed
set belongs to $\COV(\N)$ we construct a set of witnesses for that.

Let $\<l_\alpha, r_\alpha : \alpha<\omega_1\>$ be a sequence of reals
such that 
\begin{enumerate}
\item $l_\alpha, r_\alpha \in M_{\alpha+1}$,
\item $l_\alpha $ is ${\mathbf {LT}}$-generic over $M_\beta $ for $
  \beta \leq \alpha $,
\item $r_\alpha $ is ${\mathbf {B}}$-generic over $M_\alpha[l_\beta] $
  for all $   \beta < \omega_1 $.
\end{enumerate}
To meet the condition (3) we need the following result:
\begin{theorem}\label{mainth}
Suppose that $N \prec \HH(\lambda)$ is a countable model of $\ZFCa$.
Let $S \in N \cap {\mathbf {LT}}$ and  let $x$ be a  random real 
over $N$.
 There exists $T
 \geq
S$ such that  $T$ is $N$-generic and 
$T \forces_{{\mathbf {LT}}} x
\text{ is random over } N[\dot{G}].$ 
\end{theorem} 
\begin{proof}
  See \cite{JudShel90Kun}, \cite{PawLav} or \cite{BJbook}.
\end{proof}

Let $X=\{l_\alpha : \alpha<\omega_1\}$.
The difference between this and the previous construction is that
we define the set of witnesses $\{r_\alpha : \alpha <\omega_1\}$ that
$X \in \COV(\N)$.

\bigskip

$X \in \COV(\N)$. Let $H \subseteq \omom \times \twoom $ be
a Borel set with null sections.
Find $ \alpha $ such that $H \in M_\alpha $. 
Note that
$$r_\alpha \not \in \bigcup_{\beta<\omega_1} (H)_{l_\beta},$$
since $r_\alpha $ is random over $M_\alpha[l_\beta]$ for all $
\beta $ and $(H)_{l_\beta} \in M_\alpha[l_\beta]$.

\bigskip

$X \in \UNIF(\N)$. Let $F: X \longrightarrow \twoom $ be a Borel
mapping. Find $ \alpha $ such that $F$ is coded in $M_\alpha $.
 Let $B=\bigcup \{A: A \in \N \cap M_\alpha \}$. Since $M_\alpha $ is
 countable, $B$ is a null set.
By Lemma~\ref{laver}(3) for every $ \beta > \alpha $, $F(l_\alpha) \in
B$.  

\bigskip

$X \not\in \D$. This is obvious, by Lemma~\ref{laver}(1), for every $
\alpha $ 
$$ \forall f \in M_\alpha \cap \omom \ \forall^\infty n \ \lft1(f(n)
< l_\alpha(n)\rgt1).$$
\end{proof}

The method of constructing  counterexamples to the Cicho\'n diagram 
described above is very elegant
and effective but assumes a rather large body of knowledge involving
forcing, preservation theorems etc.  
We will conclude this section with a sketch of an alternative method of
constructing examples of small sets which is also quite general but more
direct. Along the way translate the forcing 
results that we have used into statements about sets of reals.

Suppose that ${\mathcal P} $ is a forcing notion and conditions of
${\mathcal P}$ are sets of reals. Note that all forcing notions
associated with The Cichon Diagram are of this form. For a description 
of much larger class of forcing of that kind see \cite{RoSh470}.

Let 
$$I_{{\mathcal P}}=\left\{X \subseteq \reals: \forall p \in {\mathcal P} \
\exists q \geq p \ \lft1(q \cap X=\emptyset\rgt1)\right\}.$$

The following lemma lists the obvious observations about $I_{\mathcal
  P} $. 
\begin{lemma}
  \begin{enumerate}
  \item $I_{{\mathcal P}}$ is an ideal,
  \item $X \in I_{{\mathcal P}}$ iff there exists a maximal antichain
    ${\mathcal A} \subseteq {\mathcal P} $ such that $X \cap \bigcup
    {\mathcal A} =\emptyset$,
  \item $\forall p \in {\mathcal P} \ (p \not\in I_{{\mathcal P}})$.
  \end{enumerate}
\end{lemma}

Suppose that ${\mathcal P} $ is a forcing notion satisfying Axiom ${\mathsf A}$.
Let
$$I^\omega_{{\mathcal P}}=\left\{X \subseteq \reals: \forall p \in
{\mathcal P} \ \forall n \in \omega \
\exists q \geq_n p \ \lft1(q \cap X=\emptyset\rgt1)\right\}.$$

Note that $I^\omega_{{\mathcal P}}$ is a $ \sigma $-ideal contained in 
$I_{{\mathcal P}}$.

If $ {\mathcal P} $ satisfies ccc, then we can witness that $ {\mathcal
  P} $ satisfies Axiom ${\mathsf A}$ by putting $p \leq_0 q$ if $p \leq q$, and
for $n>0$, $p \leq_n q$ if $p=q$.
In this case $I^\omega_{{\mathcal P}}=\emptyset$.
However, for non-ccc forcings as well as some ccc posets (like random
real algebra ${\mathbf B}$) we can define $\leq_n$'s in such a way that $I^\omega_{{\mathcal
    P}}=I_{{\mathcal P}}$.

First we will describe how to translate the forcing theorems.

\begin{lemma}
Suppose that $ {\mathcal P} $ 
is a forcing notion such that 
\begin{enumerate}
\item $ {\mathcal P} $ is proper,
\item for every $\V$-generic filter $G \subseteq {\mathcal P} $ there
  exists a real $x_G$ such that $\V[G]=\V[x_G]$,
\item conditions of $ {\mathcal P} $ are Borel sets of reals, ordered
  by inclusion,
\item every countable antichain in ${\mathcal P} $ can be represented
  by a countable family of pairwise disjoint elements of ${\mathcal P} 
  $.
\end{enumerate}
Then for every $ {\mathcal P} $-name $\dot{x}$ such that
$\forces_{\mathcal P} \dot{x} \in \twoom$ and $p \in {\mathcal P} $ there exists a Borel
function $F \in \V$, $F: \twoom \longrightarrow \twoom$ and a $q \geq
p$ such that $q \forces_{\mathcal P} \dot{x}=F(x_{\dot{G}})$.
\end{lemma}
\begin{proof}
  Fix $\dot{x}$ and let ${\mathcal A}_n$ be a maximal antichain of
  conditions deciding $\dot{x} \rest n$. Use properness to find $q
  \geq p$ such that each $ {\mathcal A}_n'=\{r\in {\mathcal A}_n: r
  \text{ is compatible with }q\}$ is countable. By 
  the assumption we can assume that elements of ${\mathcal A}_n'$ are
  pairwise disjoint.
Define $F_n:q \longrightarrow 2^n$ as
$$F_n(x)=s \text{ if } x \in r \in {\mathcal A}_n' \text{ and } r
\forces_{\mathcal P} \dot{x} \rest n=s.$$
Note that $F=\lim_n F_n$ is the function we are looking for.
\end{proof}

Let $ {\mathcal P} $ be a forcing notion satisfying the assumptions of 
the above lemma.

\begin{itemize}
\item ${\mathcal P} $ does not add random
  reals if for every ${\mathcal P} $-name $ \dot{x}$ for an element of 
  $\twoom$ and every $p \in {\mathcal P} $ there is $q \geq p$ and $H
  \in \V \cap \N$ such that $q \forces_{\mathcal P} \dot{x} \in H$.
\item ${\mathcal P} $ is $\omom$-bounding if for every ${\mathcal P} $-name $ \dot{f}$ for an element of 
  $\omom$ and every $p \in {\mathcal P} $ there is $q \geq p$ and $g
  \in \V \cap \omom$ such that $q \forces_{\mathcal P} \dot{f}
  \leq^\star g$.
\item ${\mathcal P} $ preserves outer measure if for every set of
  positive outer measure $X
  \subseteq \twoom$, $X \in \V$  and every $\dot{F}$, a  ${\mathcal P} 
  $-name for a Borel function from $\twoom$ to $\omom$ and $p \in
  {\mathcal P} $ there is $q \geq p$  such that $q \forces_{\mathcal
    P} X \setmin (N)_{\dot{F}(x_{\dot{G}})})\neq \emptyset$.
\end{itemize}

These statements translate as:
\begin{itemize}
\item (not adding random reals) For every Borel fuction $F:
  \twoom\longrightarrow\twoom$ and $p\in {\mathcal P} $ there is a set 
  $H \in \N$, $q \geq p$ and $A \in I_{{\mathcal P}}$ such that
  $F{``}(q\setmin A) \subseteq H$.
\item (${\mathcal P}$ is $\omom$-bounding) For every Borel fuction $F:
  \twoom\longrightarrow\omom$ and $p\in {\mathcal P} $ there is a function 
  $f \in \omom$, $q \geq p$ and $A \in I_{{\mathcal P}}$ such that
  $F{``}(q\setmin A) \leq^\star f$.
\item (${\mathcal P} $ preserves outer measure) for every set of
  positive outer measure $X
  \subseteq \twoom$, and every Borel function $F:
  \twoom\longrightarrow \omom$ and $p \in
  {\mathcal P} $ there is $q \geq p$  and $A \in I_{{\mathcal P}}$
  such that $X \setmin \bigcup_{x \in q \setmin A} (N)_{F(x)}\neq \emptyset$.
\end{itemize}

If in addition ${\mathcal P} $ satisfies axiom ${\mathsf A}$ and
$I_{{\mathcal P}}=I^\omega_{{\mathcal P}}$, then we can put $A=\emptyset$.

\bigskip

\begin{proof}[Second proof of \ref{exam1}]
  For $p,q \in \mathbf {EE}$ and $n \in \omega$ we define
$p \geq_n q$ if $ p \geq q$ and first $n$ elements of $\omega
  \setmin \dom(p)$ and $\omega
  \setmin \dom(q)$ are the same.

For $p,q \in \mathbf B $ and $ n \in \omega $ let $p \geq_n q$ if $p
\geq q$ and $\mu(q \setmin p) \leq 2^{-n} \cdot \mu(q)$.

The forcing notions $\mathbf {EE}$, $ \mathbf B $ (and the remaining ones
as well) can be represented as the collections of perfect subsets of
$\twoom $ (or $ \omom $). This is not critical for the construction,
but it makes it more natural.

In case of $\mathbf {EE}$ for $ n \in \omega $ let $k_n=2^{n+1}-1$.
Consider sets $P \subseteq \twoom $ of form $\bigcap_{n \in \omega}
[C_n]$, where $\{C_n: n \in \omega\}$ satisfies the following
conditions:
\begin{enumerate}
\item $C_n \subseteq {^{[k_n,k_{n+1})}}2$,
  \item for every $n$, $|C_n|=1$ or $|C_n|=2^n$ (so
    $C_n={^{[k_n,k_{n+1})}}2$),
  \item $ \exists^\infty n \ |C_n|=2^n$.
\end{enumerate}
  It is clear that every condition $p \in \mathbf {EE}$ corresponds to
  a set $P$ as above and vice versa.
Therefore from now on we identify $\mathbf {EE}$ with these sets.

Let $ {\mathbf B} \star \mathbf {EE}$ be the collection of 
subsets $H \subseteq 2^\omega \times 2^\omega $ such that
\begin{enumerate}
\item $H$ is Borel and $\dom(H)=\{x : (H)_x \neq \emptyset\} \in
  \mathbf B $,
\item $\forall x \ \left( (H)_x\neq \emptyset \rightarrow  (H)_x \in \mathbf {EE}\right)$.
\end{enumerate}
The
elements of $ {\mathbf B} \star \mathbf {EE}$ are $ {\mathbf B}
$-names for the elements of $\mathbf {EE}$.
Thus, the set $ {\mathbf B} \star \mathbf {EE}$ indeed corresponds 
to the iteration of ${\mathbf B} $ and $\mathbf {EE} $. 
For $H_1, H_2 \in {\mathbf B} \star \mathbf {EE}$ and $ n \in \omega $ 
let $H_1 \geq H_2$ mean that
\begin{enumerate}
\item $\dom(H_1) \geq_n \dom(H_2)$,
\item $ \forall x \in \dom(H_1) \ \lft1((H_1)_x \geq_n (H_2)_x\rgt1).$
\end{enumerate}
Note that $\geq_n$ on ${\mathbf B} \star \mathbf {EE}$ witnesses that
it satisfies Axiom ${\mathsf A}$.

Let $\<x_\alpha: \alpha<\omega_1\>$ be an enumeration
of $\twoom $, and $\<F_\alpha: \alpha<\omega_1\>$ of
$\borel(\twoom\times \twoom,\omom)$, and $\<f_\alpha: \alpha
<\omega_1\>$ an enumeration of $\omom$.
We will build an $\omega_1$-tree ${\mathbf A}$ 
of elements of $ {\mathbf B} \star \mathbf {EE}$.
Let ${{\mathbf A}}_\alpha $ denote the $ \alpha $-th level of ${{\mathbf A}}$.
The tree ${{\mathbf A}}$ satisfies the following inductive conditions:
\begin{enumerate}
\item $ \forall \beta > \alpha \ \forall n \ \forall H \in {{\mathbf A}}_\alpha \
  \exists H' \in {{\mathbf A}}_\beta \ \lft1(H' \geq_n H\rgt1).$
\item $ \exists f \in \omom \ \forall H \in {{\mathbf A}}_{\alpha+1}\
  \lft1({F_\alpha}{``}(H) \leq^\star f\rgt1)$,
\item $\forall H \in {{\mathbf A}}_{\alpha+1}\
  \lft1(\dom(H) \cap (N)_{f_\alpha} = \emptyset\rgt1)$,
\item $\forall H \in {{\mathbf A}}_{\alpha+1}\
  \forall x \in \dom(H) \ \exists^\infty n \
  \lft2( |C_n^x|=1 \ \& \ C_n^x \subseteq x_\alpha \rgt2)$, where  
$(H)_x = \bigcap_n [C_n^x]$.

\end{enumerate}

{\sc  Case 1}. \hspace{0.1in} $ \alpha =\beta+1$. We will describe how to build a set
of immediate successors of an element $H \in {{\mathbf A}}_\beta $. 
Given $H \in {{\mathbf A}}_\beta $ and $ n \in \omega $  
find $H_n' \geq_n H$ satisfying conditions (3) and (4).
By further shinking we can ensure that (2) holds as well. Condition
(2) is just the statement that the iteration of $ {\mathbf {EE}}$ and
  ${\mathbf B} $ is $\omom$-bounding.

\bigskip

{\sc Case 2}. \hspace{0.1in}
$ \alpha $ is limit.
Suppose that $H \in {{\mathbf A}}_{\beta_0}$ for  some $ \beta_0 \in \omega $ and
that $n \in \omega $ is given.
Fix an increasing 
 sequence $\<\beta_k: k \in \omega\>$ such that $ \beta_k
\rightarrow \alpha $.
Choose a sequence $\<H_k: k\in \omega \>$ such that 
\begin{enumerate}
\item $H_0=H_1=\dots =H_n=H$,
\item for $ k \geq 0$, $H_{n+k+1} \geq_n H_{n+k}$,
\item $H_{k+n} \in {{\mathbf A}}_{\beta_k}$.
\end{enumerate}
Use Axiom ${\mathsf A}$ to find $H'$ such that $H' \geq_k H_k$. Level
${\mathbf A}_\alpha $ will consist of elements selected in this way.

Let $X=\{(x_\alpha,y_\alpha): \alpha < \omega_1\}$ be a selector from elements of
${{\mathbf A}}$. 
Note that $\pi_1(X)=\{x_\alpha: \alpha < \omega_1\} \not\in \UNIF(\N)$
(by (3)), 
$\pi_2(X)=\{y_\alpha: \alpha < \omega_1\} \not\in \COV(\N)$ (by (4))
and $X \in \D$ (by (2)).
\end{proof}

\bigskip

Now let us look at the set constructed in Theorem \ref{exam3}.

\bigskip

\begin{proof}[Second proof of \ref{exam3}]
For every $T \in {\mathbf {LT}}$ and $s \in {}^{<\omega}\omega$ define a node
$T(s)$ in the following way:
$T(\emptyset)=\stem(T)$ and for $n \in \omega$ 
let $T(s^\frown n)$ be the $n$-th 
element of $\suc_T\lft1(T(s)\rgt1)$.

For $T,T' \in \mathbf{LT}$  and $n \in \omega$ define
$T \geq_n T'$ if  $T \geq T' \ \&\ 
\forall s \in n^{\leq n}\ 
\lft1(T(s)=T'(s)\rgt1).$
In particular, $T \geq_0 T'$ is equivalent to $T \geq T'$ and
$\stem(T)=\stem(T')$. 
It is easy to check that Laver forcing satisfies Axiom ${\mathsf A}$.

Suppose that 
\begin{enumerate}
\item $\<f_\alpha: \alpha<\omega_1\>$ is an enumeration of 
$ \omom $, 
\item $\<F_\alpha: \alpha<\omega_1\>$ is an enumeration of 
$\borel(\omom,\twoom)$,
\item $\<G_\alpha: \alpha<\omega_1\>$ is an enumeration of 
$\borel(\omom,\omom)$.
\end{enumerate}
We build an $\omega_1$-tree ${\mathbf A}$ satisfying the following
inductive 
conditions:
\begin{enumerate}
\item $ \forall \beta > \alpha \ \forall n \ \forall T \in {{\mathbf A}}_\alpha \
  \exists S \in {{\mathbf A}}_\beta \ \lft1(S \geq_n T\rgt1),$
\item $ \forall T \in {{\mathbf A}}_{\alpha+1} \ \forall x \in [T] \
  f_\alpha \leq^\star x$, (${\mathbf {LT}}$ adds a dominating real (\ref{laver}(1)))
\item for every $T \in {\mathbf A}_{\alpha+1}$, $F_\alpha{``}(T)$
  has measure zero, (${\mathbf {LT}}$ does not add random reals (\ref{laver}(3)))
\item $ \forall T \in {\mathbf A}_{\alpha+1} \ \lft2({}^\omega 2 \setmin  \bigcup_{x
    \in [T]} (N)_{G_\alpha(x)}$ is uncountable $\rgt2)$. (${\mathbf {LT}}$
  preserves outer measure (\ref{laver}(2))),

\end{enumerate}

%The proof that a selector $X$ from the elements of ${\mathbf A}$ has the
%required properties is similar to the previous proof.
%Condition (2) of the definition of ${\mathbf A}$ yields that $X
%\not\in \D$ while conditions (3) and (4) guarantee that $X \in
%\COV(\N) \cap \UNIF(\N)$.
Next we want to chose a selector $X$ from elements of $\mathbf
A$. Condition (2) will guarantee that $X \not \in \D$ and (3) that
$X\in \UNIF(\N)$.
Unfortunately (4) does not suffice to show that $X \in \COV(\N)$. It
is conceivable that $2^\omega = \bigcup_{T \in {\mathbf A}_{\alpha+1}} 
\bigcup_{x    \in [T]} (N)_{G_\alpha(x)}$, because $\COV(\N)$ is not a 
$ \sigma $-ideal.
Therefore we need stronger property:

\begin{enumerate}
\item[4'.] For every Borel function $F: \omom  \longrightarrow
  \omom $ and a sequence $\{T_n: n \in \omega\}$ of conditions 
  in $\mathbf {LT}$ there exists an uncountable set $Y \subseteq
  2^\omega $ such that for each $x \in Y$we can find sequence
  $\{S^n_k: n,k \in \omega \}$ such that $S^n_k \geq_k T_n$ and 
$y \not\in \bigcup_{n,k} \bigcup_{x \in [S^n_k]} (N)_{F(n)}$.
\end{enumerate}
Property (4') is a translation of \ref{mainth}.

Now we construct $X$ along with ${\mathbf A}$. On the step $\alpha $
we have $X_\alpha $ and ${\mathbf A}_\alpha $. Let ${\mathbf
  A}_\alpha =\{T_n: n \in \omega\}$ and pick  $y
\not\in X_\alpha $ together with  $\{S^n_k: n,k \in \omega \}={\mathbf 
  A}_{\alpha+1}$ as in
(4').
 
 \end{proof}

\bigskip

{\bf Historical remarks}
Parts (1) and (4) of Theorem~\ref{4chars} are due to Truss
\cite{Tru77Set} and \cite{Tru88Con} and parts (2) and (3) to Solovay
\cite{Sol70Mod}. 
Theorem~\ref{boundCor} and other preservation results are due to
Shelah \cite{she:pforcing}. Various presentations of these results
appear in \cite{GT}, \cite{JuReIteration} and most generally in
\cite{shelahbook}. 
Models for the Cicho\'n diagram were constructed by Miller in
\cite{Mil81Som}, more in \cite{JudShel90Kun} and the latest ones in
\cite{BarJuShelCi}. Theorem \ref{laver}(2) is due to Judah and Shelah
(\cite{JudShel90Kun}), the remaining parts are due to Laver
\cite{Lav76Con}.  The forcing $\mathbf {EE}$ and Lemma~\ref{miller}
are due to Miller \cite{Mil81Som}. 
Brendle \cite{Bre95Gen}   constructed the counterexamples for the $
\subsetneq$ for the 
families of small sets.  
This type of constructions were already considered in \cite{FreMil88Som}.
The technique of ``Aronszajn tree of perfect
sets'' was invented by Todorcevic (see \cite{GalMil84Gam}).   
Theorem~\ref{mainth} is due to Judah and Shelah.

\bigskip

{\bf Acknowledgement}: I would like to thank Andrzej Ros{\l}anowski
for careful reading early versions of this article and many helpful
suggestions and remarks.

\end{document}